\documentclass[11pt]{article}    
\usepackage{epsfig,amsmath,amssymb,amsthm,multirow,float}
\usepackage[section]{placeins}    
\newtheorem{remark}{Remark}[section]

\begin{document}

\title{Unconditionally stable time splitting methods for the electrostatic analysis of solvated biomolecules}

\author{Leighton Wilson and Shan Zhao\footnote{
Corresponding author. Tel: 1-205-3485303, Fax: 1-205-3487067,
Email: szhao@ua.edu}\\
Department of Mathematics, University of Alabama, Tuscaloosa, AL 35487, USA.
}

\date{\today} 
\maketitle

\begin{abstract}

This work introduces novel unconditionally stable operator splitting methods for solving the time dependent nonlinear Poisson-Boltzmann (NPB) equation for the electrostatic analysis of solvated biomolecules. In a pseudo-transient continuation solution of the NPB equation, a long time integration is needed to reach the steady state. This calls for time stepping schemes that are stable and accurate for large time increments. The existing alternating direction implicit (ADI) methods for the NPB equation are known to be conditionally stable, although being fully implicit. To overcome this difficulty, we propose several new operator splitting schemes, in both multiplicative and additive styles, including locally one-dimensional (LOD) schemes and additive operator splitting (AOS) schemes. The nonlinear term is integrated analytically in these schemes, while standard discretizations with finite differences in space and implicit time integrations are used. The proposed schemes become much more stable than the ADI methods, and some of them are indeed unconditionally stable in dealing with solvated proteins with source singularities and non-smooth solutions. Numerically, the orders of convergence in both space and time are found to be one. Nevertheless, the precision in calculating the electrostatic free energy is low, unless a small time increment is used. Further accuracy improvements are thus considered, through constructing a Richardson extrapolation procedure and a tailored recovery scheme that replaces the fast Fourier transform method by the operator splitting method in the vacuum case. After acceleration, the optimized LOD method can produce a reliable energy estimate by integrating for a small and fixed number of time steps. Since one only needs to solve a tridiagonal linear system in each independent one dimensional process, the overall computation is very efficient. The unconditionally stable LOD method scales linearly with respect to the number of atoms in the protein studies, and is over 20 times faster than the conditionally stable ADI methods.



\end{abstract}

\section{Introduction}
Analysis of the underlying biomolecular solvation is critical when carrying out quantitative descriptions of various important biological processes at the atomic level, such as protein folding and protein ligand bonding, DNA recognition, transcription, and translation. From a biological perspective, solvation analysis is concerned with interactions between a solute macromolecule and surrounding solvent ions. From a mathematical perspective, these solute-solvent interactions may be represented via solvation energies with contributions from polar and nonpolar sources. The polar portion arises from electrostatic interactions, which may be represented with the Poisson-Boltzmann (PB) model \cite{Davis90,Honig95}. The PB model provides a framework by which to model the distribution of electrostatic potential along the surface of a solute macromolecule within a surrounding solvent with a particular ionic concentration.

In the PB model, the PB equation governing electrostatic potentials takes the form of a nonlinear elliptic equation on multiple domains with discontinuous dielectric coefficients across the molecular surface or solute-solvent interface \cite{Bates1, Bates2}. The PB equation cannot be solved analytically for molecules with complex geometries, only admitting analytical solutions for shapes such as spheres or rods \cite{Holst, Kirkwood}. However, solving the PB equation numerically also presents significant difficulties because of the discontinuous dielectric coefficients, singularities in the source-term, non-smoothness of the solution, and significant nonlinearity when strong ionic effects are present.

Recently, a pseudo-transient continuation approach has been proposed for solving the nonlinear PB (NPB) equation \cite{Sayyed, Shestakov, Zhao11}, which creates a different way to tackle the nonlinear term of the NPB equation. In classical finite difference and finite element solutions of the NPB equation, a nonlinear algebraic system is typically formed through the discretization of the boundary value problem. A nonlinear relaxation method \cite{Im98,Rocchia01} or inexact Newton method \cite{Holst95} can be employed to solve such a nonlinear system. In the pseudo-transient continuation approach \cite{Sayyed, Shestakov, Zhao11}, a pseudo-time derivative is added to the NPB equation so that one solves an initial-boundary value problem now. The steady state solution of this problem gives rise to the solution to the original boundary value problem. Numerically, it is desired that a large time step can be used so that the long time integration can be computed quickly. Thus, the efficiency of a pseudo-time NPB solver is directly related to its stability, which critically depends on the nonlinear term of the NPB equation -- a hyperbolic sine function that could be exponentially large.

Several time stepping schemes have been considered for solving the time dependent NPB equation. The explicit Euler solution is straightforward, but invokes a severe stability constraint \cite{Zhao11}. Implicit time integrations have also been studied \cite{Sayyed, Shestakov}, for which care has to be exercised in handling the nonlinear term. In \cite{Sayyed}, a linearization technique based on the first order Taylor expansion is proposed so that a linear system is formed at each step of the implicit Euler integration. This linearization essentially evaluates the nonlinear term at the previous time instant. Similarly, by treating the nonlinear term explicitly, a regularized alternating direction implicit (ADI) method has been introduced in \cite{Shestakov}. Since the Thomas algorithm \cite{Mitchell} can be employed to solve the tridiagonal finite difference systems in this time splitting method, the efficiency is greatly improved. However, a very large time increment is still prohibited in these methods, because these implicit schemes are of semi-implicit nature.

More recently, we have successfully developed two fully-implicit ADI schemes \cite{Zhao13,Zhao14} for solving the time dependent NPB equation. The success lies in an analytical integration of the nonlinear term, and the use of a time splitting framework. This completely suppresses the nonlinear instability, so that these fully-implicit ADI schemes are unconditionally stable in solving benchmark problems with smooth solutions. However, for the solvation analysis in applications to real biomolecules, these ADI schemes are still conditionally stable, probably because the underlying electrostatic potentials are non-smooth. Nevertheless, these fully-implicit ADI schemes are still found to be quite efficient, due to their high temporal accuracy and much relaxed Courant-Friedrich-Lewy (CFL) constraints  \cite{Zhao13,Zhao14}.

The objective of this paper is to introduce and investigate various new operator splitting methods for solving the time dependent NPB equation. The present investigation will be conducted in the same pseudo-transient continuation framework as in  \cite{Zhao13,Zhao14}. In particular, in spatial approximation, central finite difference will be employed to discretize the nonhomogeneous diffusion operator, which yields tridiagonal matrices in all dimension-splitting schemes. In temporal approximation, the analytical integration will be carried out for the split nonlinear subsystem. The ultimate goal of this work is to eventually construct unconditionally stable operator splitting methods for solving the NPB equation in practical biomolecular applications.

In the existing fully-implicit operator splitting NPB solvers \cite{Zhao13,Zhao14}, the ADI schemes are constructed based on the classical ADI schemes \cite{DougPeace,Douglas,PeaceRach}. In the literature, there exists another family of multiplicative operator splitting methods, i.e., the fractional step methods or locally one-dimensional (LOD) methods, originally introduced by Russian mathematicians \cite{DYakonov,Yanenko63,Yanenko71}. In the present setup, the ADI discretizations can be written as some perturbations of multidimensional discretizations of the implicit methods, such as the Crank-Nicolson and backward Euler. In the proposed LOD methods, the similar splitting in alternating directions will be considered, but before the time discretization. In other words, the Crank-Nicolson and backward Euler discretizations will be conducted after the differential equations are split. Thus, the discretizations of the LOD methods are truly one-dimensional. This is different from the ADI methods, in which the right-hand side of the discretized equations contains derivatives from the other directions. Comparing with the ADI methods, the LOD methods tend to be more stable, but have a larger splitting error \cite{Mitchell}. Moreover, a previous work on the application of ADI and LOD schemes to Maxwell's equations in isotropic and lossless media found that the LOD schemes were approximately 20\% computationally less expensive than the ADI schemes \cite{Ahmed}. Other comparisons of LOD schemes and traditional ADI schemes as applied to Maxwell's equations have also found lower computational costs for the LOD methods over ADI methods \cite{Gan,Yang}. This motivates us to develop several LOD schemes for solving the time dependent NPB equation in this work.

Besides the ADI and LOD methods, there exists another family of operator splitting methods, i.e., the additive operator splitting (AOS) methods \cite{Lu91,Lu92,Tai91}. In the AOS methods, the time integration of each split subsystem can be independently carried out, while in the multiplicative operator splitting methods, the integration of the present fractional step needs the solution from the previous fractional steps. Thus, the AOS methods are well suited to be implemented in parallel processors. The stability and accuracy of the AOS methods have been analyzed in \cite{Daoud}. The application of additive operator splitting (AOS) schemes has proved popular for the purposes of nonlinear diffusion filtering for image processing \cite{Barash, Weickert}. The AOS schemes are usually efficient and stable for relatively large $\Delta t$ values, although potentially less accurate than the LOD schemes. However, their overall efficiency and simplicity merit further exploration. In the present paper, different formulations and discretizations of the AOS methods for solving time dependent NPB equation will be considered. Previous work in image processing has also explored the properties of schemes combining multiplicative and additive operator splitting, finding higher levels of accuracy with such schemes than with simpler AOS methods \cite{Barash}. We thus will also explore the construction and application of a multiplicative-additive operator splitting (MAOS) scheme for our current problem.

The rest of this paper is organized as follows. In Section \ref{sec.problem}, we first introduce the physical background of the problem. The pseudo-transient PB model will then be presented. The existing ADI schemes will be described and the remaining difficulties are discussed. The proposed operator splitting schemes will then be formulated for solving the time dependent NPB equation in Section \ref{sec.method}. Numerical validations of the proposed schemes through benchmark examples with analytical solutions are considered in Section \ref{sec.validation}, while various benchmark biological systems and subsequent free energy of solvation calculations are presented in Section \ref{sec.application}. Finally, this paper ends with a conclusion.

\section{Mathematical models and existing algorithms}\label{sec.problem}
\subsection{Governing equation}
Consider a solute macromolecule in space surrounded by a solvent aqueous solution. As defined in \cite{Richards}, the molecular surface, or solute-solvent boundary $\Gamma$, of the macromolecule divides $\mathbb{R}^3$ into the closed domain of the solute molecule $\Omega_m$ and solvent domain $\Omega_s$. The electrostatic interaction of such a solute-solvent system is governed by the nonlinear Poisson-Boltzmann (NPB) equation
\begin{equation}\label{NPB}
- \nabla \cdot ( \epsilon({\bf r}) \nabla u({\bf r}) + \bar{\kappa}^2 ({\bf r}) \sinh(u({\bf r})) = \rho_m ({\bf r}),
\end{equation} 
where $u$ is the electrostatic potential and the source term $\rho_m$ is defined as
\begin{equation}\label{RhoEq}
\rho_m({\bf r}) = \frac{4\pi e^2_c}{k_b T} \displaystyle\sum\limits_{j=1}^{N_m} q_j \delta ({\bf r} - {\bf r}_j).
\end{equation} 
The dielectric constant $\epsilon$ is piecewise such that $\epsilon({\bf r}) = \epsilon_m$ for ${\bf r} \in \Omega_m$ and $\epsilon({\bf r}) = \epsilon_s$ for ${\bf r} \in \Omega_s$. Here $k_b$ is the Boltzmann constant, $e_c$ is the fundamental charge, and $q_j$, in the same units as $e_c$, is the partial charge on the $j$th atom of the solute macromolecule locate at position ${\bf r}_j$. The modified Debye-Huckel parameter $\bar{\kappa}$ is  defined as
\begin{equation}\label{DebyeHuckel}
\bar{\kappa}^2 ({\bf r}) = \left( \frac{2N_A e^2_c}{1000k_b T} \right) I_s
\end{equation}
for ${\bf r} \in \Omega_s$ and $\bar{\kappa}=0$ for ${\bf r} \in \Omega_m$, where $N_A$ is Avogadro's Number and $I_s$ is the molar ionic strength. At room temperature ($T = 298K$), $\bar{\kappa}^2 = 8.486902807$\AA$^{-2} I_s$ from \cite{Holst}. The NPB equation varies on the choice of units, and in the present setting we utilize the dimensionless form presented in \cite{Holst95}. We may convert the dimensionless electrostatic potential $u$ to units kcal/mol/$e_c$ at room temperature ($T = 298K$) through simply multiplying the potentials by 0.592183 \cite{Holst}. 
Our numerical computations must take place on a finite domain $\Omega$. We may assign values along the boundary $\partial \Omega$ according to the approximate analytical condition
\begin{equation}\label{BoundaryCond}
 u ( {\bf r} ) = \frac{e^2_c}{k_b T} \displaystyle\sum\limits_{i=1}^{N_m} \frac{q_i}{\epsilon_s | {\bf r} - {\bf r}_i |} e^{-\bar{\kappa} \frac{| {\bf r} - {\bf r}_i |}{\sqrt{\epsilon_s}}}.
\end{equation}
When $\partial \Omega$ is of sufficient distance from the macromolecule subdomain $\Omega_m$, Eq. (\ref{BoundaryCond}) can be utilized to approximate the results for potentials found from Eq. (\ref{NPB}). We note here that Eq. (\ref{BoundaryCond}), for a collection of $N_m$ partial charges $q_i$ at positions ${\bf r}_i$, is simply a linear superposition of Coulomb's Law. For simplicity, the boundary $\partial \Omega$ is assumed to be of a cubic shape. 
\subsection{Pseudo-time dependent NPB equation}
The computational simulation of the NPB equation is of great importance for biomolecular modeling, but is quite challenging. Numerous approaches have been developed in the literature; see for example recent reviews \cite{Lu08,Li13} and references therein. What are related to the present study is the pseudo-transient continuation approach for solving the NPB equation \cite{Sayyed, Shestakov, Zhao11}. In this approach, the NPB equation (\ref{NPB}) will be converted from its time independent form to a time dependent form by introducing a pseudo-transient variation, thus producing the time dependent NPB:
\begin{equation}\label{TimeDepNPB}
\alpha \frac{\partial u}{\partial t} ({\bf r},t) = \nabla \cdot (\epsilon ({\bf r}) \nabla u({\bf r},t)) - \bar{\kappa}^2 ({\bf r}) \sinh(u({\bf r},t)) + \rho_m({\bf r}).
\end{equation}
where $\alpha$ is a scaling parameter. To solve this initial boundary value problem, we must first specify an initial solution for our potential values, which may be either the trivial solution $u = 0$, or the electrostatic potential solved from the linear Poisson-Boltzmann equation \cite{Zhao11}. Equation (\ref{TimeDepNPB}) is then numerically integrated over a sufficiently large time period to reach a steady state solution, thus essentially recovering the solution to the original NPB equation (\ref{NPB}).
\subsection{Previous schemes and remaining difficulties}
Several time integration methods have been developed for solving the time dependent NPB equation (\ref{TimeDepNPB}). A very small time increment $\Delta t$ is usually required in the explicit Euler solution \cite{Zhao11}, which is inefficient for practical use. However, the construction of a fully implicit integration of the NPB equation is hindered by the presence of the nonlinear hyperbolic sine term in the NPB equation. By evaluating the nonlinear term at the previous time step, the regular implementation of the implicit schemes usually gives rise to a semi-implicit integration \cite{Sayyed, Shestakov}.

A breakthrough has been made in our recent study \cite{Zhao14}, in which fully implicit
alternating direction implicit (ADI) methods are developed for the first time in the literature. A generalized NPB equation is considered in \cite{Zhao14}, while the similar ADI schemes for the classical NPB equation (\ref{NPB})  are introduced later in \cite{Zhao13}. The success of these ADI schemes lies in an analytical integration of the nonlinear term in a time splitting framework. With this analytical treatment of the nonlinear term, unconditionally stable results are obtained in benchmark problems with smooth solutions. Unfortunately, these operator splitting ADI schemes are still conditionally stable when they are applied to real protein systems, because of various complex features of the biomolecular solvation analysis, including complicated molecular surfaces, singular source charges, discontinuous dielectric coefficients, and nonsmooth electrostatic potentials \cite{Zhao13,Zhao14}. We thus are motivated to develop novel numerical schemes which may be both computationally less expensive as well as more stable for significantly larger $\Delta t$ values than the previously presented ADI schemes.

For a comparison, two ADI schemes presented \cite{Zhao13} will also be examined in this work. We thus briefly review  these two methods here. In the first order ADI scheme (\textit{ADI1}) \cite{Zhao13}, at each time step from $t_n$ to $t_{n+1}$, the time dependent NPB equation (\ref{TimeDepNPB}) is solved by a first order time splitting method in two stages. An analytical integration of the nonlinear term is conducted in the first stage. The nonhomogeneous diffusion equation of the second stage is first discretized by the implicit Euler scheme in time and the finite difference method in space. The multidimensional system is then solved by a classical Douglas-Rachford ADI scheme, by splitting the system into many one-dimensional (1D) systems of a tridiagonal structure.  Similarly, in the second order ADI scheme (\textit{ADI2}) \cite{Zhao13}, a second order time splitting with three stages is considered with analytical treatments in the first and last stages. The Crank-Nicholson scheme is considered in the second stage, and is realized via a Douglas ADI implementation. The \textit{ADI2} scheme is generally more accurate than the \textit{ADI1} scheme, though it is computationally more expensive \cite{Zhao13}.

\section{Towards the development of unconditionally stable schemes}\label{sec.method}
The time splitting or operator splitting methods are powerful tools for solving time dependent partial differential equations. Besides the ADI schemes, there exist many other operator splitting methods in the literature that can reduce a multidimensional problem to sets of independent 1D problems. With the central difference approximation in space, these 1D systems will have tridiagonal structures, so that they can be efficiently solved using the Thomas algorithm  \cite{Mitchell}. Several commonly used operator splitting methods will be explored in this paper, for the purpose of  eventually developing unconditionally stable schemes for solving the time dependent NPB equation. 

We adopt the following notations in this work. Following the tradition of the biomolecular simulation, we consider a uniform mesh with grid spacing $h$ in the $x$, $y$, and $z$ directions, with $N_x$, $N_y$, $N_z$ being the number of grip points in each direction, respectively. The time increment is denoted as $\Delta t$. The notation $u^n_{i,j,k} = u(x_i, y_j, z_k, t_n)$ denotes the electrostatic potential at a node $(x_i, y_j, z_k)$ and a time level $t=t_n$. Thus, the vector ${\bf U}^n = \{ u^n_{i,j,k} \}$ for $ i=1,...,N_x$, $j=1,...,N_y$, and $k=1,...,N_z$ denotes all nodal values of $u$ at  $t_n$. All proposed schemes will be focused on updating ${\bf U}^n$ at a time level $t_n$ to ${\bf U}^{n+1}$ at next time level $t_{n+1} = t_n + \Delta t$.
\subsection{Locally-one-dimensional (LOD) schemes}
The fractional step methods or locally one-dimensional (LOD) methods were first developed by  Russian mathematicians \cite{DYakonov,Yanenko63,Yanenko71}. Being multiplicative operator splitting methods too, the LOD methods adopt the similar splitting in alternating directions as in the ADI methods, but before the numerical discretization. 

In the present study, we propose several LOD schemes to solve the time dependent NPB equation (\ref{TimeDepNPB}). The spatial and temporal discretization will be detailed for one LOD scheme, while the discretization of other schemes can be similarly constructed. In particular, in the first LOD scheme, at each time step from $t_n$ to $t_{n+1}$, the time dependent NPB equation will be solved by a multiplicative operator splitting procedure consisting of five stages:
\begin{equation}\label{LOD1_1}
\alpha \frac{\partial w}{\partial t} = - \bar{\kappa}^2 \sinh(w),  
\quad \mbox{with}~ {\bf W}^n  = {\bf U}^n,
\quad t\in[t_n,t_{n+1}], 
\end{equation}
\begin{equation}\label{LOD1_2}
\alpha \frac{\partial v}{\partial t}  = \frac{\partial}{\partial x} \left( \epsilon \frac{\partial v}{\partial x} \right),
\quad \mbox{with}~ {\bf V}^{n}  = {\bf W}^{n+1},
\quad t\in[t_n,t_{n+1}], 
\end{equation}
\begin{equation}\label{LOD1_3}
\alpha \frac{\partial p}{\partial t}  = \frac{\partial}{\partial y} \left( \epsilon \frac{\partial p}{\partial y} \right),
\quad \mbox{with}~ {\bf P}^{n}  = {\bf V}^{n+1},
\quad t\in[t_n,t_{n+1}], 
\end{equation}
\begin{equation}\label{LOD1_4}
\alpha \frac{\partial q}{\partial t}  = \frac{\partial}{\partial z} \left( \epsilon \frac{\partial q}{\partial z} \right),
\quad \mbox{with}~ {\bf Q}^{n}  = {\bf P}^{n+1},
\quad t\in[t_n,t_{n+1}], 
\end{equation}
\begin{equation}\label{LOD1_5}
\alpha \frac{\partial r}{\partial t}  = \rho_m, 
\quad \mbox{with}~ {\bf R}^{n}  = {\bf Q}^{n+1},
\quad t\in[t_n,t_{n+1}], 
\end{equation}
We finally have ${\bf U}^{n+1} = {\bf R}^{n+1}$. The global splitting error of the present LOD scheme is of order one, i.e., $O(\Delta t)$. Obviously, no discretization is involved in the present splitting, which is different from the previous ADI methods \cite{Zhao13,Zhao14}. 

In the LOD schemes, the numerical approximation is introduced after the splitting. Here an analytical integration can be conducted for the nonlinear equation (\ref{LOD1_1}), as in our previous studies \cite{Zhao13,Zhao14}. In particular, Eq. (\ref{LOD1_1}) is simply a separable ordinary differential equation in this context:
\begin{equation*}
\frac{d w}{ \sinh(w)} = - \bar{\kappa}^2 \frac{dt}{\alpha} .
\end{equation*}
Thus, the integration can be carried out analytically:
\begin{equation*}
\int \frac{1}{ \sinh(w)} \,dw =
\int - \bar{\kappa}^2 \frac{dt}{\alpha} .
\end{equation*}
\begin{equation}\label{ana1}
-2 \tanh^{-1} ( e^w) = - \bar{\kappa}^2 \frac{t}{\alpha} + C
\end{equation}
Evaluating Eq. (\ref{ana1}) at both $t_n$ and $t_{n+1}$, we have
\begin{equation*}
\tanh^{-1} ( \exp ({\bf W}^{n+1})) - \tanh^{-1} ( \exp ({\bf W}^{n}))
=  \frac{1}{2\alpha} \bar{\kappa}^2 \Delta t
\end{equation*}
Taking $\tanh( \cdot)$ of both sides and re-arranging the terms, we have
\begin{equation}\label{ana2}
\exp ({\bf W}^{n+1}) = \tanh( \tanh^{-1} ( \exp ({\bf W}^{n})) 
+ \frac{1}{2\alpha} \bar{\kappa}^2 \Delta t )
\end{equation}
Equation (\ref{ana2}) can then be simplified to be
\begin{equation}\label{ana3}
{\bf W}^{n+1}  = \ln \left( 
\frac{ \cosh(\frac{1}{2\alpha} \bar{\kappa}^2 \Delta t )
+ \exp (-{\bf W}^{n}) \sinh (\frac{1}{2\alpha} \bar{\kappa}^2 \Delta t )}
{ \exp (-{\bf W}^{n}) \cosh(\frac{1}{2\alpha} \bar{\kappa}^2 \Delta t )
+ \sinh (\frac{1}{2\alpha} \bar{\kappa}^2 \Delta t )}
\right)
\end{equation}
Therefore, given ${\textbf W}^{n}$ at $t_n$, ${\textbf W}^{n+1}$ can be calculated analytically according to (\ref{ana3}) so that the difficulties associated with the $\sinh(\cdot)$ nonlinear term of the NPB equation may be bypassed.

Different LOD schemes can be obtained if we consider different time stepping procedures. We first propose the use of implicit Euler integration in time and central differencing in space to discretize Equations (\ref{LOD1_2}) through (\ref{LOD1_4}). Combining with our analytical solution for the $\sinh( \cdot )$ term (\ref{ana3}), we formulate our first scheme, \textit{LODIE1}:
\begin{equation}\label{LODIE1_1}
w_{i,j,k}  = \ln \left( 
\frac{ \cosh(\frac{1}{2\alpha} \bar{\kappa}^2 \Delta t )
+ \exp (-u^{n}_{i,j,k}) \sinh (\frac{1}{2\alpha} \bar{\kappa}^2 \Delta t )}
{ \exp (-u^{n}_{i,j,k}) \cosh(\frac{1}{2\alpha} \bar{\kappa}^2 \Delta t )
+ \sinh (\frac{1}{2\alpha} \bar{\kappa}^2 \Delta t )}
\right),
\end{equation}
\begin{equation}\label{LODIE1_2}
(1-\frac{\Delta t}{\alpha}\delta_x^2)v_{i,j,k}  = w_{i,j,k},
\end{equation}
\begin{equation}\label{LODIE1_3}
(1-\frac{\Delta t}{\alpha}\delta_y^2)p_{i,j,k}  = v_{i,j,k},
\end{equation}
\begin{equation}\label{LODIE1_4}
(1-\frac{\Delta t}{\alpha}\delta_z^2)q_{i,j,k}  = p_{i,j,k},
\end{equation}
\begin{equation}\label{LODIE1_5}
u^{n+1}_{i,j,k}  = q_{i,j,k} + \frac{\Delta t}{\alpha} Q(x_i,y_j,z_k),
\end{equation}
where $\delta_x^2$, $\delta_y^2$, and $\delta_z^2$ are the central difference operators in the $x$, $y$, and $z$ directions, respectively:
\begin{equation}\label{SpatialOperators} 
\begin{aligned}
\delta_x^2 v_{i,j,k}^{n} & = \frac{1}{h^2} \Big(
\epsilon(x_{i+\frac{1}{2}},y_j,z_k) (  v^{n}_{i+1,j,k}-  v^{n}_{i,j,k}) 
+ \epsilon(x_{i-\frac{1}{2}},y_j,z_k)(  v^{n}_{i-1,j,k}- v^{n}_{i,j,k}) \Big), \\
\delta_y^2 v_{i,j,k}^{n} & = \frac{1}{h^2} \Big(
\epsilon(x_i,y_{j+\frac{1}{2}},z_k)(  v^{n}_{i,j+1,k}-  v^{n}_{i,j,k})
+ \epsilon(x_i,y_{j-\frac{1}{2}},z_k)(  v^{n}_{i,j-1,k}-  v^{n}_{i,j,k}) \Big), \\
\delta_z^2 v_{i,j,k}^{n} & = \frac{1}{h^2} \Big(
\epsilon(x_i,y_j,z_{k+\frac{1}{2}})(  v^{n}_{i,j,k+1}-  v^{n}_{i,j,k})
+ \epsilon(x_i,y_j,z_{k-\frac{1}{2}})(  v^{n}_{i,j,k-1}-  v^{n}_{i,j,k}) \Big).
\end{aligned}
\end{equation}
Furthermore, $Q(x_i,y_j,z_k)$ is the distribution of all source charges in the source term $\rho_m$ from Equation (\ref{RhoEq}), distributed by a trilinear interpolation. The value of $\epsilon$ is determined by its location which is either on/inside or outside the molecular surface $\Gamma$. Specifically, $\epsilon (x,y,z) = \epsilon_m$ if $(x,y,z) \in \Omega_m$, and $\epsilon (x,y,z) = \epsilon_s$ if $(x,y,z) \in \Omega_s$, where the molecular surface $\Gamma$ is determined by a commonly used software package: MSMS \cite{Sanner}. Both $Q(x_i,y_j,z_k)$ and $\epsilon(x_i,y_j,z_k)$ are time independent and need to be calculated only once at the beginning. We utilize Dirichlet boundary conditions in our numerical simulations, and we assume the same boundary values for $w$, $v$, $p$, $q$, and $u$. We utilize similar assumptions for our numerical simulations with all following methods.

We further formulate a second scheme based on the LOD1 scheme that utilizes Crank-Nicholson integration in time and central differencing in space. This scheme, \textit{LODCN1}, replaces Equations (\ref{LODIE1_2}) through (\ref{LODIE1_4}) from \textit{LODIE1} with the following:
\begin{equation}\label{LODCN1_2}
(1-\frac{1}{2\alpha} \Delta t\delta_x^2)v_{i,j,k}  = (1+\frac{1}{2\alpha} \Delta t\delta_x^2)w_{i,j,k},
\end{equation}
\begin{equation}\label{LODCN1_3}
(1-\frac{1}{2\alpha} \Delta t\delta_y^2)p_{i,j,k}  = (1+\frac{1}{2\alpha} \Delta t\delta_y^2)v_{i,j,k},
\end{equation}
\begin{equation}\label{LODCN1_4}
(1-\frac{1}{2\alpha} \Delta t\delta_z^2)q_{i,j,k}  = (1+\frac{1}{2\alpha} \Delta t\delta_z^2)p_{i,j,k}
\end{equation}

Minor changes can be resulted if we alter the order of the subsystems in the LOD splitting. To illustrate this, we also propose a second set of LOD schemes that solve the five stages presented in (\ref{LOD1_1}) through (\ref{LOD1_5}) in a modified order: the source term equation in (\ref{LOD1_5}) is treated first, while the nonlinear equation in (\ref{LOD1_1}) is treated last. This LOD method, similar to the first one, has a form utilizing implicit Euler integration in time and another form utilizing Crank-Nicholson integration in time. \textit{LODIE2}, then, simply consists of equations (\ref{LODIE1_5}), (\ref{LODIE1_2}), (\ref{LODIE1_3}), (\ref{LODIE1_4}), (\ref{LODIE1_1}), respectively. \textit{LODCN2} consists of (\ref{LODIE1_5}), (\ref{LODCN1_2}), (\ref{LODCN1_3}), (\ref{LODCN1_4}), (\ref{LODIE1_1}), respectively.

We conclude this subsection by presenting some theoretical results of the proposed LOD schemes. 
\begin{remark}
The proposed \textit{LODIE1}, \textit{LODIE2}, \textit{LODCN1}, and \textit{LODCN2} schemes are of first order accuracy in time, because the underlying LOD splittings are first order and the discretization error of the implicit Euler or Crank-Nicolson integration is at least first order. 
\end{remark}
\begin{remark}
If the solution $u$ is sufficiently smooth, the proposed \textit{LODIE1}, \textit{LODIE2}, \textit{LODCN1}, and \textit{LODCN2} schemes are unconditionally stable, because each individual implicit Euler or Crank-Nicolson time integration with central difference approximation is unconditionally stable.
\end{remark}

\subsection{Additive operator splitting (AOS) schemes}
We next propose a series of schemes utilizing an additive operator splitting (AOS) formulation \cite{Lu91,Lu92,Tai91}. Unlike the ADI and LOD schemes, in which the subsystems have to be solved sequentially, the split equations can be solved concurrently in the AOS schemes. Thus, the AOS methods are well suited to be implemented in parallel processors. Since the present work focuses on the stability investigation, a simple series implementation is still conducted for the AOS schemes.

In the first AOS scheme, at each time step from $t_n$ to $t_{n+1}$, we solve the time dependent NPB equation in the following stages:
\begin{equation}\label{AOS1_1}
\alpha \frac{\partial w}{\partial t}  = - 4 \bar{\kappa}^2  \sinh(w),  
\quad \mbox{with}~ {\bf W}^n  = {\bf U}^n,
\quad t\in[t_n,t_{n+1}], 
\end{equation}
\begin{equation}\label{AOS1_2}
\alpha \frac{\partial v}{\partial t}  = 4\frac{\partial}{\partial x} \left( \epsilon \frac{\partial v}{\partial x} \right)+4 \rho_m,
\quad \mbox{with}~ {\bf V}^{n}  = {\bf U}^{n},
\quad t\in[t_n,t_{n+1}], 
\end{equation}
\begin{equation}\label{AOS1_3}
\alpha \frac{\partial p}{\partial t}  = 4\frac{\partial}{\partial y} \left( \epsilon \frac{\partial p}{\partial y} \right),
\quad \mbox{with}~ {\bf P}^{n}  = {\bf U}^{n},
\quad t\in[t_n,t_{n+1}], 
\end{equation}
\begin{equation}\label{AOS1_4}
\alpha \frac{\partial q}{\partial t}  = 4\frac{\partial}{\partial z} \left( \epsilon \frac{\partial q}{\partial z} \right),
\quad \mbox{with}~ {\bf Q}^{n}  = {\bf U}^{n},
\quad t\in[t_n,t_{n+1}], 
\end{equation}
\begin{equation}\label{AOS1_5}
{\bf U}^{n+1}  = \frac{1}{4}\left({\bf W}^{n+1}+{\bf V}^{n+1}+{\bf P}^{n+1}+{\bf Q}^{n+1}\right).
\end{equation}
In other words, the NPB system is split into four parts, i.e., one nonlinear subsystem and three subsystems along Cartesian directions. The time independent source term is assigned to the $x$-direction subsystem. Then an arithmetic average is carried out to advance ${\bf U}^{n}$ to ${\bf U}^{n+1}$. 

Similar to our LOD schemes, we will formulate two types of time discretizations. Using implicit Euler integration in time and central differencing in space to discretize Equations (\ref{AOS1_2}) through (\ref{AOS1_4}), we propose first the scheme \textit{AOSIE1}:
\begin{equation}\label{AOSIE1_1}
w_{i,j,k}  = 4\ln \left( 
\frac{ \cosh(\frac{1}{2\alpha} \bar{\kappa}^2 \Delta t )
+ \exp (-u^{n}_{i,j,k}) \sinh (\frac{1}{2\alpha} \bar{\kappa}^2 \Delta t )}
{ \exp (-u^{n}_{i,j,k}) \cosh(\frac{1}{2\alpha} \bar{\kappa}^2 \Delta t )
+ \sinh (\frac{1}{2\alpha} \bar{\kappa}^2 \Delta t )}
\right),
\end{equation}
\begin{equation}\label{AOSIE1_2}
(1-\frac{4}{\alpha}\Delta t\delta_x^2)v_{i,j,k}  = u^n_{i,j,k}+\frac{4}{\alpha}\Delta t \rho_m,
\end{equation}
\begin{equation}\label{AOSIE1_3}
(1-\frac{4}{\alpha}\Delta t\delta_y^2)p_{i,j,k}  = u^n_{i,j,k},
\end{equation}
\begin{equation}\label{AOSIE1_4}
(1-\frac{4}{\alpha}\Delta t\delta_z^2)q_{i,j,k}  = u^n_{i,j,k},
\end{equation}
\begin{equation}\label{AOSIE1_5}
u^{n+1}_{i,j,k} =\frac{1}{4}\left(w_{i,j,k}+v_{i,j,k}+p_{i,j,k}+q_{i,j,k}\right)
\end{equation}
We also formulate a second scheme based on the AOS1 scheme that utilizes Crank-Nicholson integration in time and central differencing in space. This scheme, \textit{AOSCN1}, replaces Equations (\ref{AOSIE1_2}) through (\ref{AOSIE1_4}) from \textit{AOSIE1} with the following:
\begin{equation}\label{AOSCN1_2}
(1-\frac{2}{\alpha} \Delta t\delta_x^2)v_{i,j,k}  = (1+\frac{2}{\alpha} \Delta t\delta_x^2)u^n_{i,j,k}+\frac{4}{\alpha}\Delta t \rho_m,
\end{equation}
\begin{equation}\label{AOSCN1_3}
(1-\frac{2}{\alpha} \Delta t\delta_y^2)p_{i,j,k}  = (1+\frac{2}{\alpha} \Delta t\delta_y^2)u^n_{i,j,k},
\end{equation}
\begin{equation}\label{AOSCN1_4}
(1-\frac{2}{\alpha} \Delta t\delta_z^2)q_{i,j,k}  = (1+\frac{2}{\alpha} \Delta t\delta_z^2)u^n_{i,j,k}
\end{equation}

For the AOS schemes, since each subsystem is solved independently, any change in their order will not alter the numerical outcome. In the present study, we consider a different variation for the proposed AOS schemes. In the \textit{AOSCN1} scheme, the source term is imposed only along the $x$ direction. We are interested in a more symmetric splitting by introducing a second set of AOS schemes. At each time step from $t_n$ to $t_{n+1}$, we solve the time dependent NPB equation in the following stages:
\begin{equation}\label{AOS2_1}
\alpha \frac{\partial w}{\partial t} = - 4 \bar{\kappa}^2 \sinh(w),  
\quad \mbox{with}~ {\bf W}^n  = {\bf U}^n,
\quad t\in[t_n,t_{n+1}], 
\end{equation}
\begin{equation}\label{AOS2_2}
\alpha \frac{\partial v}{\partial t}  = 4\frac{\partial}{\partial x} \left( \epsilon \frac{\partial v}{\partial x} \right)+\frac{4}{3} \rho_m,
\quad \mbox{with}~ {\bf V}^{n}  = {\bf U}^{n},
\quad t\in[t_n,t_{n+1}], 
\end{equation}
\begin{equation}\label{AOS2_3}
\alpha \frac{\partial p}{\partial t}  = 4\frac{\partial}{\partial y} \left( \epsilon \frac{\partial p}{\partial y} \right)+\frac{4}{3} \rho_m,
\quad \mbox{with}~ {\bf P}^{n}  = {\bf U}^{n},
\quad t\in[t_n,t_{n+1}], 
\end{equation}
\begin{equation}\label{AOS2_4}
\alpha \frac{\partial q}{\partial t}  = 4\frac{\partial}{\partial z} \left( \epsilon \frac{\partial q}{\partial z} \right)+\frac{4}{3} \rho_m,
\quad \mbox{with}~ {\bf Q}^{n}  = {\bf U}^{n},
\quad t\in[t_n,t_{n+1}], 
\end{equation}
\begin{equation}\label{AOS2_5}
{\bf U}^{n+1}  = \frac{1}{4}\left({\bf W}^{n+1}+{\bf V}^{n+1}+{\bf P}^{n+1}+{\bf Q}^{n+1}\right)
\end{equation}
It is obvious that we split the source term into three Cartesian subsystems so that certain symmetry is maintained in three directions. It is noted that we do not
distribute the source term into the nonlinear subsystem, because the time integration of the resulting nonlinear process becomes quite involved. 

We similarly formulate two schemes based on the AOS2 method, one utilizing implicit Euler integration in time and another utilizing Crank-Nicholson integration in time. Both of these schemes largely take a similar form as the \textit{AOSIE1} scheme; for \textit{AOSIE2}, Equations (\ref{AOSIE1_2}) through (\ref{AOSIE1_4}) are replaced with the following:
\begin{equation}\label{AOSIE2_2}
(1-\frac{4}{\alpha}\Delta t\delta_x^2)v_{i,j,k}  = u^n_{i,j,k}+\frac{4}{3\alpha}\Delta t  \rho_m,
\end{equation}
\begin{equation}\label{AOSIE2_3}
(1-\frac{4}{\alpha}\Delta t\delta_y^2)p_{i,j,k}  = u^n_{i,j,k}+\frac{4}{3\alpha}\Delta t  \rho_m,
\end{equation}
\begin{equation}\label{AOSIE2_4}
(1-\frac{4}{\alpha}\Delta t\delta_z^2)q_{i,j,k}  = u^n_{i,j,k}+\frac{4}{3\alpha}\Delta t  \rho_m,
\end{equation}
and similarly for the \textit{AOSCN2} scheme, we utilize the following series of equations:
\begin{equation}\label{AOSCN2_2}
(1-\frac{2}{\alpha} \Delta t\delta_x^2)v_{i,j,k}  = (1+\frac{2}{\alpha} \Delta t\delta_x^2)u^n_{i,j,k}+\frac{4}{3\alpha}\Delta t  \rho_m,
\end{equation}
\begin{equation}\label{AOSCN2_3}
(1-\frac{2}{\alpha} \Delta t\delta_y^2)p_{i,j,k}  = (1+\frac{2}{\alpha} \Delta t\delta_y^2)u^n_{i,j,k}+\frac{4}{3\alpha}\Delta t  \rho_m,
\end{equation}
\begin{equation}\label{AOSCN2_4}
(1-\frac{2}{\alpha} \Delta t\delta_z^2)q_{i,j,k}  = (1+\frac{2}{\alpha} \Delta t\delta_z^2)u^n_{i,j,k}+\frac{4}{3\alpha}\Delta t  \rho_m.
\end{equation}

The error and stability analyses of the general AOS schemes can be found in \cite{Daoud}. For the present AOS schemes, we have the following results. 
\begin{remark}
The proposed \textit{AOSIE1}, \textit{AOSIE2}, \textit{AOSCN1}, and \textit{AOSCN2} schemes are of first order accuracy in time, because the underlying AOS splittings are first order and the discretization error of the implicit Euler or Crank-Nicolson integration is at least first order. 
\end{remark}
\begin{remark}
If the solution $u$ is sufficiently smooth, the proposed \textit{AOSIE1}, \textit{AOSIE2}, \textit{AOSCN1}, and \textit{AOSCN2} schemes are unconditionally stable, because each individual implicit Euler or Crank-Nicolson time integration with central difference approximation is unconditionally stable.
\end{remark}

\subsection{Multiplicative-additive operator splitting (MAOS) schemes}
Finally, we propose two hybrid schemes combining multiplicative and additive operator splitting stages. Our intention is to treat the nonlinear subsystem separately from the three linear Cartesian subsystems. So, at each time step from $t_n$ to $t_{n+1}$, the time dependent NPB equation will be solved by a two-stage multiplicative operator splitting scheme, while the second of which is solved with an additive operator splitting scheme:
\begin{equation}\label{MAOS1_1}
\alpha \frac{\partial w}{\partial t}  = - \bar{\kappa}^2  \sinh(w),  
\quad \mbox{with}~ {\bf W}^n  = {\bf U}^n,
\quad t\in[t_n,t_{n+1}], 
\end{equation}
\begin{equation}\label{MAOS1_2}
\alpha \frac{\partial v}{\partial t}  = 3\frac{\partial}{\partial x} \left( \epsilon \frac{\partial v}{\partial x} \right)+ \rho_m,
\quad \mbox{with}~ {\bf V}^{n}  = {\bf W}^{n},
\quad t\in[t_n,t_{n+1}], 
\end{equation}
\begin{equation}\label{MAOS1_3}
\alpha \frac{\partial p}{\partial t}  = 3\frac{\partial}{\partial y} \left( \epsilon \frac{\partial p}{\partial y} \right)+ \rho_m,
\quad \mbox{with}~ {\bf P}^{n}  = {\bf W}^{n},
\quad t\in[t_n,t_{n+1}], 
\end{equation}
\begin{equation}\label{MAOS1_4}
\alpha \frac{\partial q}{\partial t}  = 3\frac{\partial}{\partial z} \left( \epsilon \frac{\partial q}{\partial z} \right)+ \rho_m,
\quad \mbox{with}~ {\bf Q}^{n}  = {\bf W}^{n},
\quad t\in[t_n,t_{n+1}], 
\end{equation}
\begin{equation}\label{MAOS1_5}
{\bf U}^{n+1}  = \frac{1}{3}\left({\bf V}^{n+1}+{\bf P}^{n+1}+{\bf Q}^{n+1}\right)
\end{equation}
Again, just as with the previous schemes, this MAOS method is discretized using implicit Euler integration in time, forming the \textit{MAOSIE} scheme, and also discretized using Crank-Nicholson integration in time, forming the \textit{MAOSCN} scheme.

By combining the previous studies, we have the following results for the MAOS schemes. 
\begin{remark}
The proposed \textit{MAOSIE} and \textit{MAOSCN} schemes are of first order accuracy in time, because the underlying MAOS splitting is first order and the discretization error of the implicit Euler or Crank-Nicolson integration is at least first order. 
\end{remark}
\begin{remark}
If the solution $u$ is sufficiently smooth, \textit{MAOSIE} and \textit{MAOSCN} schemes are unconditionally stable, because each individual implicit Euler or Crank-Nicolson time integration with central difference approximation is unconditionally stable.
\end{remark}

\section{Numerical validation}\label{sec.validation}
In this section, we validate the proposed schemes numerically by solving the NPB equation on a sphere -- a case with an analytical solution. We will explore the stability as well as spatial and temporal convergence of the proposed time splitting schemes and compare these results to our previous ADI methods. Such studies will help us to identify  well-performed NPB solvers to be used in real biomolecular simulations. All simulations are compiled with the Intel Fortran Compiler and run on an early-2011 MacBook Pro with an i7-2820QM 2.3GHz GPU and 8GB memory.

\subsection{Benchmark problem}

The NPB equation only possesses analytical solutions for certain simple geometries, such as a sphere. For verification of our schemes, we introduce the following solution for the case of a spherical cavity, based on a case from \cite{Sauter}:
\begin{equation}\label{NPB_analytical_1}
u({\bf r}) = \begin{cases} \frac{1}{\epsilon R} - \frac{1}{R} + \frac{1}{\| {\bf r} \|}, & \| {\bf r} \| < R \\
	\frac{1}{\epsilon \| {\bf r} \|}, & \| {\bf r} \| > R \end{cases}
\end{equation}
\begin{equation}\label{NPB_analytical_2}
\rho_m({\bf r}) = \begin{cases} 4\pi \epsilon_m \delta ({\bf r}), & \| {\bf r} \| < R \\
	\bar{\kappa}^2 \sinh \left( \frac{1}{\epsilon \| {\bf r} \|} \right), & \| {\bf r} \| > R \end{cases}
\end{equation}
where $\epsilon = \epsilon_s/\epsilon_m$ and $R$ is the radius of the spherical cavity. We note here that because of the singularity in the source term defined by (\ref{NPB_analytical_2}) and the non-smoothness due to the interface jump conditions, the accuracy of the finite difference spatial discretization is normally reduced and additional instability may be introduced in the time-stepping. In our numerical validations, we choose $R=1${\AA} with a single centered charge of 1 $e_c$. Our dielectric constants are $\epsilon_s=80$ and $\epsilon_m=1$, and we set our nonlinear constant $\bar{\kappa}=1$. Furthermore, for our scaling parameter, we select the standard value $\alpha=1$.

A cubic domain $[-3,3] \times [-3,3] \times [-3,3]$ with the same spacing in all three directions $h=\Delta x = \Delta y = \Delta z$ is used in our computations. A Dirichlet boundary condition is assumed on all boundaries with the boundary data being given by the analytical solution. For the purpose of investigating the convergence and stability of our new schemes, we will examine the sensitivity of the methods to the nonlinear term of the NPB. To this end, we construct a family of nontrivial initial solutions:
\begin{equation}\label{InitSolutions}
u({\bf r},0) = H \cos(\frac{\pi}{6} x) \cos(\frac{\pi}{6} y) \cos(\frac{\pi}{6} z) + \frac{1}{\epsilon \sqrt{x^2 + y^2 + z^2}}
\end{equation}
in which we employ different magnitudes of $H$. Note that this initial solution satisfies the Dirichlet boundary condition. With such an initial solution at $t=0$, the time-stepping will be carried out until a stopping time $t=T$ with a time increment $\Delta t$. Denoting $u_h$ as the numerical solution, the following measures are used to estimate relative errors:
\begin{equation}\label{ErrorMeasure}
L_{\infty} = \frac{\max | u-u_h |}{\max | u |} , \ L_2 = \sqrt{\frac{\Sigma_{i,j,k} | u-u_h |^2}{\Sigma_{i,j,k} |u|^2}}.
\end{equation}

\subsection{Stability}
We first investigate the stability. For smooth solutions, the proposed time splitting methods are unconditionally stable. Nevertheless, the present analytical solution $u$ is just $C^0$ continuous across the circular interface. It is thus of interest to test the stability of our newly constructed schemes for nonsmooth solutions. 

Previously, the stability of several explicit and implicit schemes has been examined for solving the time dependent NPB equation, by using this benchmark example \cite{Zhao13}. All tested schemes in that paper were found to be stable with a constraint $\Delta t \le \frac{h^2}{m}$ for some $m$. However, when the parameter $H$ changes from $1$ to $20$, which represents a strong nonlinearity, all explicit and semi-implicit time integration schemes become unconditionally unstable, while the fully implicit ADI schemes still remain to be stable with a similar $m$ value bounded by $m < 20$ \cite{Zhao13}.

\begin{figure}[!tb]
\centering
\begin{tabular}{cc}
\includegraphics[width=0.5\linewidth]{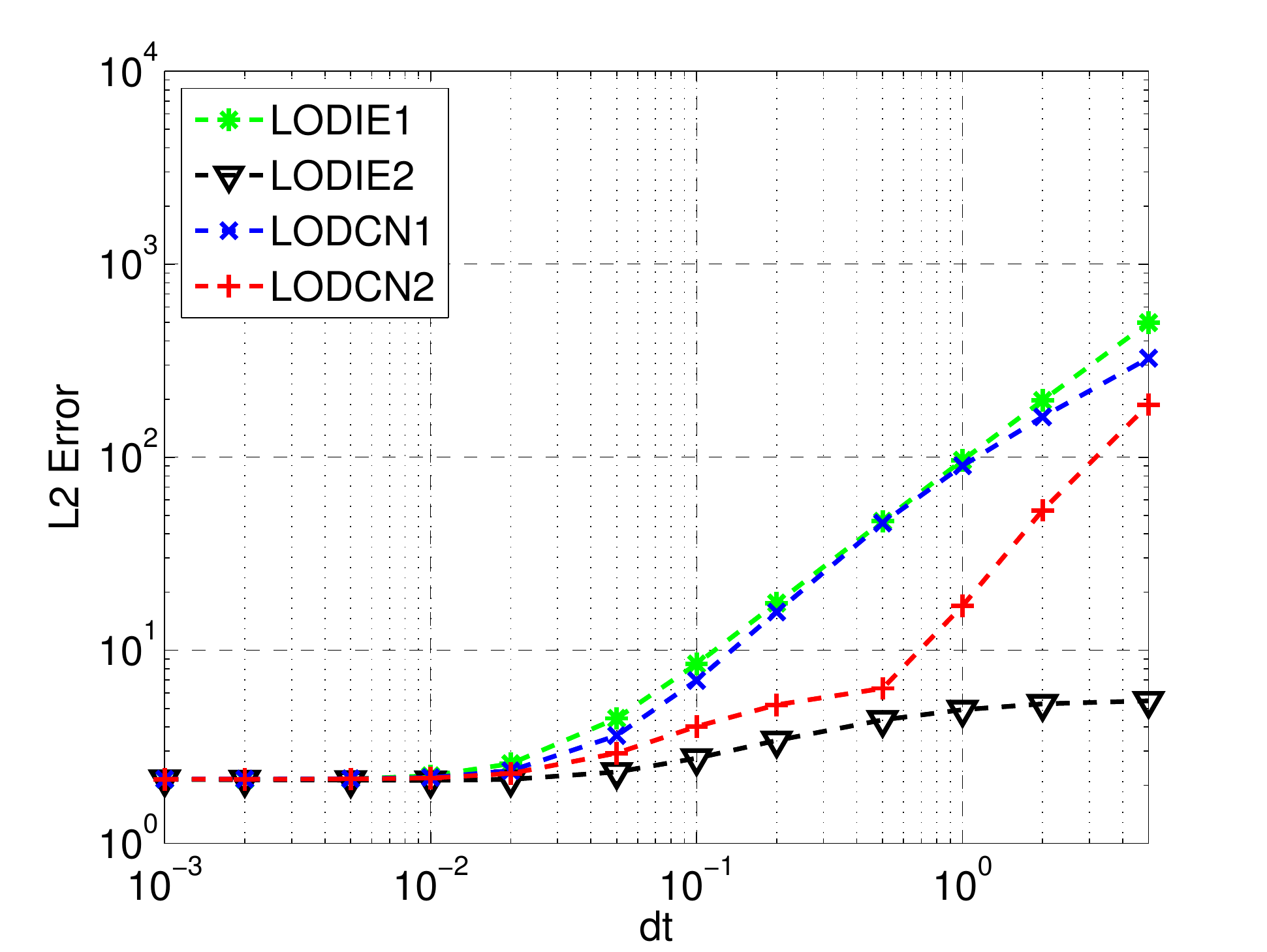} &
\includegraphics[width=0.5\linewidth]{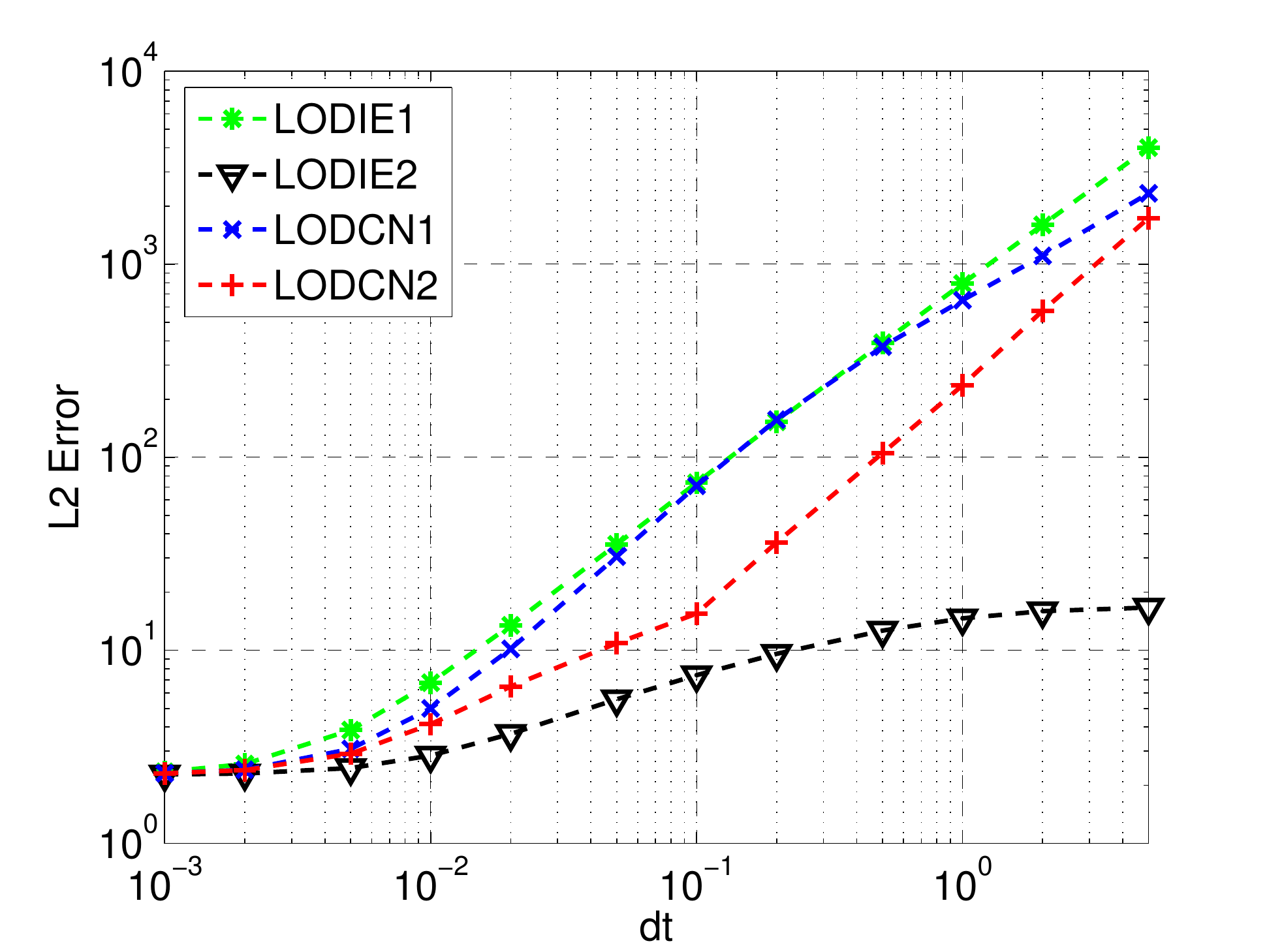}\\
(a) & (b)
\end{tabular}
\caption{Stability verification of the LOD schemes with $H=1$. (a) $h = 0.5$; (b) $h = 0.25$.}
\label{fig.LOD_AnaStabVer1}
\end{figure}
\begin{figure}[!tb]
\centering
\begin{tabular}{cc}
\includegraphics[width=0.5\linewidth]{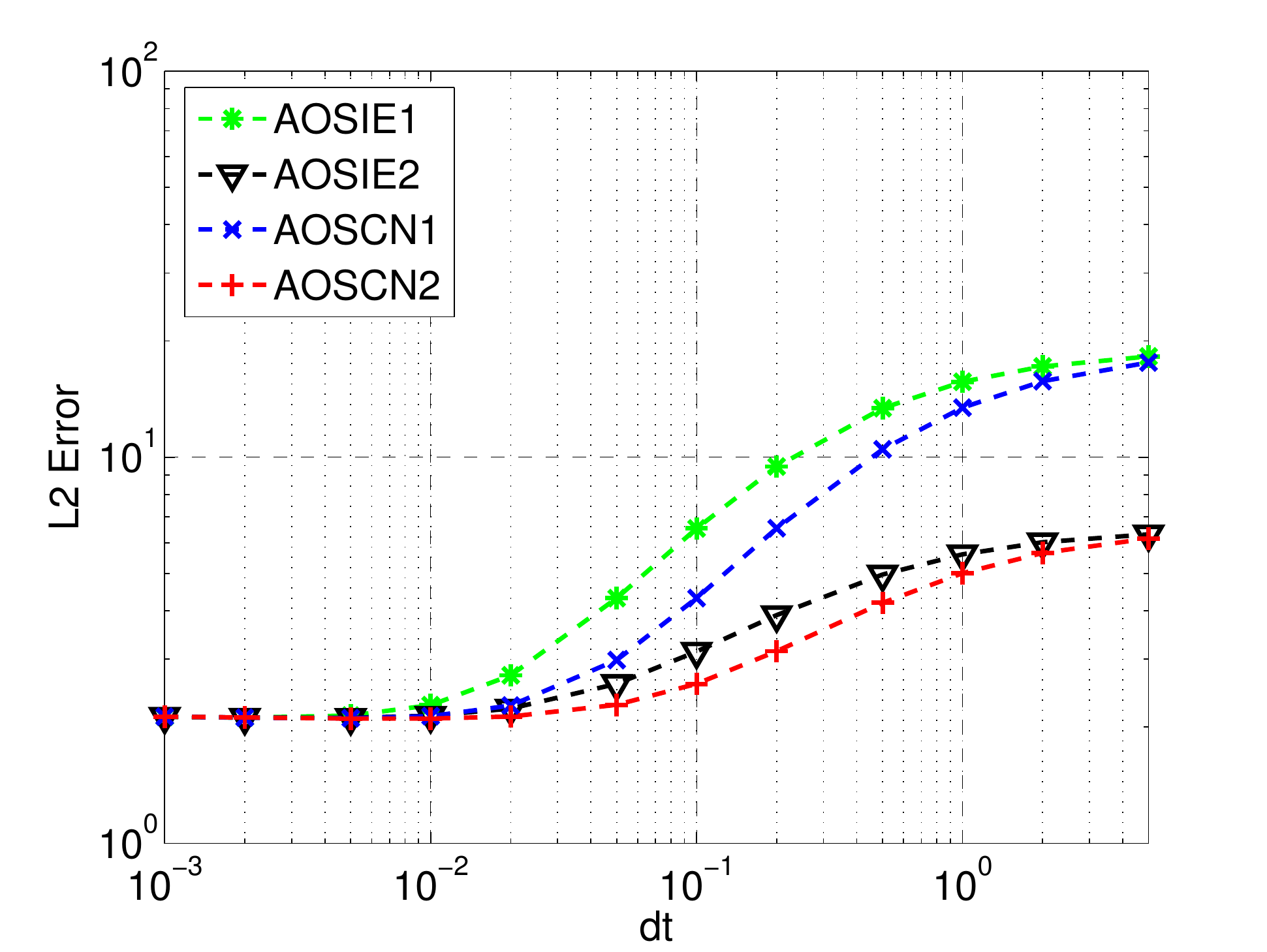} &
\includegraphics[width=0.5\linewidth]{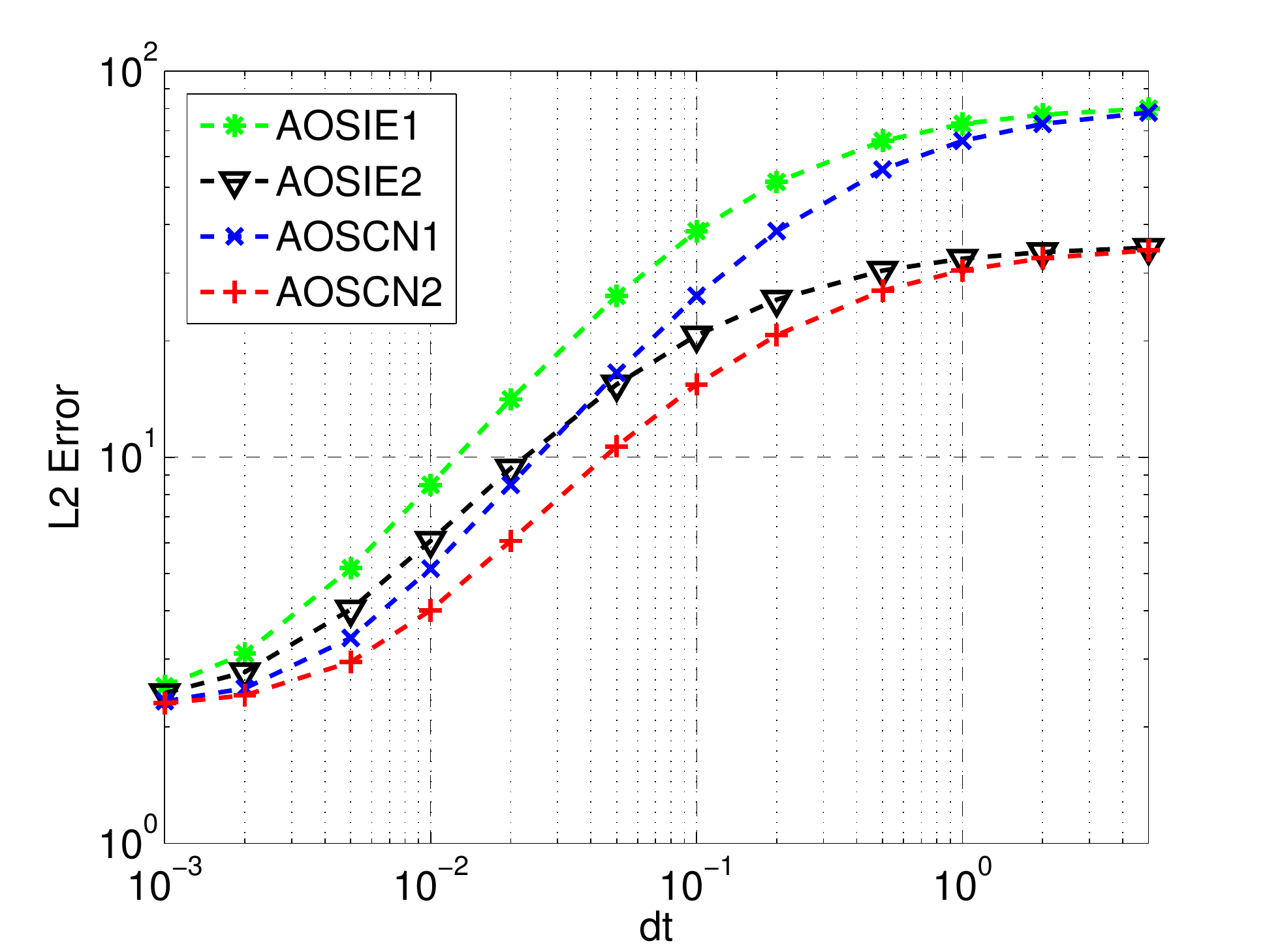}\\
(a) & (b)
\end{tabular}
\caption{Stability verification of the AOS schemes with $H=1$. (a) $h=0.5$; (b) $h=0.25$.}
\label{fig.AOS_AnaStabVer1}
\end{figure}
\begin{figure}[!tb]
\centering
\begin{tabular}{cc}
\includegraphics[width=0.5\linewidth]{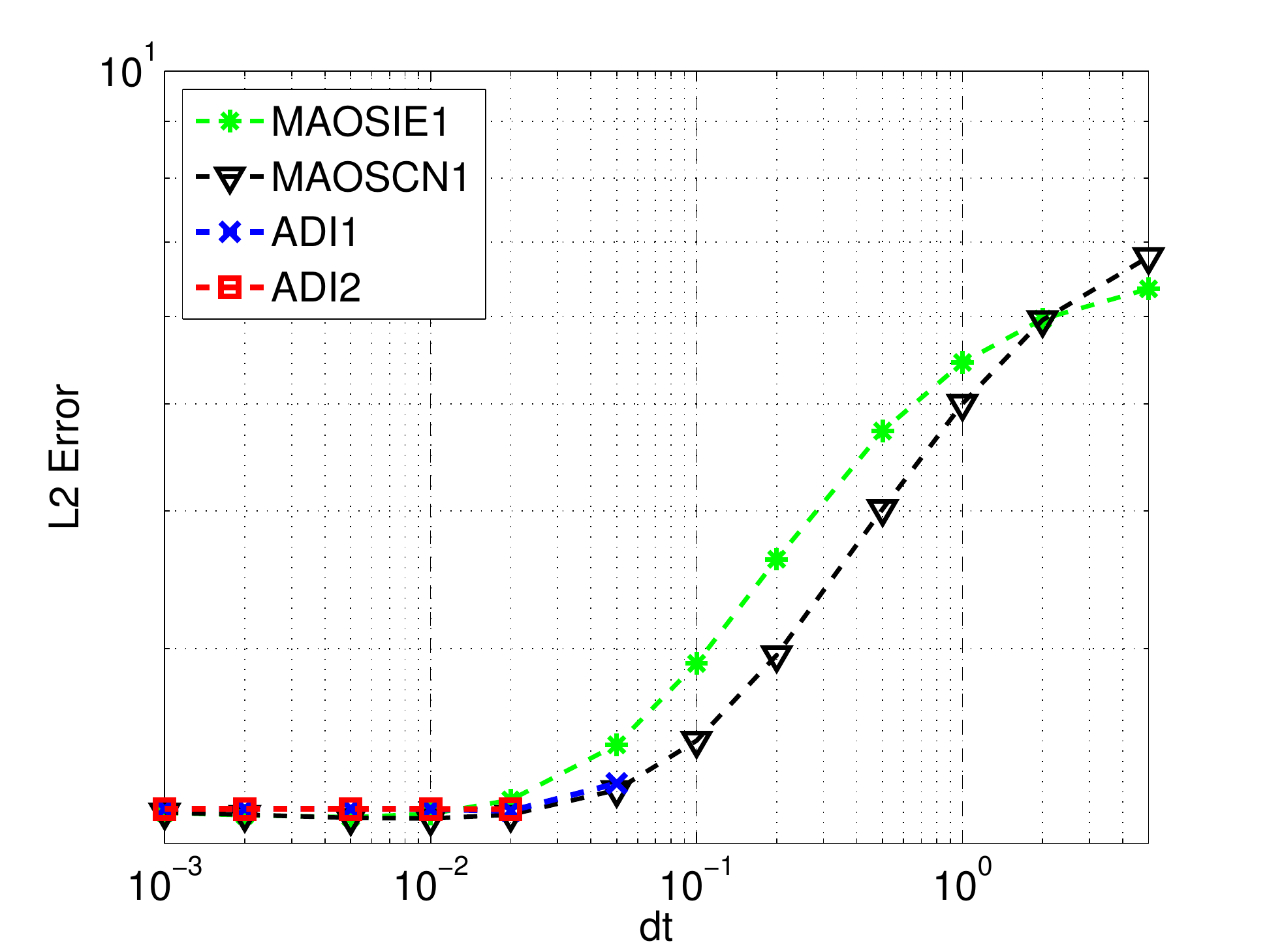} &
\includegraphics[width=0.5\linewidth]{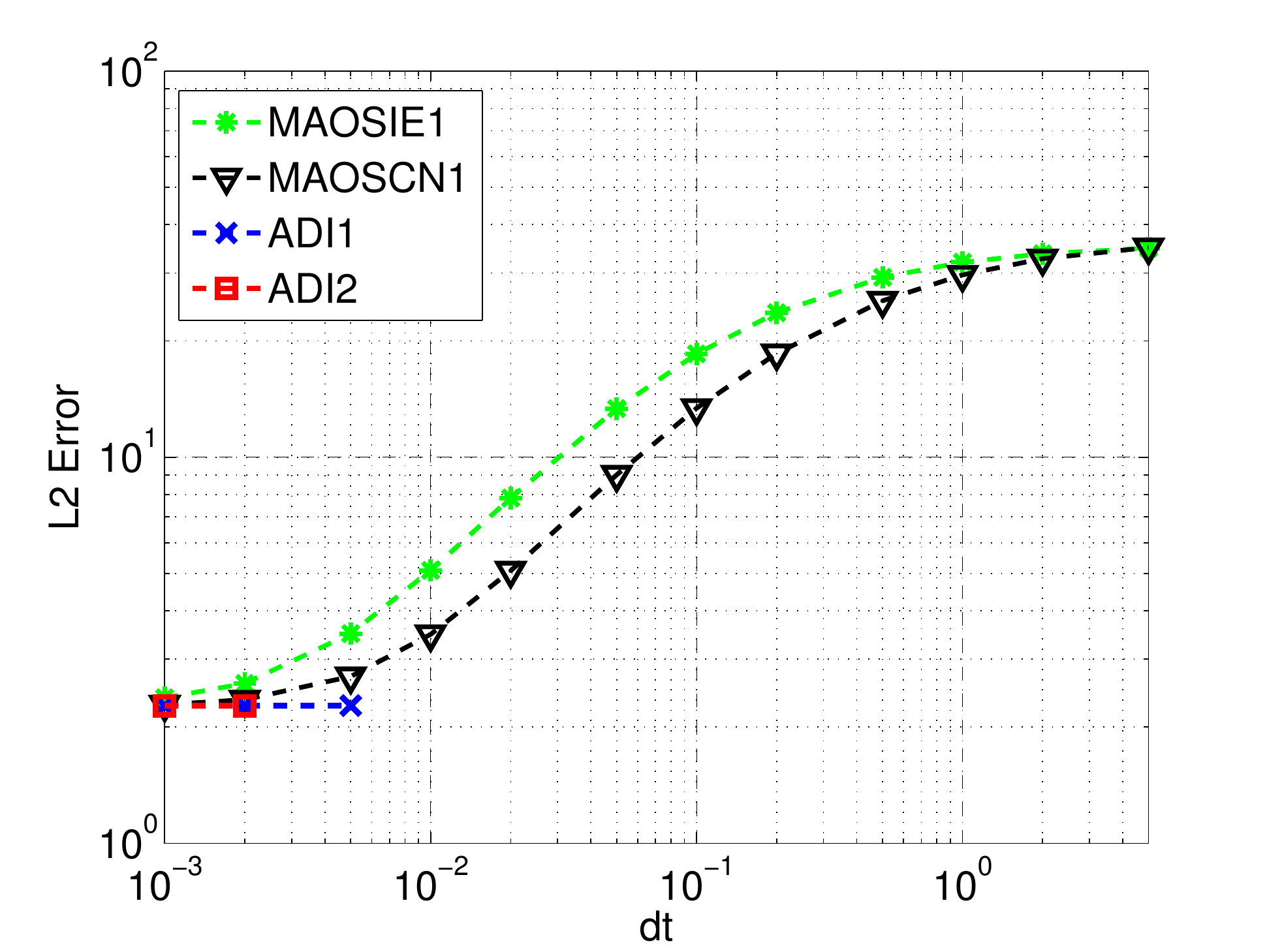}\\
(a) & (b)
\end{tabular}
\caption{Stability verification of the MAOS and ADI schemes  with $H=1$. (a) $h=0.5$; (b) $h=0.25$.}
\label{fig.MAOS_ADI_AnaStabVer1}
\end{figure}
\begin{figure}[!tb]
\centering
\begin{tabular}{cc}
\includegraphics[width=0.5\linewidth]{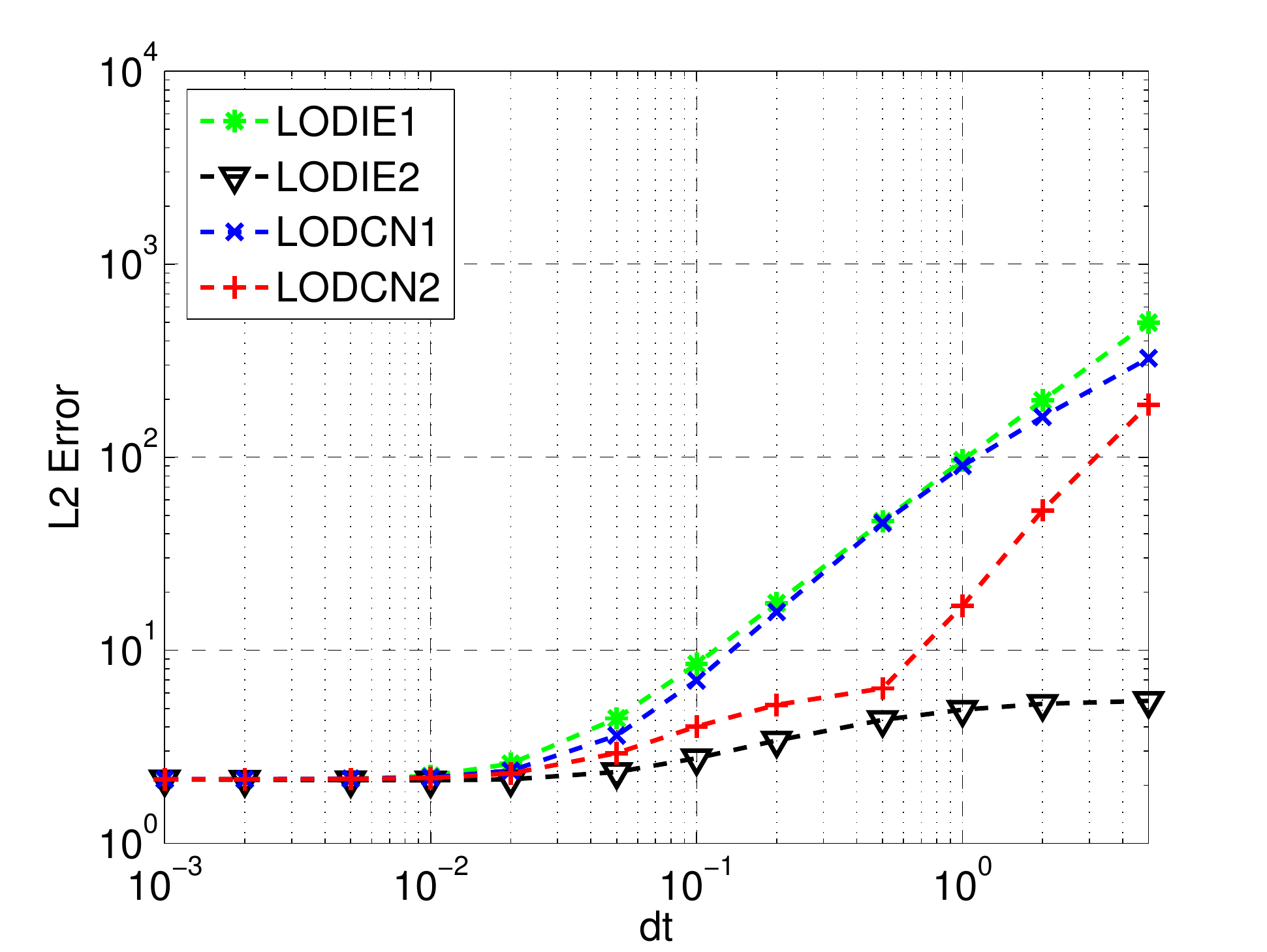} &
\includegraphics[width=0.5\linewidth]{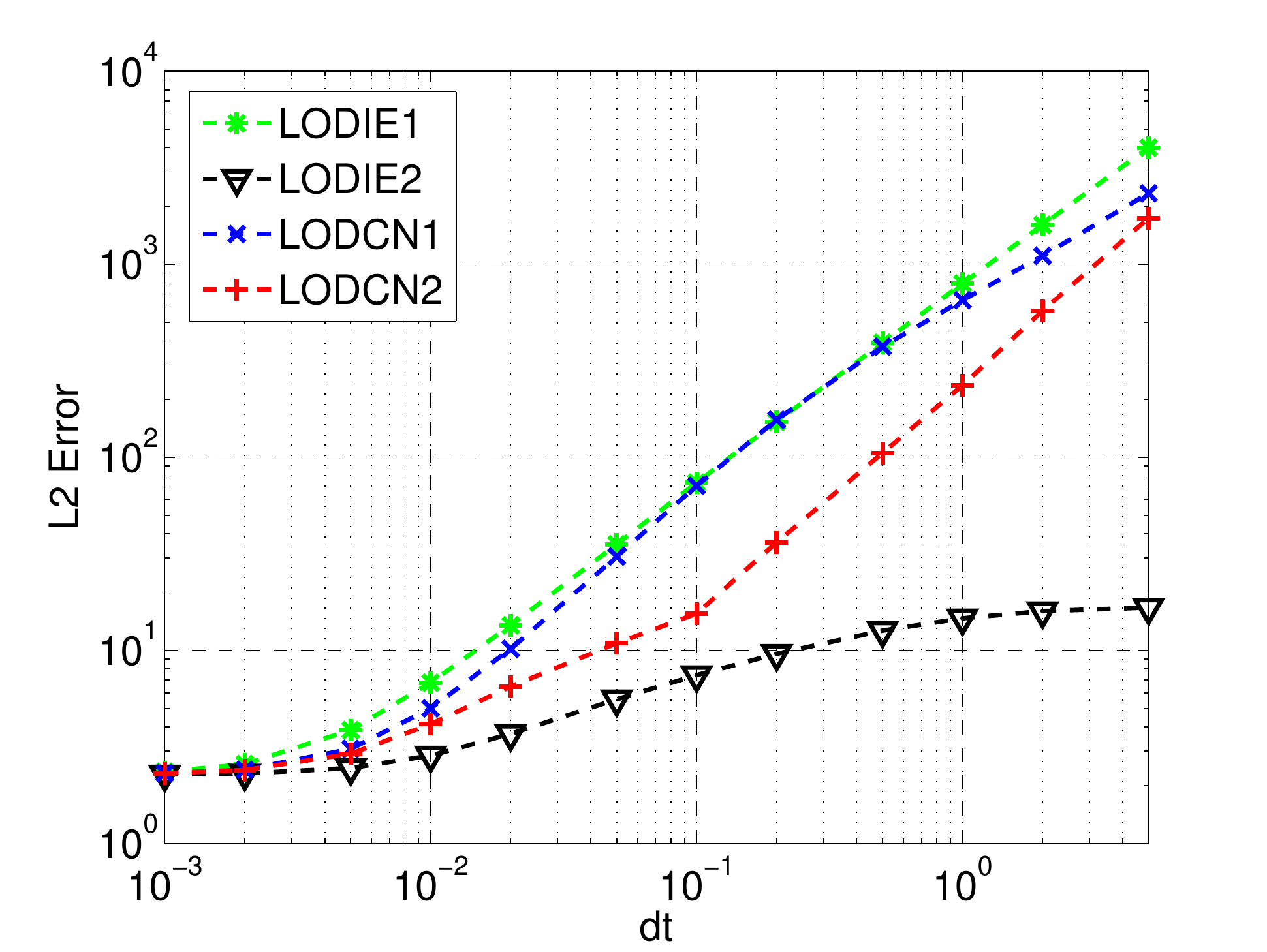}\\
(a) & (b)
\end{tabular}
\caption{Stability verification of the LOD schemes with $H=20$. (a) $h = 0.5$; (b) $h = 0.25$. }
\label{fig.LOD_AnaStabVer20}
\end{figure}
\begin{figure}[!tb]
\centering
\begin{tabular}{cc}
\includegraphics[width=0.5\linewidth]{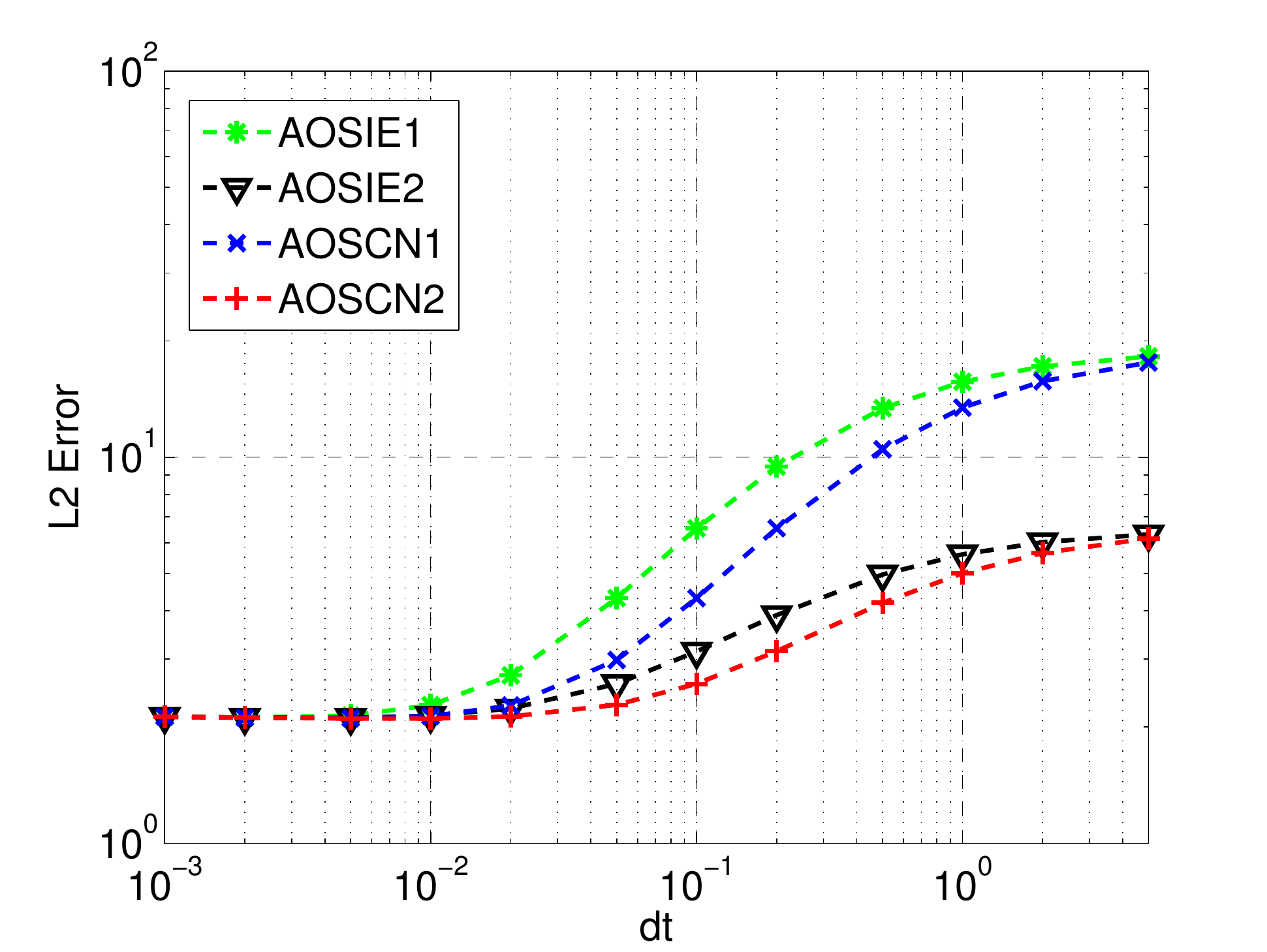} &
\includegraphics[width=0.5\linewidth]{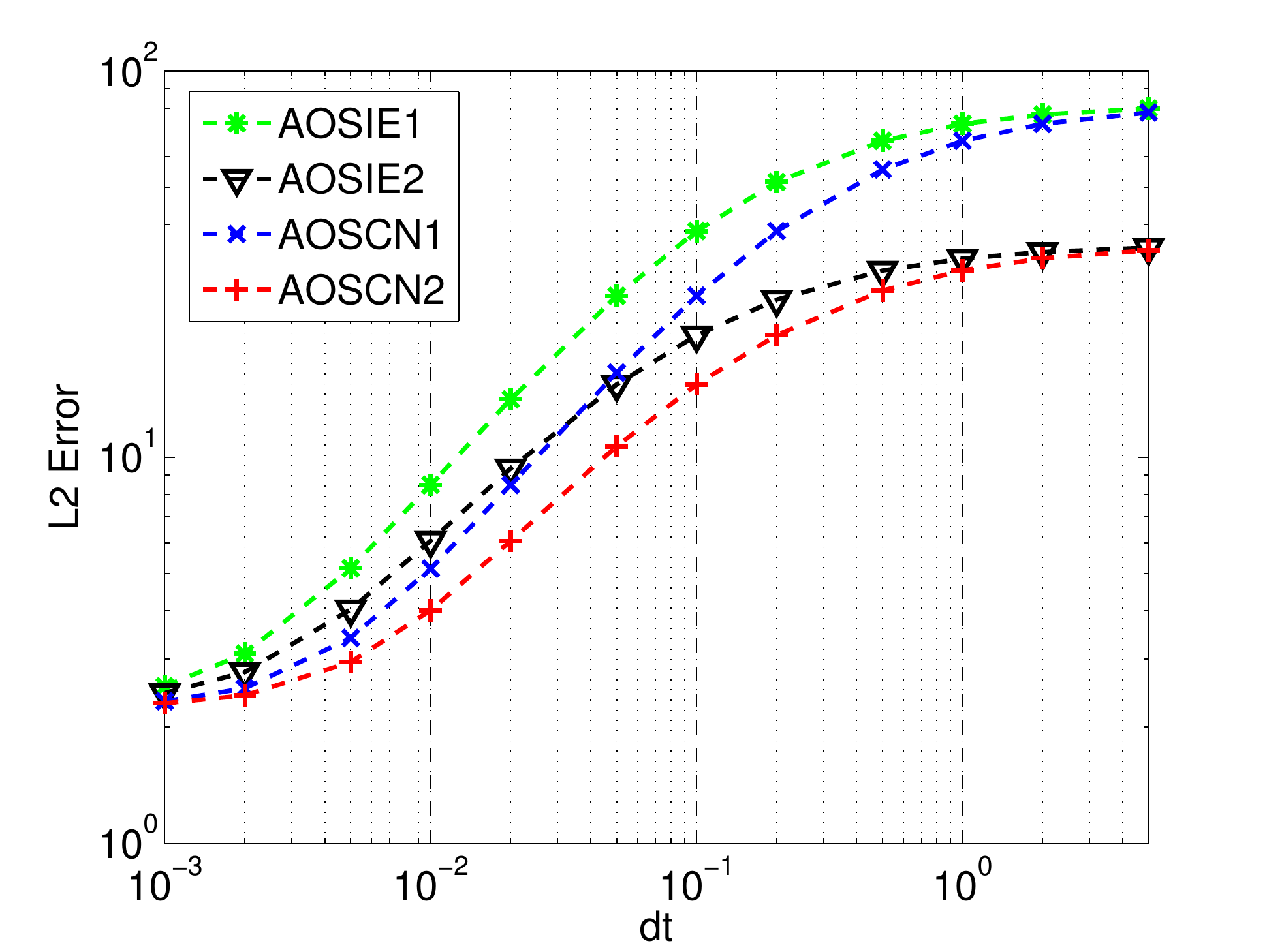}\\
(a) & (b)
\end{tabular}
\caption{Stability verification of the AOS schemes with $H=20$. (a) $h=0.5$; (b) $h=0.25$.}
\label{fig.AOS_AnaStabVer20}
\end{figure}
\begin{figure}[!tb]
\centering
\begin{tabular}{cc}
\includegraphics[width=0.5\linewidth]{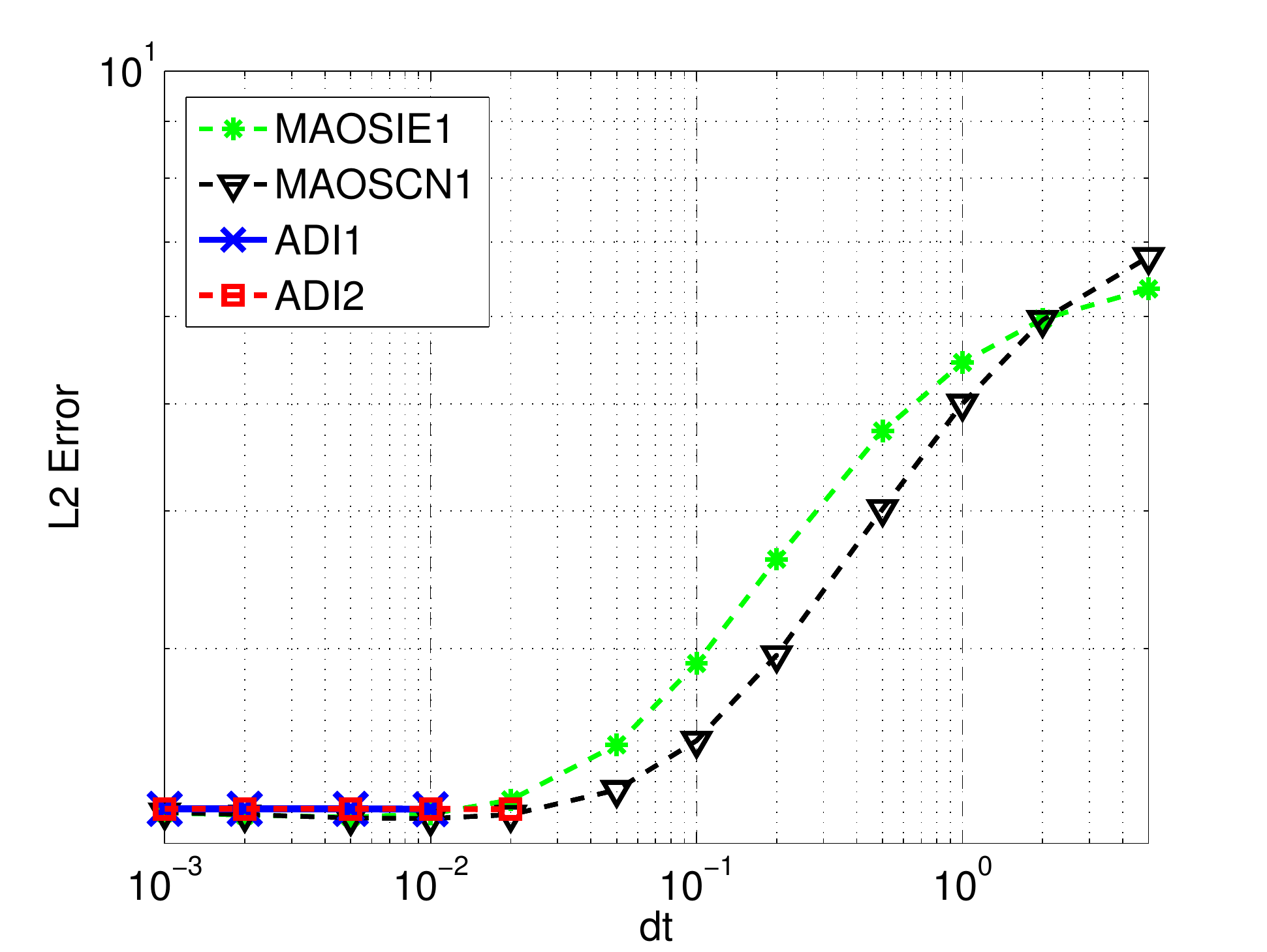} &
\includegraphics[width=0.5\linewidth]{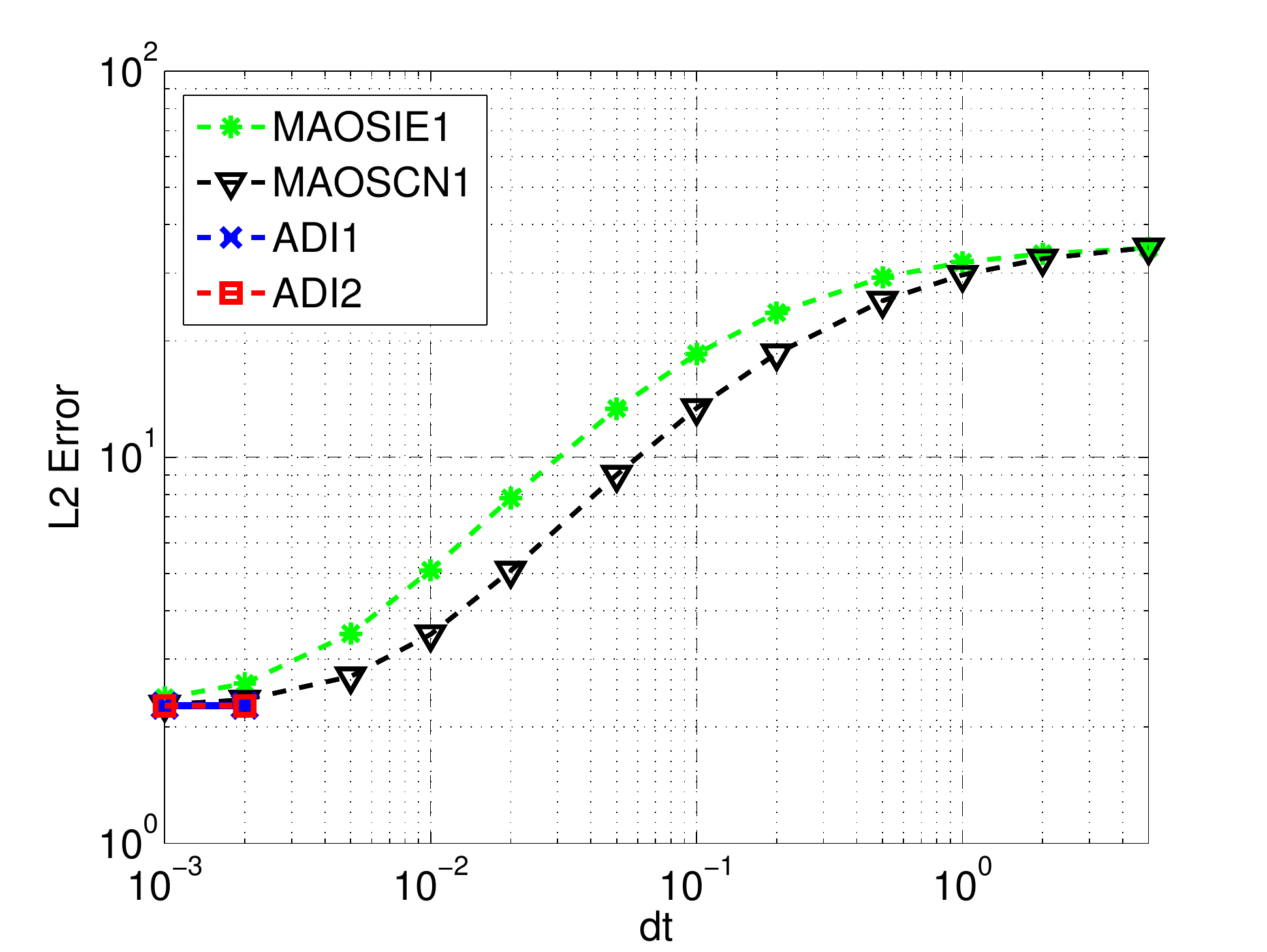}\\
(a) & (b)
\end{tabular}
\caption{Stability verification of the MAOS and ADI schemes with $H=20$. (a) $h=0.5$; (b) $h=0.25$.}
\label{fig.MAOS_ADI_AnaStabVer20}
\end{figure}

In the present study, it is found numerically that the proposed time splitting methods are all unconditionally stable for all tested $H$ values in the range of $H \in [1, 20]$. To illustrate this, we consider the following $\Delta t$ samplings: 
$\Delta t \in \{0.001, 0.002, 0.005, 0.01, 0.02, 0.05, 0.1, 0.2, 0.5, 1, 2, 5\}$. In our tests, the stopping time $T$ is chosen as $T=10^4 \Delta t$ so that enough accumulations are experienced, and the $L_2$ errors will be reported. Of course, when $\Delta t$ is as large as $\Delta t=5$, the numerical error is meaningless. But as long as this error remains to be finite, this demonstrates the stability of the underlying time integration. 

By considering two mesh sizes ($h=0.5$ or $h=0.25$), and two $H$ values ($H=1$ or $H=20$), the numerical results of the proposed LOD, AOS, and MAOS schemes are depicted in Fig. \ref{fig.LOD_AnaStabVer1} through Fig. \ref{fig.MAOS_ADI_AnaStabVer20}. By producing finite error values for all tested $\Delta t$ and $h$ values, all new schemes are demonstrated to be unconditionally stable for the present example. Moreover, such an unconditional stability is not affected by the strong nonlinearity with $H=20$. For a comparison, the results of the ADI schemes \cite{Zhao13} are also shown in Fig. \ref{fig.MAOS_ADI_AnaStabVer1} and Fig. \ref{fig.MAOS_ADI_AnaStabVer20}. The missing points in both figures denote $\Delta t$ values for which the ADI methods were unstable. The present stability results of the ADI schemes are consistent with those reported in \cite{Zhao13}.

\begin{table}[!tb]
\caption{Stability range of all methods for the sampled $\Delta t$ values.}
\label{table.stab}
\begin{center}
\begin{tabular}{lllll}
\hline
 & \multicolumn{2}{c}{$H=1$} & \multicolumn{2}{c}{$H=20$} \\
\cline{2-3}\cline{4-5}
Method & $h=0.5$ & $h=0.25$   & $h=0.5$ & $h=0.25$\\
\hline
LOD     	& $[0.001,5]$ 	& $[0.001,5]$ 	& $[0.001,5]$ 	& $[0.001,5]$\\
AOS   	& $[0.001,5]$ 	& $[0.001,5]$ 	& $[0.001,5]$ 	& $[0.001,5]$\\
MAOS     	& $[0.001,5]$ 	& $[0.001,5]$ 	& $[0.001,5]$ 	& $[0.001,5]$\\
ADI1  	& $[0.001,0.05]$ 	& $[0.001,0.005]$ 	& $[0.001,0.01]$ 	& $[0.001,0.002]$\\
ADI2		& $[0.001,0.02]$	& $[0.001,0.002]$	& $[0.001,0.02]$ 	& $[0.001,0.002]$\\
\hline
\end{tabular}
\end{center}
\end{table}

We summarize the stability results shown in Fig. \ref{fig.LOD_AnaStabVer1} through Fig. \ref{fig.MAOS_ADI_AnaStabVer20} by reporting the stability ranges of all tested methods in Table \ref{table.stab}. For the LOD, AOS, and MAOS schemes, the range is always $[0.001,5]$, because they are stable for all sampled $\Delta t$ values. The upbounds of stability ranges of the ADI schemes are all significantly less than $5$, indicating that the ADI schemes are unstable for large $\Delta$ values. We note that the reported stability range in Table \ref{table.stab} is just a subset of the actual stability interval of each scheme, because such a range is concluded based on selected $\Delta t$ values. For example, all reported schemes are stable for $\Delta t < 0.001$. But the interval $(0,0.001)$ is not included in our ranges. Similarly, the proposed LOD, AOS, and MAOS schemes could be stable for some $\Delta t > 5$. But such large $\Delta t$ values are of little numerical interest.

\begin{table}[!tb]
\caption{Spatial convergence of the LODIE and LODCN schemes in solving the time dependent NPB equation on a sphere.}
\label{table.LODspat}
\begin{center}
\begin{tabular}{lllllllll}
\hline
 & \multicolumn{4}{c}{LODIE1 and LODIE2} & \multicolumn{4}{c}{LODCN1 and LODCN2} \\
\cline{2-5}\cline{6-9}
$h$ & $L_{2}$ & Order   & $L_\infty$ & Order &  $L_{2}$ & Order &  $L_\infty$ & Order\\
\hline
$1$     	& $4.58\text{\sc{e}-}01$ 	& $ 1.18$ 	& $1.73\text{\sc{e}+}00$ 	& $ 3.43$ 	& $9.82\text{\sc{e}-}01$ 	& $2.36$ 	& $2.79\text{\sc{e}+}00$ 	& $4.23$ \\
$0.5$     	& $2.02\text{\sc{e}-}01$ 	& $ 0.84$ 	& $1.61\text{\sc{e}-}01$ 	& $ 0.60$ 	& $1.91\text{\sc{e}-}01$ 	& $0.88$ 	& $1.48\text{\sc{e}-}01$ 	& $0.72$ \\
$0.25$     	& $1.13\text{\sc{e}-}01$ 	& $ 0.92$ 	& $1.06\text{\sc{e}-}01$ 	& $ 0.06$ 	& $1.04\text{\sc{e}-}01$ 	& $1.00$ 	& $8.97\text{\sc{e}-}02$ 	& $0.03$ \\
$0.125$   	& $5.98\text{\sc{e}-}02$ 	& $ 0.63$ 	& $1.02\text{\sc{e}-}01$ 	& $ 0.08$	& $5.17\text{\sc{e}-}02$ 	& $0.65$ 	& $8.80\text{\sc{e}-}02$ 	& $0.07$ \\
$0.0625$  & $3.85\text{\sc{e}-}02$	& 	 	& $9.64\text{\sc{e}-}02$ 	& 		& $3.29\text{\sc{e}-}02$ 	&		& $8.36\text{\sc{e}-}02$ 	&	       \\
\hline
\end{tabular}
\end{center}
\end{table}
\begin{figure}[!tb]
\centering
\begin{tabular}{cc}
\includegraphics[width=0.5\linewidth]{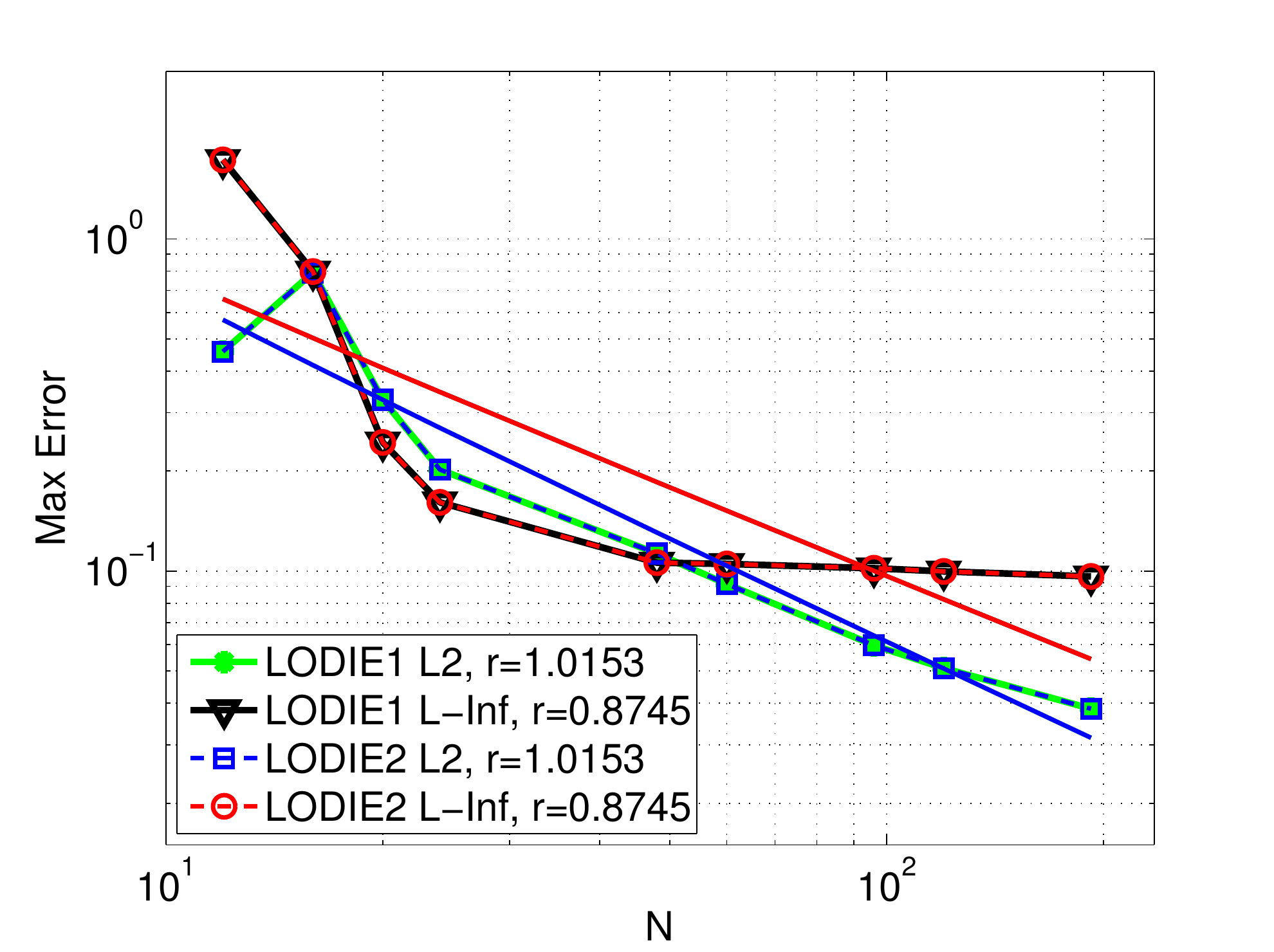} &
\includegraphics[width=0.5\linewidth]{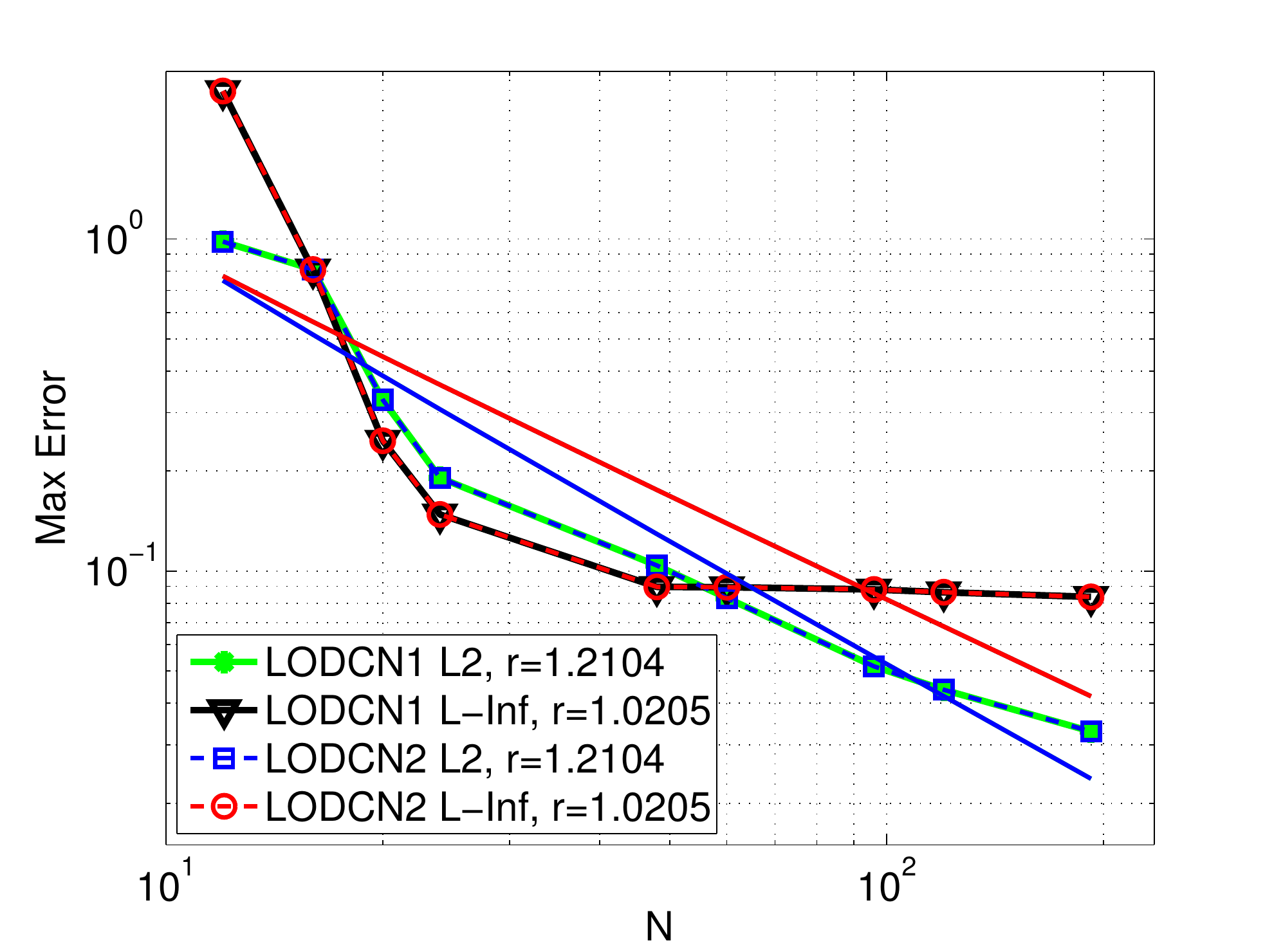} \\
(a) & (b)
\end{tabular}
\caption{Spatial convergence in solving the time dependent NPB equation on a sphere. (a) The LODIE schemes; (b) The LODCN schemes. In both charts, the dashed lines are error plots and the solid lines are linear least square fittings.}
\label{fig.LOD_Spat}
\end{figure}
\begin{table}[!tb]
\caption{Spatial convergence of the AOSIE schemes in solving the time dependent NPB equation on a sphere.}
\label{table.AOSIEspat}
\begin{center}
\begin{tabular}{lllllllll}
\hline
 & \multicolumn{4}{c}{AOSIE1} & \multicolumn{4}{c}{AOSIE2} \\
\cline{2-5}\cline{6-9}
$h$ & $L_{2}$ & Order   & $L_\infty$ & Order &  $L_{2}$ & Order &  $L_\infty$ & Order\\
\hline
$1$     	& $6.34\text{\sc{e}-}01$ 	& $ 0.64$ 	& $1.18\text{\sc{e}+}00$ 	& $ 0.67$ 	& $5.97\text{\sc{e}-}01$ 	& $ 1.09$ 	& $5.80\text{\sc{e}-}01$ 	& $ 1.19$ \\
$0.5$     	& $4.06\text{\sc{e}-}01$ 	& $ 0.66$ 	& $7.37\text{\sc{e}-}01$ 	& $ 0.37$ 	& $2.80\text{\sc{e}-}01$ 	& $ 0.51$ 	& $2.54\text{\sc{e}-}01$ 	& $ 0.26$ \\
$0.25$     	& $2.58\text{\sc{e}-}01$ 	& $ 0.71$ 	& $5.71\text{\sc{e}-}01$ 	& $ 0.18$	& $1.96\text{\sc{e}-}01$ 	& $ 0.71$ 	& $2.11\text{\sc{e}-}01$ 	& $ 0.13$ \\
$0.125$   	& $1.58\text{\sc{e}-}01$ 	& $ 0.59$ 	& $5.03\text{\sc{e}-}01$ 	& $ 0.09$	& $1.20\text{\sc{e}-}01$ 	& $ 0.58$ 	& $1.93\text{\sc{e}-}01$ 	& $ 0.09$ \\
$0.0625$  & $1.05\text{\sc{e}-}01$	& 	 	& $4.71\text{\sc{e}-}01$ 	& 		& $8.00\text{\sc{e}-}02$ 	&		& $1.82\text{\sc{e}-}01$ 	&	       \\
\hline
\end{tabular}
\end{center}
\end{table}
\begin{table}[!tb]
\caption{Spatial convergence of the AOSCN schemes in solving the time dependent NPB equation on a sphere.}
\label{table.AOSCNspat}
\begin{center}
\begin{tabular}{lllllllll}
\hline
 & \multicolumn{4}{c}{AOSCN1} & \multicolumn{4}{c}{AOSCN2} \\
\cline{2-5}\cline{6-9}
$h$ & $L_{2}$ & Order   & $L_\infty$ & Order &  $L_{2}$ & Order &  $L_\infty$ & Order\\
\hline
$1$     	& $3.84\text{\sc{e}-}01$ 	& $ 0.42$ 	& $1.04\text{\sc{e}+}00$ 	& $ 1.12$ 	& $3.65\text{\sc{e}-}01$ 	& $ 0.62$ 	& $6.18\text{\sc{e}-}01$ 	& $ 1.55$ \\
$0.5$     	& $2.88\text{\sc{e}-}01$ 	& $ 0.75$ 	& $4.76\text{\sc{e}-}01$ 	& $ 0.46$ 	& $2.38\text{\sc{e}-}01$ 	& $ 0.70$ 	& $2.12\text{\sc{e}-}01$ 	& $ 0.45$ \\
$0.25$     	& $1.71\text{\sc{e}-}01$ 	& $ 0.76$ 	& $3.46\text{\sc{e}-}01$ 	& $ 0.17$	& $1.47\text{\sc{e}-}01$ 	& $ 0.79$ 	& $1.55\text{\sc{e}-}01$ 	& $ 0.11$ \\
$0.125$   	& $1.01\text{\sc{e}-}01$ 	& $ 0.61$ 	& $3.08\text{\sc{e}-}01$ 	& $ 0.09$	& $8.52\text{\sc{e}-}02$ 	& $ 0.61$ 	& $1.44\text{\sc{e}-}01$ 	& $ 0.08$ \\
$0.0625$  & $6.61\text{\sc{e}-}02$	& 	 	& $2.89\text{\sc{e}-}01$ 	& 		& $5.57\text{\sc{e}-}02$ 	&		& $1.36\text{\sc{e}-}01$ 	&	       \\
\hline
\end{tabular}
\end{center}
\end{table}
\begin{figure}[!tb]
\centering
\begin{tabular}{cc}
\includegraphics[width=0.5\linewidth]{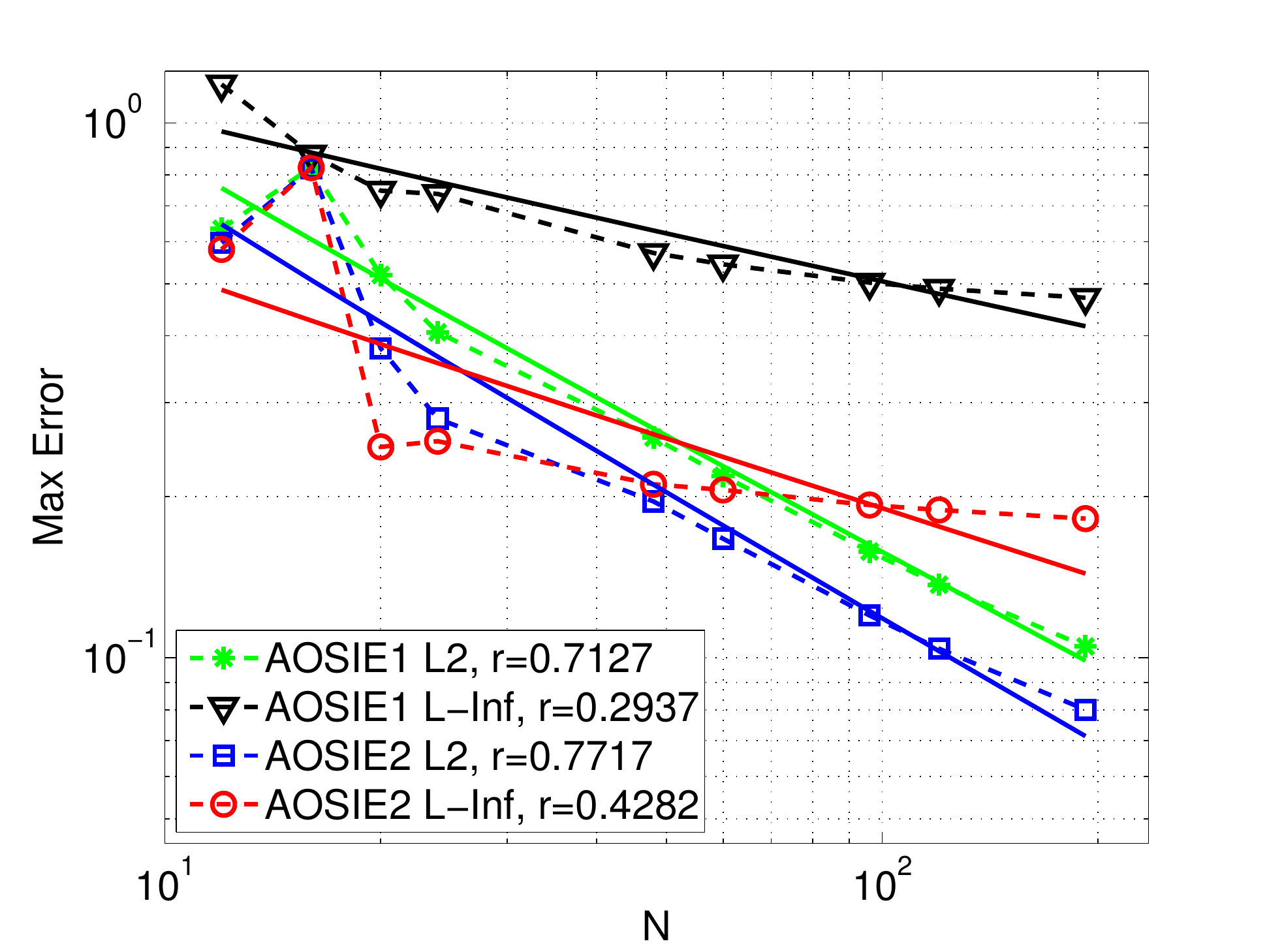} &
\includegraphics[width=0.5\linewidth]{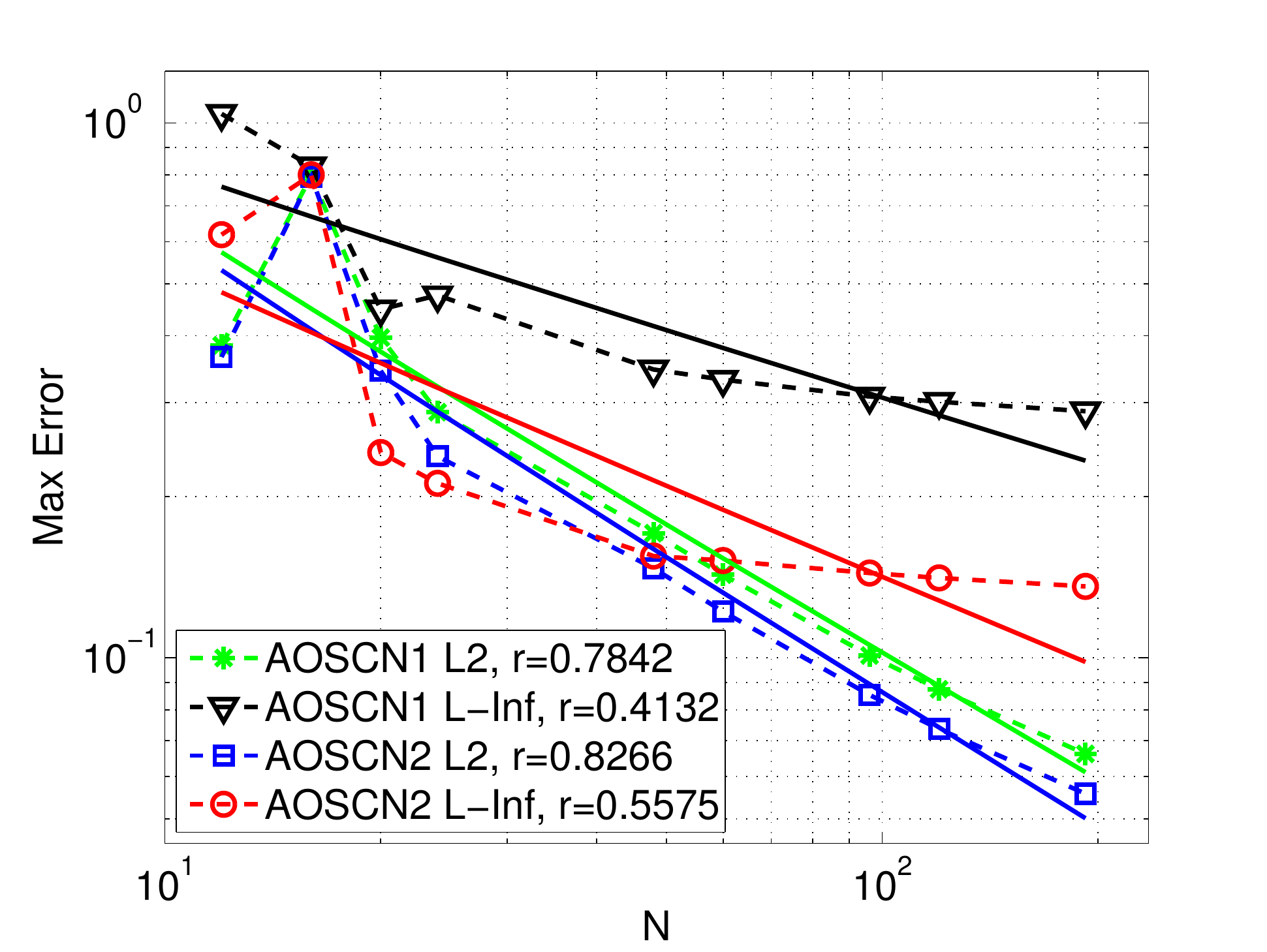} \\
(a) & (b)
\end{tabular}
\caption{Spatial convergence in solving the time dependent NPB equation on a sphere. (a) The AOSIE schemes; (b) The AOSCN schemes. In both charts, the dashed lines are error plots and the solid lines are linear least square fittings.}
\label{fig.AOS_Spat}
\end{figure}
\begin{table}[!tb]
\caption{Spatial convergence of the MAOS schemes in solving the time dependent NPB equation on a sphere.}
\label{table.MAOSspat}
\begin{center}
\begin{tabular}{lllllllll}
\hline
 & \multicolumn{4}{c}{MAOSIE1} & \multicolumn{4}{c}{MAOSCN1} \\
\cline{2-5}\cline{6-9}
$h$ & $L_{2}$ & Order   & $L_\infty$ & Order &  $L_{2}$ & Order &  $L_\infty$ & Order\\
\hline
$1$     	& $5.53\text{\sc{e}-}01$ 	& $ 1.08$ 	& $5.33\text{\sc{e}-}01$ 	& $ 1.19$ 	& $3.41\text{\sc{e}-}01$ 	& $ 0.58$ 	& $5.90\text{\sc{e}-}01$ 	& $ 1.58$ \\
$0.5$     	& $2.61\text{\sc{e}-}01$ 	& $ 0.60$ 	& $2.33\text{\sc{e}-}01$ 	& $ 0.35$ 	& $2.27\text{\sc{e}-}01$ 	& $ 0.75$ 	& $1.97\text{\sc{e}-}01$ 	& $ 0.49$ \\
$0.25$     	& $1.73\text{\sc{e}-}01$ 	& $ 0.74$ 	& $1.84\text{\sc{e}-}01$ 	& $ 0.12$	& $1.36\text{\sc{e}-}01$ 	& $ 0.82$ 	& $1.40\text{\sc{e}-}01$ 	& $ 0.10$ \\
$0.125$   	& $1.04\text{\sc{e}-}01$ 	& $ 0.58$ 	& $1.69\text{\sc{e}-}01$ 	& $ 0.09$	& $7.70\text{\sc{e}-}02$ 	& $ 0.62$ 	& $1.31\text{\sc{e}-}01$ 	& $ 0.08$ \\
$0.0625$  & $6.91\text{\sc{e}-}02$	& 	 	& $1.59\text{\sc{e}-}01$ 	& 		& $5.02\text{\sc{e}-}02$ 	&		& $1.24\text{\sc{e}-}01$ 	&	       \\
\hline
\end{tabular}
\end{center}
\end{table}
\begin{figure}[!tb]
\centering
\begin{tabular}{cc}
\includegraphics[width=0.5\linewidth]{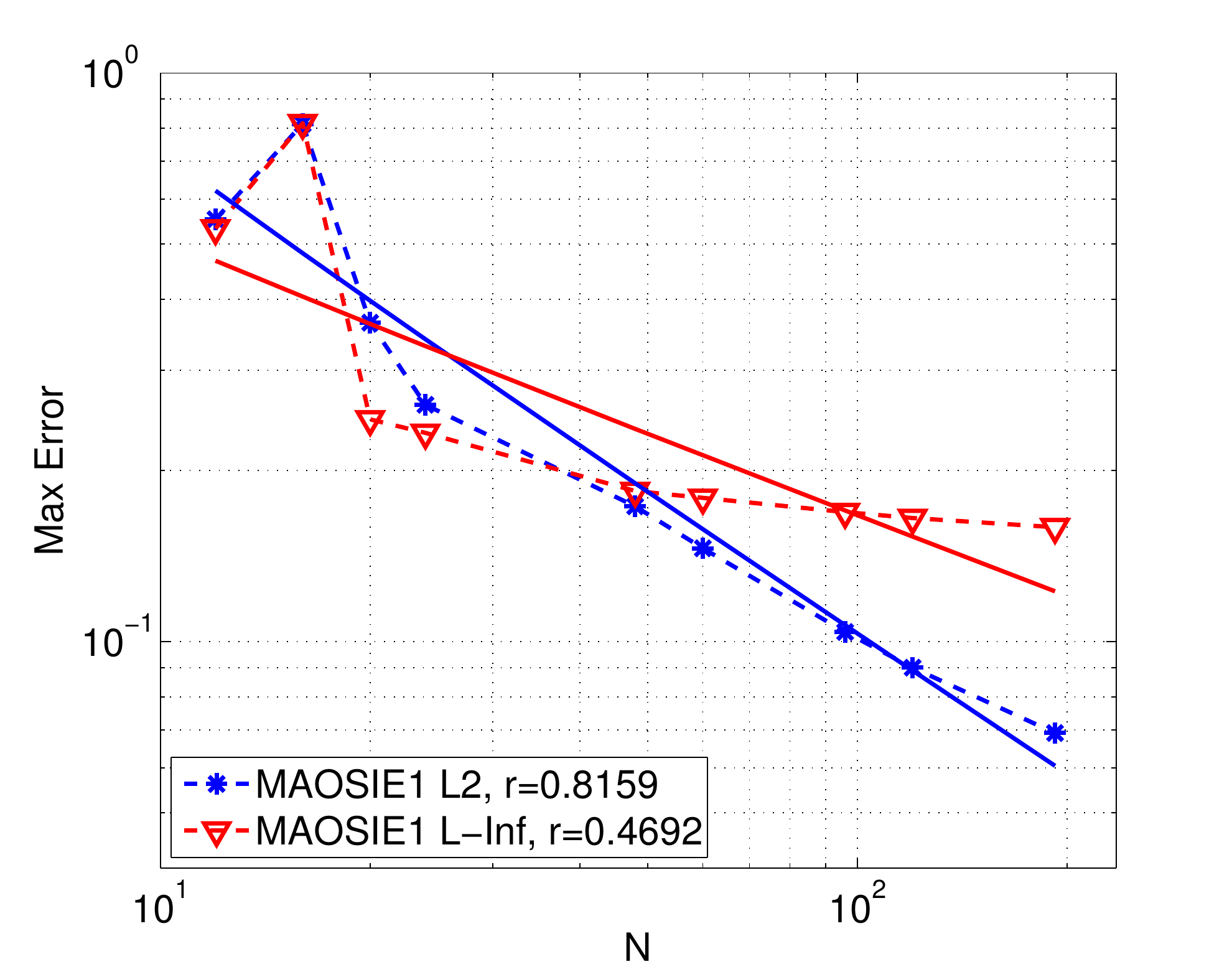} &
\includegraphics[width=0.5\linewidth]{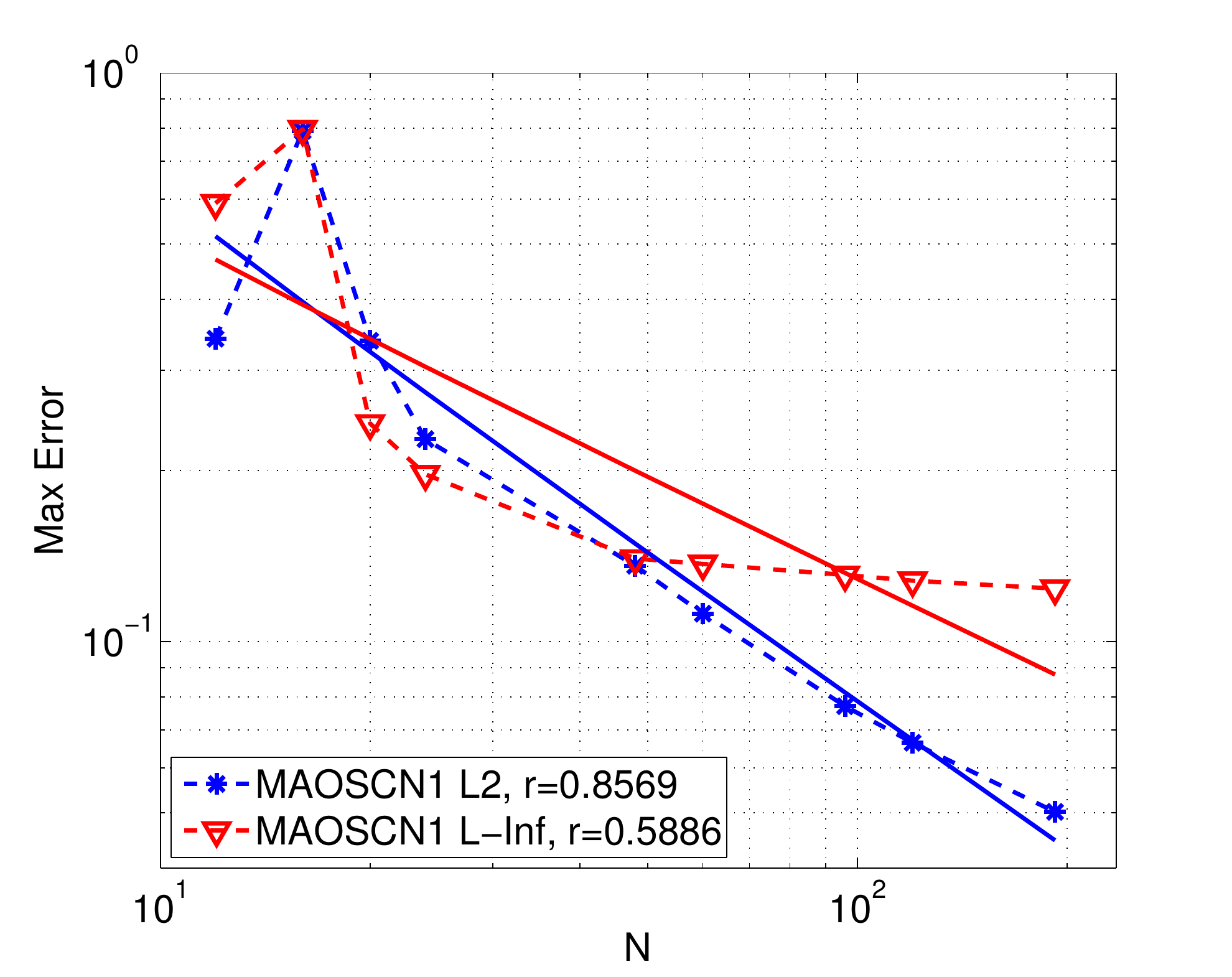} \\
(a) & (b)
\end{tabular}
\caption{Spatial convergence in solving the time dependent NPB equation on a sphere. (a) The MAOSIE schemes; (b) The MAOSCN schemes. In both charts, the dashed lines are error plots and the solid lines are linear least square fittings.}
\label{fig.MAOS_Spat}
\end{figure}

\subsection{Spatial Convergence}
We next study the spatial convergence of the proposed time splitting schemes. For this purpose, we vary the spatial element size $h$ from $1$ to $0.0625$, and take the time increment to be sufficiently small in all cases, i.e., $\Delta t = \frac{h^2}{20}$. We take the magnitude $H=1$ in the initial solution and the stopping time to be $T=10$ for the rest of the studies in this section.  The $L_2$ and $L_\infty$ errors of the LOD, AOS, and MAOS schemes are listed in Tables \ref{table.LODspat} through \ref{table.MAOSspat}, and are depicted in Figures \ref{fig.LOD_Spat} through \ref{fig.MAOS_Spat} as well.

By considering two consecutive mesh refinements, the numerical orders can be calculated for each refinement and are reported in Tables \ref{table.LODspat} through \ref{table.MAOSspat}. A similar pattern can be seen for all tested schemes. In particular, the numerical orders become smaller when $h$ decreases. The orders in $L_2$ norms are higher than those in $L_\infty$ norms and the $L_2$ errors are also smaller than the $L_\infty$ errors. To calculate an overall order of accuracy, the errors based on more $h$ values are shown in Figures \ref{fig.LOD_Spat} through \ref{fig.MAOS_Spat} in logarithmic scales. Here the mesh size $N$ is used, which is inversely proportional to $h$ by $N=6/h$. In all plots, the solid lines represent the linear least-squares fittings of the errors. The slope $r$ of these lines reflects the average convergence rate. It can be seen that $L_\infty$ orders are all below one. As pointed out in \cite{Zhao13}, such accuracy reduction is due to the singular solution at the center of the sphere and the non-smoothness of the solution across the interface $\Gamma$, which is treated in an approximate sense in the present finite difference discretization. The non-uniform convergence in the $L_\infty$ norm is particularly affected by the charge singularity at the center, because large errors exist in the calculation of potential near the center. However, the overall $L_2$ orders of all schemes are found to be within $0.7$ to $1.2$. This suggests that these schemes display approximately first-order convergence in the $L_2$ norm.

We are more interested in a comparison among these schemes. Interestingly, both LODIE schemes produce the exact same errors for this simple spherical case, as do the LODCN schemes. For other schemes, some minor changes will affect the orders slightly. In general, the orders of the AOS and MAOS schemes appear significantly worse than those exhibited by the LOD schemes.

\begin{table}[!tb]
\caption{Temporal convergence of the LODIE and LODCN schemes in solving the time dependent NPB equation on a sphere.}
\label{table.LODtemp}
\begin{center}
\begin{tabular}{lllllllll}
\hline
 & \multicolumn{4}{c}{LODIE1 and LODIE2} & \multicolumn{4}{c}{LODCN1 and LODCN2} \\
\cline{2-5}\cline{6-9}
$dt$ & $L_{2}$ & Order   & $L_\infty$ & Order &  $L_{2}$ & Order &  $L_\infty$ & Order\\
\hline
$8.00\text{\sc{e}-}04$     	& $8.04\text{\sc{e}-}02$ 	& $1.06$ 	& $1.55\text{\sc{e}-}01$ 	& $ 1.06$ 	& $5.46\text{\sc{e}-}02$ 	& $1.07$ 	& $1.06\text{\sc{e}-}01$ 	& $1.07$ \\
$4.00\text{\sc{e}-}04$     	& $3.85\text{\sc{e}-}02$ 	& $1.11$ 	& $7.44\text{\sc{e}-}02$ 	& $ 1.11$ 	& $2.59\text{\sc{e}-}02$ 	& $1.11$ 	& $5.03\text{\sc{e}-}02$ 	& $1.11$ \\
$2.00\text{\sc{e}-}04$  	& $1.79\text{\sc{e}-}02$ 	& $1.23$ 	& $3.45\text{\sc{e}-}02$ 	& $ 1.23$	& $1.20\text{\sc{e}-}02$ 	& $1.23$ 	& $2.32\text{\sc{e}-}02$ 	& $1.23$ \\
$1.00\text{\sc{e}-}04$	& $7.65\text{\sc{e}-}03$ 	& $1.59$ 	& $1.48\text{\sc{e}-}02$ 	& $ 1.59$	& $5.11\text{\sc{e}-}03$ 	& $1.59$ 	& $9.91\text{\sc{e}-}03$ 	& $1.59$ \\
$5.00\text{\sc{e}-}05$	& $2.55\text{\sc{e}-}03$    & 	 	& $4.91\text{\sc{e}-}03$ 	& 	 	& $1.70\text{\sc{e}-}03$ 	&	 	& $3.30\text{\sc{e}-}03$ 	&	\\
$2.50\text{\sc{e}-}05$	& $0.00\text{\sc{e}+}00$   & 	 	& $0.00\text{\sc{e}+}00$ 	& 	 	& $0.00\text{\sc{e}+}00$ 	&	 	& $0.00\text{\sc{e}+}00$ 	&	\\
\hline
\end{tabular}
\end{center}
\end{table}
\begin{table}[!tb]
\caption{Temporal convergence of the AOSIE schemes in solving the time dependent NPB equation on a sphere.}
\label{table.AOSIEtemp}
\begin{center}
\begin{tabular}{lllllllll}
\hline
 & \multicolumn{4}{c}{AOSIE1} & \multicolumn{4}{c}{AOSIE2} \\
\cline{2-5}\cline{6-9}
$dt$ & $L_{2}$ & Order   & $L_\infty$ & Order &  $L_{2}$ & Order &  $L_\infty$ & Order\\
\hline
$8.00\text{\sc{e}-}04$     	& $1.20\text{\sc{e}-}01$ 	& $0.98$ 	& $1.16\text{\sc{e}-}01$ 	& $ 0.97$ 	& $8.36\text{\sc{e}-}02$ 	& $0.99$ 	& $1.16\text{\sc{e}-}01$ 	& $0.97$ \\
$4.00\text{\sc{e}-}04$     	& $6.09\text{\sc{e}-}02$ 	& $1.06$ 	& $5.94\text{\sc{e}-}02$ 	& $ 1.06$ 	& $4.22\text{\sc{e}-}02$ 	& $1.06$ 	& $5.94\text{\sc{e}-}02$ 	& $1.06$ \\
$2.00\text{\sc{e}-}04$  	& $2.92\text{\sc{e}-}02$ 	& $1.20$ 	& $2.86\text{\sc{e}-}02$ 	& $ 1.20$	& $2.02\text{\sc{e}-}02$ 	& $1.20$ 	& $2.86\text{\sc{e}-}02$ 	& $1.20$ \\
$1.00\text{\sc{e}-}04$	& $1.27\text{\sc{e}-}02$ 	& $1.57$ 	& $1.25\text{\sc{e}-}02$ 	& $ 1.57$	& $8.77\text{\sc{e}-}03$ 	& $1.57$ 	& $1.25\text{\sc{e}-}02$ 	& $1.57$ \\
$5.00\text{\sc{e}-}05$	& $4.27\text{\sc{e}-}03$    & 	 	& $4.19\text{\sc{e}-}03$ 	& 	 	& $2.95\text{\sc{e}-}03$ 	&	 	& $4.19\text{\sc{e}-}03$ 	&	\\
$2.50\text{\sc{e}-}05$	& $0.00\text{\sc{e}+}00$   & 	 	& $0.00\text{\sc{e}+}00$ 	& 	 	& $0.00\text{\sc{e}+}00$ 	&	 	& $0.00\text{\sc{e}+}00$ 	&	\\
\hline
\end{tabular}
\end{center}
\end{table}
\begin{table}[!tb]
\caption{Temporal convergence of the AOSCN schemes in solving the time dependent NPB equation on a sphere.}
\label{table.AOSCNtemp}
\begin{center}
\begin{tabular}{lllllllll}
\hline
 & \multicolumn{4}{c}{AOSCN1} & \multicolumn{4}{c}{AOSCN2} \\
\cline{2-5}\cline{6-9}
$dt$ & $L_{2}$ & Order   & $L_\infty$ & Order &  $L_{2}$ & Order &  $L_\infty$ & Order\\
\hline
$8.00\text{\sc{e}-}04$     	& $6.31\text{\sc{e}-}02$ 	& $1.01$ 	& $6.16\text{\sc{e}-}02$ 	& $ 1.00$ 	& $4.36\text{\sc{e}-}02$ 	& $1.01$ 	& $6.16\text{\sc{e}-}02$ 	& $1.00$ \\
$4.00\text{\sc{e}-}04$     	& $3.13\text{\sc{e}-}02$ 	& $1.08$ 	& $3.07\text{\sc{e}-}02$ 	& $ 1.08$ 	& $2.16\text{\sc{e}-}02$ 	& $1.08$ 	& $3.07\text{\sc{e}-}02$ 	& $1.08$ \\
$2.00\text{\sc{e}-}04$  	& $1.49\text{\sc{e}-}02$ 	& $1.21$ 	& $1.46\text{\sc{e}-}02$ 	& $ 1.21$	& $1.02\text{\sc{e}-}02$ 	& $1.21$ 	& $1.46\text{\sc{e}-}02$ 	& $1.21$ \\
$1.00\text{\sc{e}-}04$	& $6.42\text{\sc{e}-}03$ 	& $1.58$ 	& $6.31\text{\sc{e}-}03$ 	& $ 1.58$	& $4.43\text{\sc{e}-}03$ 	& $1.58$ 	& $6.31\text{\sc{e}-}03$ 	& $1.58$ \\
$5.00\text{\sc{e}-}05$	& $2.15\text{\sc{e}-}03$    & 	 	& $2.11\text{\sc{e}-}03$ 	& 	 	& $1.48\text{\sc{e}-}03$ 	&	 	& $2.11\text{\sc{e}-}03$ 	&	\\
$2.50\text{\sc{e}-}05$	& $0.00\text{\sc{e}+}00$   & 	 	& $0.00\text{\sc{e}+}00$ 	& 	 	& $0.00\text{\sc{e}+}00$ 	&	 	& $0.00\text{\sc{e}+}00$ 	&	\\
\hline
\end{tabular}
\end{center}
\end{table}
\begin{table}[!tb]
\caption{Temporal convergence of the MAOS schemes in solving the time dependent NPB equation on a sphere.}
\label{table.MAOStemp}
\begin{center}
\begin{tabular}{lllllllll}
\hline
 & \multicolumn{4}{c}{MAOSIE1} & \multicolumn{4}{c}{MAOSCN1} \\
\cline{2-5}\cline{6-9}
$dt$ & $L_{2}$ & Order   & $L_\infty$ & Order &  $L_{2}$ & Order &  $L_\infty$ & Order\\
\hline
$8.00\text{\sc{e}-}04$     	& $6.40\text{\sc{e}-}02$ 	& $1.00$ 	& $8.96\text{\sc{e}-}02$ 	& $ 0.99$ 	& $3.31\text{\sc{e}-}02$ 	& $1.02$ 	& $4.69\text{\sc{e}-}02$ 	& $1.01$ \\
$4.00\text{\sc{e}-}04$     	& $3.20\text{\sc{e}-}02$ 	& $1.07$ 	& $4.53\text{\sc{e}-}02$ 	& $ 1.06$ 	& $1.64\text{\sc{e}-}02$ 	& $1.08$ 	& $2.33\text{\sc{e}-}02$ 	& $1.08$ \\
$2.00\text{\sc{e}-}04$  	& $1.53\text{\sc{e}-}02$ 	& $1.21$ 	& $2.16\text{\sc{e}-}02$ 	& $ 1.20$	& $7.73\text{\sc{e}-}03$ 	& $1.21$ 	& $1.10\text{\sc{e}-}02$ 	& $1.21$ \\
$1.00\text{\sc{e}-}04$	& $6.61\text{\sc{e}-}03$ 	& $1.43$ 	& $9.40\text{\sc{e}-}03$ 	& $ 1.43$	& $3.33\text{\sc{e}-}03$ 	& $1.43$ 	& $4.75\text{\sc{e}-}03$ 	& $1.43$ \\
$5.00\text{\sc{e}-}05$	& $2.45\text{\sc{e}-}03$    & 	 	& $3.48\text{\sc{e}-}03$ 	& 	 	& $1.24\text{\sc{e}-}03$ 	&	 	& $1.76\text{\sc{e}-}03$ 	&	\\
$2.50\text{\sc{e}-}05$	& $0.00\text{\sc{e}+}00$   & 	 	& $0.00\text{\sc{e}+}00$ 	& 	 	& $0.00\text{\sc{e}+}00$ 	&	 	& $0.00\text{\sc{e}+}00$ 	&	\\
\hline
\end{tabular}
\end{center}
\end{table}

\subsection{Temporal Convergence}
We finally investigate the temporal convergence of our schemes, displayed in Tables \ref{table.LODtemp} through \ref{table.MAOStemp}. Because we should expect the spatial discretization error to be much larger than the temporal discretization error, due to the charge singularity and non-smooth interface previously mentioned, we may not be able to determine temporal order of convergence by comparing solutions calculated at various time increments $\Delta t$ to the analytical solution. 

Thus, we choose a fixed spacing $h=0.125$, and use the potentials calculated at the time increment $\Delta t = 2.5\text{\sc{e}-}05$ as the reference solution. We then proceed to calculate the $L_2$ and $L_\infty$ norms, and their respective orders of convergence, in the same manner as the previous subsection. The temporal order of convergence for all schemes in $L_2$ and $L_\infty$ norms appears to be slightly higher than first-order, and the orders for the MAOS schemes in Table \ref{table.MAOStemp} appear to be the lowest. We also note that two versions of the LODIE and LODCN schemes yield exactly the same results. 

\section{Biological applications}\label{sec.application}
In this section, we further explore the stability and accuracy of the LOD, AOS, and MAOS schemes by considering the solvation analysis of real proteins. A detailed comparison among different time splitting schemes will be conducted for a particular protein system. Some corrections will be introduced to improve the accuracy of the time splitting. After identifying an optimum NPB solver, its usage on various protein systems will be considered. The CPU acceleration with respect to the existing operator splitting schemes will be examined. 

\subsection{Solvation energy and numerical setup}
To quantitatively verify and compare the accuracy and efficiency of our new schemes, we compute solvation free energy based on the electrostatic potential calculated from the NPB equation. The energy released when the solute macromolecule is dissolved in solvent is known as the free energy of solvation, or solvation free energy. This solvation free energy may be calculated in our setting by computing the difference between total free energy of the macromolecule in a vacuum and in the solvent. Because we consider here only electrostatic effects, we may define the solvation energy as
\begin{equation}
\Delta G = G_s - G_0 = \frac{1}{2} \int_{\Omega} \! \rho_m (\phi_m ({\bf r}) - \phi_0 ({\bf r})) \ \textnormal{d}{\bf r}
\end{equation}
where $\phi_m$ is the electrostatic potential, in units of $kcal/mol/e_c$, in solvent, and $\phi_0$ is the electrostatic potential in vacuum. These potential values are obtained by scaling our calculated dimensionless potentials with the constant $0.592183$, corresponding to room temperature $(298K)$. In a discrete setting, we may calculate the solvation energy as
\begin{equation}\label{deltaG}
\Delta G = \frac{1}{2} \sum\limits_i \sum\limits_j \sum\limits_k Q(x_i,y_j,z_k)(\phi_m (x_i,y_j,z_k) - \phi_0 (x_i,y_j,z_k))
\end{equation}
where $Q$ is the trilinear interpolation of the singular charges in $\rho_m$. 
The scaled potential $\phi_m$ is calculated based on the time dependent NPB equation, while $\phi_0$ is computed by solving a simple Poisson equation with a constant dielectric coefficient. As in the literature \cite{Lu08,Li13}, the Fast Fourier Transform (FFT) can be used to efficiently determine $\phi_0$. 

In all cases, a uniform mesh size $h=0.5$ is used along all three dimensions and a large enough computational domain is chosen. In our calculations, we set the dielectric constants $\epsilon_s=80$, and $\epsilon_m=1$. The ionic strength $I$ is set to 9.48955M, so the nonlinear constant $\bar{\kappa}$ is 0.1261. The scaling parameter is chosen as $\alpha=1$, unless specified otherwise. A sufficiently large stopping time $T$ will be used in each case to ensure that the steady state solution is reached. The molecular surface $\Gamma$ underlying our computations is calculated based on the MSMS package \cite{Sanner}. In the MSMS surface generation, the probe radius is set to 1.4, and the density is chosen as 10. A Lagrangian to Eulerian conversion \cite{Zhou08} is conducted to convert the triangular surface mesh of the MSMS package into 3D Cartesian grid values for determining $\epsilon(x_i,y_j,z_k)$.

\begin{figure}[!tb]
\centering
\begin{tabular}{cc}
\includegraphics[width=0.5\linewidth]{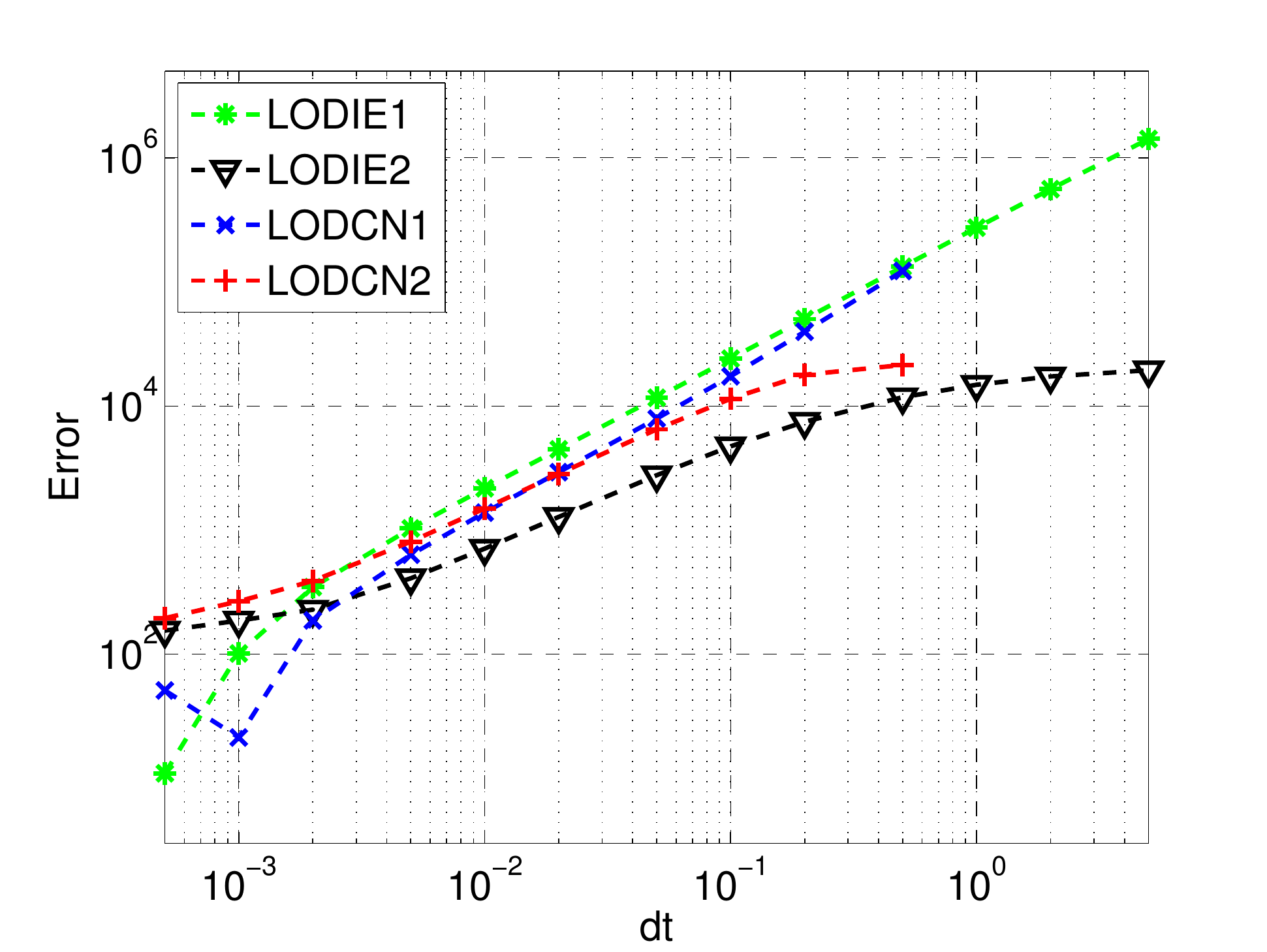} &
\includegraphics[width=0.5\linewidth]{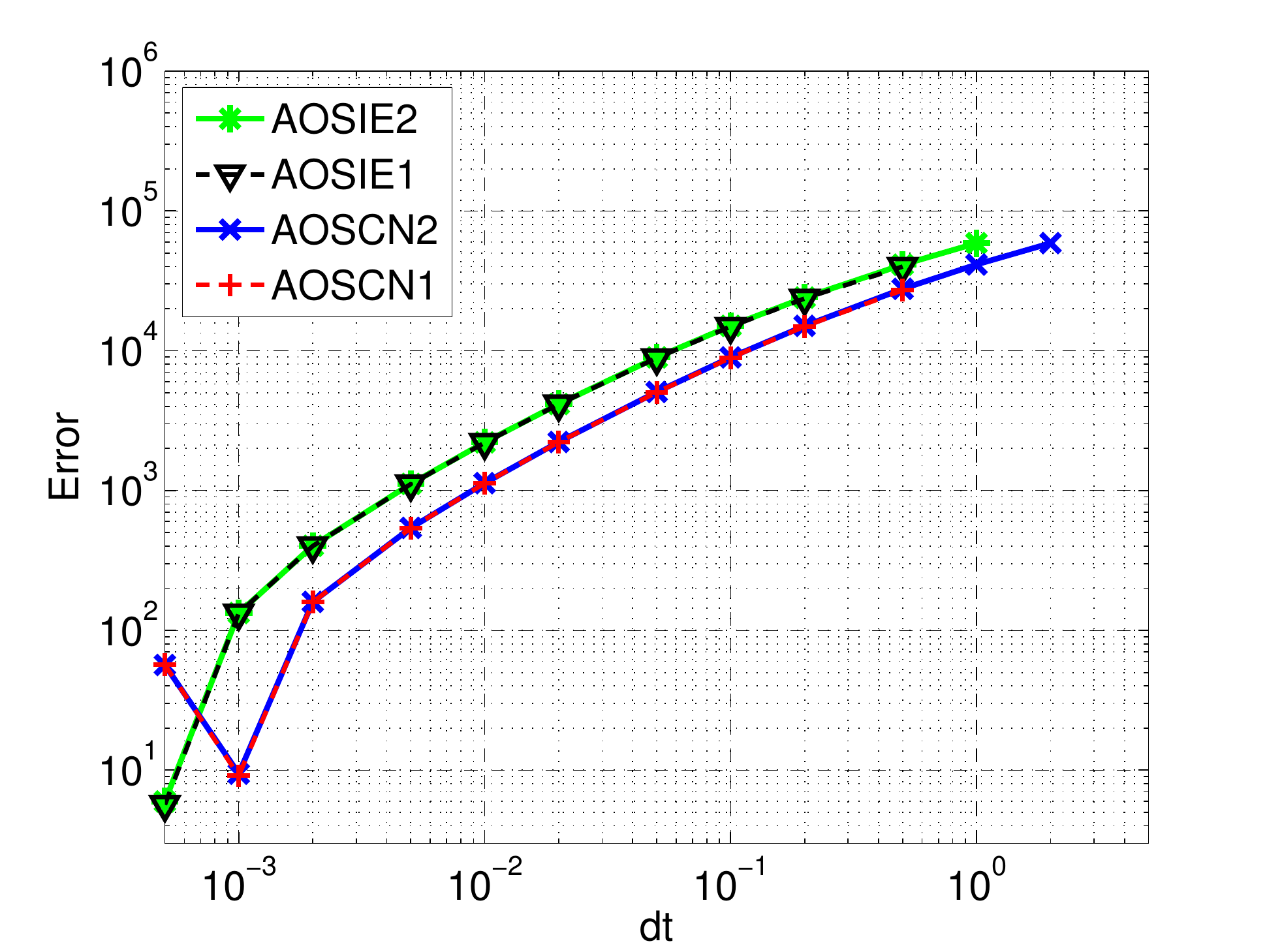} \\
(a) & (b)
\end{tabular}
\caption{Stability verification of the LOD schemes (a) and the AOS schemes (b) for protein 1ajj. }
\label{fig.LOD_AOS_StabVer1ajj}
\end{figure}
\begin{figure}[!tb]
\centering
\begin{tabular}{cc}
\includegraphics[width=0.5\linewidth]{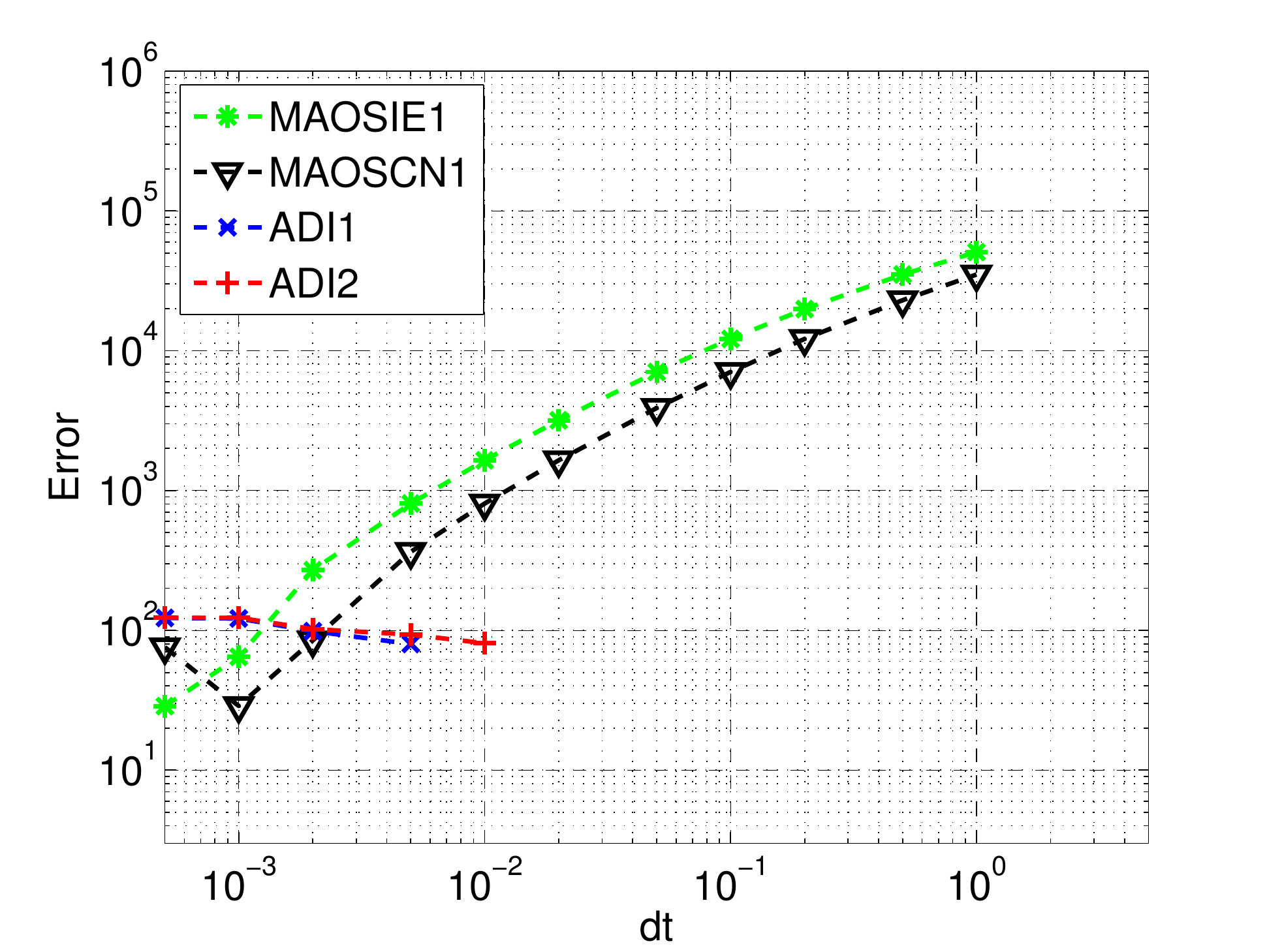}
\end{tabular}
\caption{Stability verification of the MAOS and ADI schemes for protein 1ajj.}
\label{fig.MAOS_ADI_StabVer1ajj}
\end{figure}
\begin{table}[!tb]
\caption{Stability range of all methods on 1ajj for the sampled $\Delta t$ values.}
\label{table.1ajjstabver}
\begin{center}
\begin{tabular}{lllll}
\cline{1-2}\cline{4-5}
Method & Stability Range && Method & Stability Range\\
\cline{1-2}\cline{4-5}
LODIE1     	& $[0.0005,5]$&&AOSCN1     	& $[0.0005,0.5]$\\
LODIE2     	& $[0.0005,5]$&&AOSCN2	   	& $[0.0005,2]$\\
LODCN1     	& $[0.0005,0.5]$&&MAOSIE1     	& $[0.0005,1]$\\
LODCN2     	& $[0.0005,0.5]$&&MAOSCN1     	& $[0.0005,1]$\\
AOSIE1	   	& $[0.0005,0.5]$&&ADI1  			& $[0.0005,0.005]$\\
AOSIE2	   	& $[0.0005,1]$&&ADI2			& $[0.0005,0.01]$\\
\cline{1-2}\cline{4-5}
\end{tabular}
\end{center}
\end{table}

\subsection{Stability analysis}
In the previous section, all of the proposed time splitting schemes are found to be unconditionally stable for a single atom case. It is of great interest to see if these schemes are still stable for real protein systems. We carry out this stability study by considering a protein 1ajj, a low-density lipoprotein receptor. As in the previous studies, we set $T=10^4 \Delta t $ and consider a set of sampled $\Delta t$ values, $\Delta t \in \{0.0005, 0.001, 0.002, 0.005, 0.01, 0.02, 0.05, 0.1, 0.2, 0.5, 1, 2, 5\}$. For a reference, an Euler solution is produced from the explicit Euler method with a very small time increment $\Delta t=1\text{\sc{e}-}05$ and a large enough $T=10$. The solvation energy of the Euler solution is calculated to be $\Delta G=-1209.7$, which will be taken as the reference energy value in studying pseudo-time NPB solvers.

The stability results in terms of the solvation energy errors of the proposed time splitting schemes are depicted in Figs. \ref{fig.LOD_AOS_StabVer1ajj} and \ref{fig.MAOS_ADI_StabVer1ajj}. We note that for the LODIE1 and LODIE2 methods, the finite error at each $\Delta t$ implies stability of these two methods at every tested time step value. In fact, the LODIE methods are also found to be unconditionally stable for solving other protein systems, while all other time splitting schemes are conditionally stable for real proteins. In particular, as shown in Table \ref{table.1ajjstabver}, the LODCN, AOSIE1, and AOSCN1 methods were only unstable for $\Delta t=1, 2, 5.$ The AOSIE2 and MAOS methods were only unstable for $\Delta t=2,5.$ The AOSCN2 method was only unstable for $\Delta t=5.$ The ADI methods, however, showed much reduced stability, with ADI1 only stable for $\Delta t=0.0005,0.001,0.002,0.005,$ and ADI2 only stable for $\Delta t=0.0005,0.001,0.002,0.005,0.01$. We note again that the stability range of Table \ref{table.1ajjstabver} is simply based on a set of discrete $\Delta t$ values being sampled. Such a range does not give a sharp estimate to the upper bound of the actual stability interval. For example, the LODCN methods are unstable at $\Delta t=1$ and stable at $\Delta t=0.5$. So only $\Delta t=0.5$ is included in the stability range in Table \ref{table.1ajjstabver}. The real critical $\Delta t$ is in between $0.5$ and $1$, and has to be detected separately. In this case, such a critical value is about $0.93$. In summary, the newly created schemes, particularly the LODIE schemes, exhibit much improved stability over the existing ADI schemes.

\begin{table}[!tb]
\caption{Electrostatic solvation free energies (kcal/mol) for protein 1ajj using the LOD methods.}
\label{table.1ajjLOD}
\begin{center}
\begin{tabular}{llllll}
\hline
Scheme & $\Delta t$& $\Delta G$ &  \% Error & Time (mm:ss.ms) \\
\hline
LODIE1 			& $0.0050$  		& $-\phantom{00}80.2$  	& $\phantom{0}93.37\%$	&$\phantom{0000}02:30.476$ \\
LODIE1 			& $0.0020$ 		& $-\phantom{0}759.8$ 	& $\phantom{0}37.19\%$	&$\phantom{0000}06:08.715$ \\
LODIE1 			& $0.0010$  		& $-1203.8$  			& $\phantom{0}18.57\%$	&$\phantom{0000}11:45.158$ \\
LODIE1 			& $0.0005$ 		& $-1097.5$ 			& $\phantom{00}9.28\%$	&$\phantom{0000}24:22.774$ \\
\hline
LODIE2 			& $0.0050$ 		& $-1521.2$ 			& $\phantom{0}25.75\%$	&$\phantom{0000}02:21.649$ \\
LODIE2 			& $0.0020$ 		& $-1336.1$ 			& $\phantom{0}10.46\%$	&$\phantom{0000}05:34.977$ \\
LODIE2 			& $0.0010$ 		& $-1273.3$ 			& $\phantom{00}5.26\%$ 	&$\phantom{0000}10:47.592$ \\
LODIE2 			& $0.0005$ 		& $-1241.6$ 			& $\phantom{00}2.64\%$	&$\phantom{0000}24:11.753$ \\
\hline
LODCN1 			& $0.0050$ 	 	& $-\phantom{0}477.2$  	& $\phantom{0}60.55\%$	&$\phantom{0000}02:45.387$ \\
LODCN1 			& $0.0020$  		& $-\phantom{0}919.3$  	& $\phantom{0}24.00\%$	&$\phantom{0000}06:22.503$ \\
LODCN1 			& $0.0010$  		& $-1065.0$  			& $\phantom{0}11.96\%$	&$\phantom{0000}12:52.169$ \\
LODCN1 			& $0.0005$ 		& $-1137.5$ 			& $\phantom{00}5.97\%$	&$\phantom{0000}25:34.671$ \\
\hline
LODCN2 			& $0.0050$ 		& $-1918.2$ 			& $\phantom{0}58.57\%$	&$\phantom{0000}02:39.077$ \\
LODCN2 			& $0.0020$ 		& $-1495.7$ 			& $\phantom{0}23.64\%$	&$\phantom{0000}06:12.191$ \\
LODCN2 			& $0.0010$ 		& $-1353.2$		 	& $\phantom{0}11.86\%$ 	&$\phantom{0000}12:18.983$ \\
LODCN2 			& $0.0005$ 		& $-1281.6$		 	& $\phantom{00}5.94\%$	&$\phantom{0000}24:36.294$ \\
\hline
\end{tabular}
\end{center}
\end{table}
\begin{table}[!tb]
\caption{Electrostatic solvation free energies (kcal/mol) for protein 1ajj using the AOS methods.}
\label{table.1ajjAOS}
\begin{center}
\begin{tabular}{llllll}
\hline
Scheme & $\Delta t$& $\Delta G$ &  \% Error & Time (mm:ss.ms) \\
\hline
AOSIE1 			& $0.005$ 		& $\phantom{-000}1.1$ 	& $100.09\%$ 			&$\phantom{0000}02:48.719$ \\
AOSIE1 			& $0.0020$ 		& $-\phantom{0}705.4$ 	& $\phantom{0}41.69\%$	&$\phantom{0000}06:38.573$ \\
AOSIE1 			& $0.0010$  		& $-\phantom{0}953.3$  	& $\phantom{0}21.19\%$	&$\phantom{0000}13:22.339$ \\
AOSIE1 			& $0.0005$ 		& $-1080.4$ 			& $\phantom{0}10.69\%$	&$\phantom{0000}28:00.736$ \\
\hline
AOSIE2 			& $0.005$ 		& $\phantom{-000}4.6$ 	& $100.38\%$			&$\phantom{0000}03:05.225$ \\
AOSIE2 			& $0.0020$ 		& $-\phantom{0}704.1$ 	& $\phantom{0}41.79\%$	&$\phantom{0000}06:49.526$ \\
AOSIE2 			& $0.0010$ 		& $-\phantom{0}952.7$ 	& $\phantom{0}21.24\%$ 	&$\phantom{0000}13:39.326$ \\
AOSIE2 			& $0.0005$ 		& $-1080.1$ 			& $\phantom{0}10.71\%$	&$\phantom{0000}28:00.741$ \\
\hline
AOSCN1 			& $0.005$ 		& $-\phantom{0}305.6$  	& $\phantom{0}74.74\%$	&$\phantom{0000}02:58.692$ \\
AOSCN1 			& $0.0020$  		& $-\phantom{0}947.1$  	& $\phantom{0}21.70\%$	&$\phantom{0000}08:00.691$ \\
AOSCN1 			& $0.0010$  		& $-1076.9$  			& $\phantom{0}10.97\%$	&$\phantom{0000}15:01.169$ \\
AOSCN1 			& $0.0005$ 		& $-1143.0$ 			& $\phantom{00}5.51\%$	&$\phantom{0000}31:22.666$ \\
\hline
AOSCN2 			& $0.005$ 		& $-\phantom{0}572.6$ 	& $\phantom{0}52.66\%$	&$\phantom{0000}03:25.016$ \\
AOSCN2 			& $0.0020$ 		& $-\phantom{0}946.5$ 	& $\phantom{0}21.75\%$	&$\phantom{0000}07:30.646$ \\
AOSCN2 			& $0.0010$ 		& $-1076.6$		 	& $\phantom{0}11.00\%$ 	&$\phantom{0000}14:32.417$ \\
AOSCN2 			& $0.0005$ 		& $-1142.8$		 	& $\phantom{00}5.52\%$	&$\phantom{0000}29:52.148$ \\
\hline
\end{tabular}
\end{center}
\end{table}
\begin{table}[!tb]
\caption{Electrostatic solvation free energies (kcal/mol) for protein 1ajj using the MAOS methods.}
\label{table.1ajjMAOS}
\begin{center}
\begin{tabular}{llllll}
\hline
Scheme & $\Delta t$& $\Delta G$ &  \% Error & Time (mm:ss.ms) \\
\hline
MAOSIE1 		& $0.005$ 		& $-\phantom{0}571.0$ 	& $\phantom{0}52.80\%$	&$\phantom{0000}03:18.824$ \\
MAOSIE1			& $0.0020$ 		& $-\phantom{0}838.4$ 	& $\phantom{0}30.69\%$	&$\phantom{0000}06:34.288$ \\
MAOSIE1			& $0.0010$  		& $-1022.0$		  	& $\phantom{0}15.51\%$	&$\phantom{0000}13:02.756$ \\
MAOSIE1			& $0.0005$ 		& $-1115.2$ 			& $\phantom{00}7.81\%$	&$\phantom{0000}25:55.982$ \\
\hline
MAOSCN1 		& $0.005$  		& $-\phantom{0}748.8$  	& $\phantom{0}38.10\%$	&$\phantom{0000}03:16.729$ \\
MAOSCN1		& $0.0020$ 		& $-1022.9$		 	& $\phantom{0}15.44\%$	&$\phantom{0000}06:57.526$ \\
MAOSCN1		& $0.0010$ 		& $-1115.7$ 			& $\phantom{00}7.77\%$ 	&$\phantom{0000}14:46.882$ \\
MAOSCN1		& $0.0005$ 		& $-1162.5$ 			& $\phantom{00}3.90\%$	&$\phantom{0000}29:11.307$ \\
\hline
\end{tabular}
\end{center}
\end{table}

\subsection{Accuracy improvements}
We next examine the accuracy by also considering the protein 1ajj. Tables \ref{table.1ajjLOD} through \ref{table.1ajjMAOS} show calculated solvation energies of various time splitting schemes based on $T=10$ and different $\Delta t$ values. For each case, the CPU time and the percentage error with respect to the reference value generated by the Euler solution are also reported. For each scheme, it can be observed that the accuracy improves when $\Delta t$ becomes smaller, while the CPU time also increases dramatically. Unfortunately, to produce a reasonable solvation energy value, a very small $\Delta t=0.0005$ is typically required in these time splitting schemes, which implies a very slow computation. The superior stability of the proposed time splitting schemes cannot be taken advantage of, due to such a poor accuracy. On the other hand, we note that the LODIE2 scheme clearly produced the smallest error, while also requiring the least CPU time for all $\Delta t$. Based on these observations, and the unconditional stability of the LODIE2 method, we will focus only on the LODIE2 scheme in the following studies. Several remedies will be proposed to improve the accuracy of the LODIE2 scheme.

Our main goal here is to produce a reliable, albeit not highly accurate, estimate of the solvation energy based on some ultra large $\Delta t$ values, so that the efficiency of the LODIE2 scheme can be significantly enhanced. Given the nearly first order convergence of the LODIE2 scheme shown in Table \ref{table.1ajjLOD}, we first investigate the use of the Richardson extrapolation technique to accelerate the convergence. With this process, we may cancel out high order error terms in our calculated electrostatic potentials by taking a linear combination of two results, calculated with different $\Delta t$ time increments. In particular, we propose the following pointwise Richardson approximation
\begin{equation}\label{rich}
\phi_m (x_i,y_j,z_k) = 2  \phi_m^{\Delta t /2} (x_i,y_j,z_k)  - \phi_m^{\Delta t } (x_i,y_j,z_k)
\end{equation}
where $\phi_m^{\Delta t } (x_i,y_j,z_k)$ represents the steady state electrostatic potential at a point $(x_i,y_j,z_k)$ calculated with a time increment $\Delta t$. Essentially, a linear combination of our vectors of potentials calculated with time increment $\Delta t$ and $\Delta t/2$ is taken.

\begin{table}[!tb]
\caption{Electrostatic solvation free energies (kcal/mol) for protein 1ajj using the improved LODIE2 and ADI methods.}
\label{table.1ajjOTHER}
\begin{center}
\begin{tabular}{llllll}
\hline
Scheme & $\Delta t$& $\Delta G$ &  \% Error & Time (mm:ss.ms) \\
\hline
LODIE2 RE 		& $0.0050, 0.0025$ 	& $-1212.8$ 	& $\phantom{00}0.25\%$	&$\phantom{0000}06:15.884$ \\
LODIE2 RE 		& $0.0150, 0.0075$  & $-1234.7$  	& $\phantom{00}2.07\%$	&$\phantom{0000}02:12.502$ \\
LODIE2 RE 		& $0.0500, 0.0250$ 	& $-1435.5$	& $\phantom{0}18.67\%$	&$\phantom{0000}00:48.154$ \\
LODIE2 RE+V	 	& $0.4000, 0.2000$ 	& $-1122.4$ 	& $\phantom{00}5.27\%$	&$\phantom{0000}00:27.030$ \\
\hline
ADI1 	& $0.005$ 		& $-1190.1$ 	& $\phantom{00}1.62\%$	&$\phantom{0000}02:44.651$ \\
ADI2 	& $0.005$  		& $-1203.8$  	& $\phantom{00}0.49\%$	&$\phantom{0000}03:22.330$ \\
\hline
\end{tabular}
\end{center}
\end{table}

The results of the LODIE2 scheme after applying the Richardson extrapolation are given in Table \ref{table.1ajjOTHER}. Such a new LODIE2 method is labeled as the LODIE2 RE method. The reported CPU time of the LODIE2 RE scheme includes calculations of $\Delta t$, $\Delta t /2$, and Eq. (\ref{rich}). Typically, such a CPU time is about 3 times larger than that of the LODIE method with the same $\Delta t$. Nevertheless, as can be seen from Table \ref{table.1ajjOTHER}, the LODIE2 RE scheme becomes about 100 times more accurate. For a comparison, the results of the ADI1 and ADI2 methods \cite{Zhao13} are also listed in Table \ref{table.1ajjOTHER}. It can be seen that by using the same $\Delta t=0.005$, the LODIE2 RE scheme is more accurate than the ADI methods, although it is also more expensive. We note that, being unconditionally stable, the LODIE2 RE method could be more efficient by using a larger $\Delta t$. On the other hand, the ADI methods will be unstable when $\Delta t > 0.01$. 

To further improve the energy estimate of the LODIE2 scheme, we propose to apply the LODIE method to calculate both $\phi_m$ and $\phi_0$ in Eq. (\ref{deltaG}). To our knowledge, such a treatment has never been explored in the existing pseudo-transient continuation approaches for solving the NPB equation \cite{Zhao11,Zhao13,Zhao14}. In the literature, the pseudo-time methods are only applied to solve $\phi_m$, while $\phi_0$ is still obtained via the FFT fast Poisson solver. Symbolically, we can rewrite Eq. (\ref{deltaG}) as
\begin{equation}\label{deltaG2}
\Delta G = \frac{1}{2} \sum\limits_i \sum\limits_j \sum\limits_k Q(x_i,y_j,z_k)(\phi^{\rm LOD}_m (x_i,y_j,z_k) - \phi^{\rm FFT}_0 (x_i,y_j,z_k)),
\end{equation}
in which $\phi^{\rm LOD}_m (x_i,y_j,z_k)$ is calculated by the LODIE2 method, while $\phi^{\rm FFT}_0 (x_i,y_j,z_k)$ is calculated by the FFT method. In this work, we propose to solve the Poisson equation in vacuum by using the LODIE2 method and denote the corresponding potential as $\phi^{\rm LOD}_0 (x_i,y_j,z_k)$. The electrostatic free energy is then calculated as 
\begin{equation}\label{deltaG3}
\Delta G = \frac{1}{2} \sum\limits_i \sum\limits_j \sum\limits_k Q(x_i,y_j,z_k)(\phi^{\rm LOD}_m (x_i,y_j,z_k) - \phi^{\rm LOD}_0 (x_i,y_j,z_k)).
\end{equation}
Even though the application of the LODIE2 scheme to the vacuum case would produce greater error in free energy calculation than the application of the FFT, we should expect the cancellation of some spatial-temporal discretization error between the vacuum and solvent cases. Thus, we should expect a potentially more accurate free energy of solvation value. 

By applying the LODIE2 RE method to both solvent and vacuum cases, the new method is labeled as LODIE2 RE+V. Furthermore, to achieve the best balance between the accuracy and efficiency, we have studied various combinations of time increment $\Delta t$ and the model scaling parameter $\alpha$. In practical testing of the LODIE2 RE+V method implementation, we find temporal convergence to be slow for the vacuum case. Thus, we shall choose the scaling parameter $\alpha$ other than one. The optimal practical scaling parameter through empirical testing was determined to be $\alpha = \frac{1}{25}$. We find that the resulting method produces superior results with $T=10$ and $\Delta t =0.4$. The results of such a LODIE2 RE+V method are also given in Table \ref{table.1ajjOTHER}. It can be seen that the relative error of the LODIE2 RE+V method with $\Delta t =0.4$ is about 5\%, which is acceptable in common biological simulations. However, the CPU time is roughly 15\% of that of the ADI methods.

\subsection{Solvation energies of proteins}
We finally validate the proposed unconditionally stable LODIE2 RE and LODIE2 RE+V methods by considering a series of proteins with different size and geometric structures. The conditionally stable ADI methods \cite{Zhao13} will also be tested for a comparison. In the following studies, the same stopping time $T=10$ is used in all methods. For the ADI, LODIE2 RE and LODIE2 RE+V methods, the time increment is chosen as $\Delta t=0.005$, $\Delta t=0.05$, and $\Delta t=0.4$, respectively. We note that the ADI results in the present study are different from those in \cite{Zhao13}, because a stronger nonlinearity is considered here with the ionic strength $I=9.48955$M, while $I=0.15$M in \cite{Zhao13}.

\begin{table}[!tb]
\caption{Electrostatic solvation free energies (kcal/mol) for 23 proteins. All proteins use a density parameter of 10 except for 451c and 1a7m, which use density parameters of 58 and 40, respectively.}
\label{table.proteins}
\begin{center}
\begin{tabular}{lllllll}
\hline
PDB ID & \# atoms & Euler	 	& ADI1 		& ADI2 		& LODIE2 RE 		& LODIE2 RE+V\\
\hline
1ajj & 519 	& $-1209.7$ 	& $-1190.1$ 	& $-1203.8$	&$-1435.5$ 		&$-1122.4$	\\
1bbl & 576 	& $-1269.8$  	& $-1246.6$	& $-1263.6$	&$-1517.1$ 		&$-1221.5$	\\
1bor & 832 	& $-1043.0$  	& $-1029.7$  	& $-1039.2$	&$-1351.2$		&$-1048.5$	\\
1bpi & 898 	& $-1440.6$ 	& $-1415.9$ 	&  $-1432.4$	&$-1799.9$		&$-1300.3$	\\
1cbn & 648 	& $-458.3$  	& $-449.9$  	&  $-455.7$	&$-717.1$			&$-489.7$		\\
1fca & 729 	& $-1109.9$ 	& $-1098.4$ 	&  $-1106.7$	&$-1407.3$		&$-1157.0$	\\
1frd & 1478	& $-2645.1$ 	& $-2611.9$ 	&  $-2634.9$ 	&$-3252.5$		&$-2530.3$	\\
1fxd & 824 	& $-2235.8$ 	& $-2220.6$ 	&  $-2231.5$	&$-2599.0$		&$-2342.6$	\\
1hpt & 858 	& $-1139.0$  	& $-1116.8$  	&  $-1132.3$	&$-1512.4$		&$-1096.1$	\\
1mbg & 903 	& $-1401.7$ 	& $-1383.2$ 	& $-1395.8$	&$-1763.6$		&$-1347.0$	\\
1neq & 1187	& $-2092.9$ 	& $-2064.6$ 	& $-2083.9$	&$-2577.5$		&$-1974.9$	\\
1ptq & 795 	& $-1098.5$  	& $-1079.2$  	& $-1092.0$	&$-1389.9$		&$-985.6$		\\
1r69 & 997 	& $-1302.1$ 	& $-1283.4$ 	& $-1296.1$	&$-1701.6$		&$-1230.6$	\\
1sh1 & 702 	& $-994.6$  	& $-981.1$  	&  $-990.8$	&$-1279.3$		&$-970.1$		\\
1svr & 1435	& $-2236.3$ 	& $-2205.7$ 	& $-2226.5$	&$-2812.2$		&$-2080.8$	\\
1uxc & 809 	& $-1363.8$ 	& $-1345.7$ 	&  $-1357.9$	&$-1671.3$		&$-1275.4$	\\
1vii & 596 		& $-1109.7$  	& $-1097.0$  	&  $-1105.7$	&$-1307.4$		&$-1046.5$	\\
2erl & 573 	& $-925.5$  	& $-961.6$	&  $-922.2$	&$-1163.1$		&$-942.7$		\\
2pde & 667 	& $-974.4$  	& $-962.5$  	& $-970.9$ 	&$-1243.0$		&$-897.2$		\\
451c & 1216	& $-1326.7$ 	& $-1845.8$ 	&  $-1341.8$	&$-1800.2$		&$-1295.0$	\\
1a2s & 1272	& $-1814.0$ 	& $-1794.1$ 	&  $-1808.1$	&$-2320.2$		&$-1898.3$	\\
1a63 & 2065	& $-3117.2$ 	& $-3072.2$ 	&  $-3102.4$	&$-3937.5$		&$-2931.3$	\\
1a7m & 2809	& $-2545.6$ 	& $-1845.8$	&  $-1847.4$	&$-2950.0$		&$-1921.9$	\\
\hline
\end{tabular}
\end{center}
\end{table}

We first consider a set of 23 proteins, which have been used for testing the previous solvation models \cite{Zhao11,Zhao13,Zhao14}. For a reference, the explicit Euler method with  $T=10$ and $\Delta t = 1.0\text{\sc{e}-}05$ is also employed. The electrostatic solvation free energies calculated by the tested methods are given in Table \ref{table.proteins}. As in the previous studies, the Euler solution can be regarded as the reference energy value, to which all tested time splitting methods should converge. It can be observed from Table \ref{table.proteins} that the LODIE2 RE method with $\Delta t =0.05$ yields a poor accuracy, even though it is much faster than the ADI methods. On the other hand, being even faster, the LODIE2 RE+V method produces energy estimates which are in good agreement with those of the Euler method. We note that the results of two proteins, i.e., 451c and 1a7m, appear to be inconsistent. In generating the molecular surface of the tested proteins, a density value of 10 is used in the MSMS software \cite{Sanner}. However, the ADI1 method with $\Delta t=0.005$ turns out to be unstable for these two proteins, while such an instability was not encountered in our previous study \cite{Zhao13}, because a strong nonlinearity with $I=9.48955$M is considered here. To fix the problem while still using $\Delta t=0.005$ for the ADI1, we have chosen the density parameters to be 58 and 40, respectively, for 451c and 1a7m in the MSMS package. Due to this change, the differences in energy results of the tested methods are quite large for these two proteins.

\begin{table}[!tb]
\caption{Electrostatic solvation free energies (kcal/mol) for additional proteins.}
\label{table.proteins2}
\begin{center}
\begin{tabular}{lllll}
\hline
PDB ID & \# atoms 	& ADI1 		& ADI2 		& LODIE2 RE+V\\
\hline
2lzj & 1687 	& --		 	& $-1983.2$	&$-1922.6$	\\
2l6n & 2093 	& --		 	& $-2593.5$	&$-2533.8$	\\
1hd2 & 2414 	& $-2203.6$ 	& $-2237.6$	&$-2121.8$	\\
1urm & 2425 	& --		 	& $-2261.0$	&$-2155.7$	\\ 
1qah & 3954 	& --		 	& $-3853.4$	&$-3513.6$	\\
3lup & 4363 	& $-6090.9$ 	& $-6147.5$	&$-5730.8$	\\
1eri & 4522 	& $-5482.5$ 	& $-5549.4$	&$-5134.9$	\\
1beb & 4972	& --		 	& $-6701.7$	&$-6208.2$	\\
2a33 & 5182 	& $-5074.6$ 	& $-5131.9$	&$-4914.5$	\\
1wkd & 5755 	& --		 	& $-5335.5$	&$-5069.4$	\\
1e5m & 6060 	& --		 	& $-4863.0$	&$-4528.6$	\\
2w2o & 6246	& --		 	& $-6481.0$	&$-5940.6$	\\ 
2w2q & 6253 	& $-6331.0$ 	& $-6399.2$	&$-5795.5$	\\
2w2p & 6263	& --		 	& $-6369.8$	&$-5742.3$	\\
2w2m & 6288 	& --		 	& $-6804.7$	&$-6042.3$	\\
4dn3 & 6344 	& $-3718.5$ 	& $-3721.7$	&$-3789.8$	\\
3r79 & 6737 	& $-8065.4$ 	& $-8146.4$	&$-7525.4$	\\
1bif & 6974 	& $-7137.7$ 	& $-7050.5$	&$-6381.2$	\\
1gsu & 7164 	& $-5906.5$ 	& $-5964.6$	&$-5531.7$	\\
4dn4 & 7390	& $-5947.1$ 	& $-5970.1$	&$-5826.7$	\\
3loq & 8136 	& $-12060.4$ 	& $-11494.7$	&$-10466.5$	\\
3gcw & 8505 	& $-8594.8$ 	& $-8568.3$	&$-7514.3$	\\ 
3bps & 8513 	& --		 	& $-8389.9$	&$-7517.9$	\\
1rva & 8726	& --		 	& $-9007.5$	&$-7892.4$	\\
1vng & 8808 	& $-7735.6$ 	& $-7798.1$	&$-6596.6$	\\
1vns & 8815	& --			& $-7908.3$	&$-6976.1$	\\ 
\hline
\end{tabular}
\end{center}
\end{table}

We next study a set of 26 proteins with larger atom numbers. The electrostatic solvation free energies produced by the ADI1, ADI2, and LODIE2 RE+V methods are shown in Table \ref{table.proteins2}. A density parameter of 10 is used in the MSMS package for all calculations. The proteins for which the ADI1 method was unstable for this parameter are marked with a dash in the ADI1 column. No reference solutions are computed for these proteins, because the computation of the Euler solution is too time consuming for large proteins. In comparing the LODIE2 RE+V results with those of ADI1 and ADI2, the present results suggest that the proposed time splitting method provides reliable and fairly accurate energy estimates.

\begin{figure}[!tb]
\centering
\begin{tabular}{cc}
\includegraphics[width=0.5\linewidth]{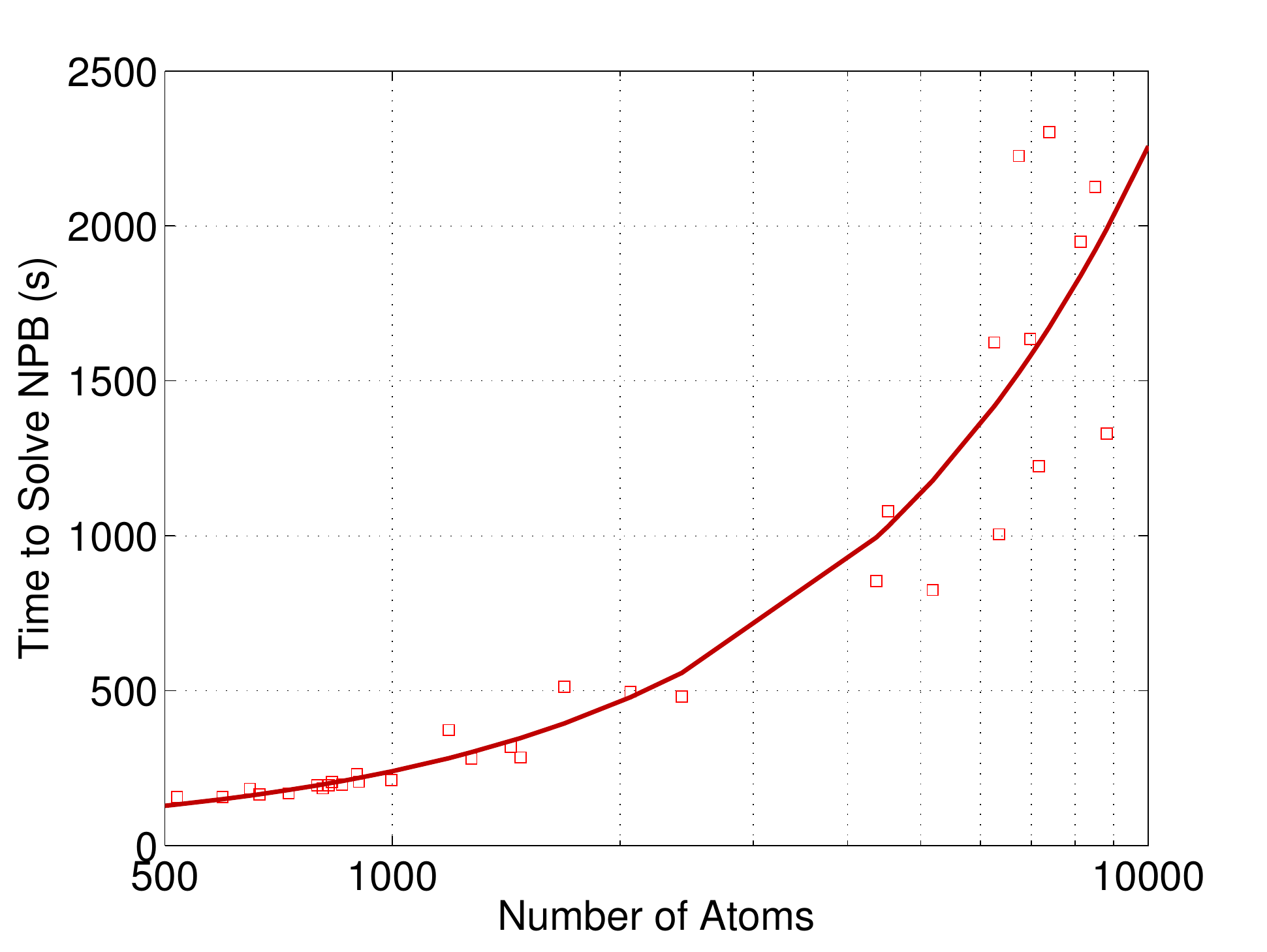} &
\includegraphics[width=0.5\linewidth]{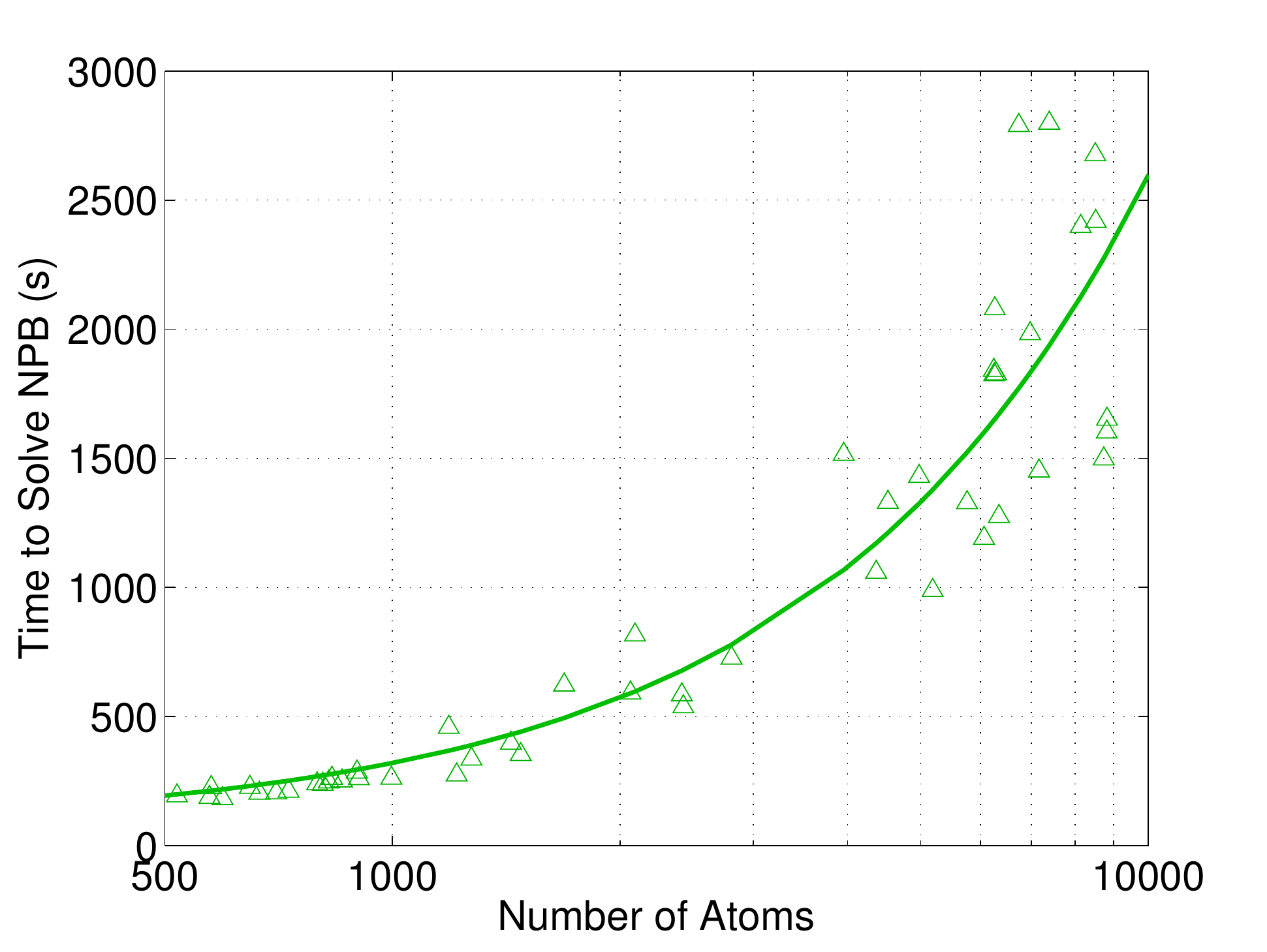} \\
(a) & (b)\\
\includegraphics[width=0.5\linewidth]{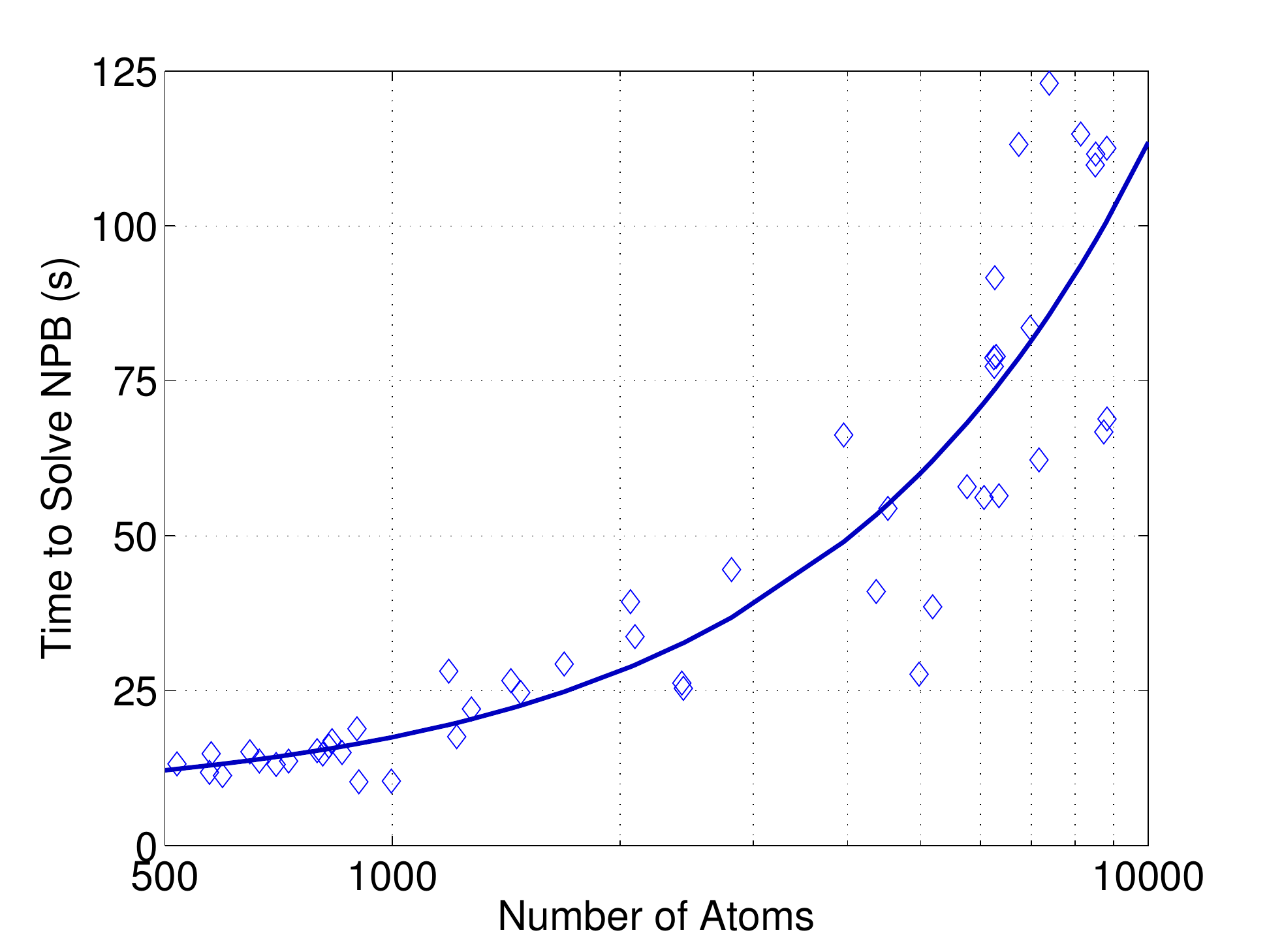} &
\includegraphics[width=0.5\linewidth]{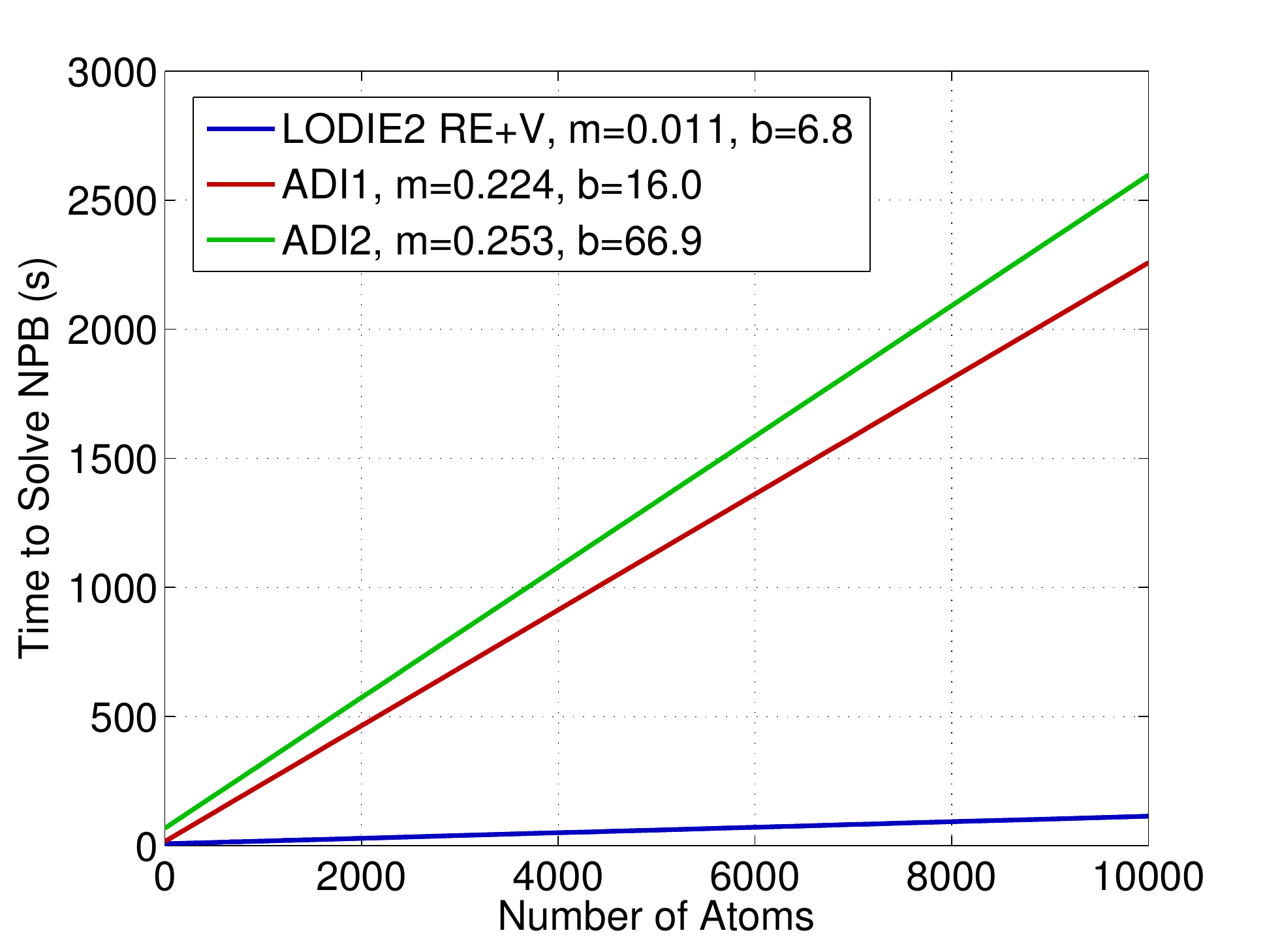} \\
(c) & (d)
\end{tabular}
\caption{The CPU time in seconds consumed by the ADI1, ADI2, and LODIE2 RE+V methods for solving the time dependent NPB equation. (a) ADI1; (b) ADI2; (c) LODIE2 RE+V; (d) A comparison of three methods. }
\label{fig.cpu}
\end{figure}

To quantitatively analyze the efficiency of the LODIE2 RE+V method, we report the CPU time in seconds of the ADI1, ADI2, and LODIE2 RE+V methods for the proteins reported in Tables \ref{table.proteins} and \ref{table.proteins2} in Fig. \ref{fig.cpu}. Note that only the CPU time consumed in solving the NPB and Poisson equations is reported here. The CPU times used for the initial numerical setup, including the trilinear interpolation of the source term, the Lagrangian to Eulerian conversion for computing the dielectric coefficient, etc., are not considered here, because they are the same in all time splitting schemes. To achieve a better understanding, we plot the CPU time against the number of atoms $N_a$ in Fig. \ref{fig.cpu}. It is clear that the CPU time is roughly a linear function of $N_a$. A least square linear fitting can be conducted to represent CPU time in seconds as a function of $N_a$: CPU$=m N_a + b$. Here, we have $(m,b)=(0.253, 66.9)$, $(0.224, 16.0)$, and $(0.011, 6.8)$, respectively, for the ADI1, ADI2, and LODIE2 RE+V methods. It is clear from Fig. \ref{fig.cpu} that the time required for the modified LODIE2 method grows much less slowly as the number of atoms increases. On average, the modified LODIE2 method is over 20 times more efficient than the ADI methods for large proteins.

\begin{figure}[!tb]
\centering
\begin{tabular}{cc}
\includegraphics[width=0.5\linewidth]{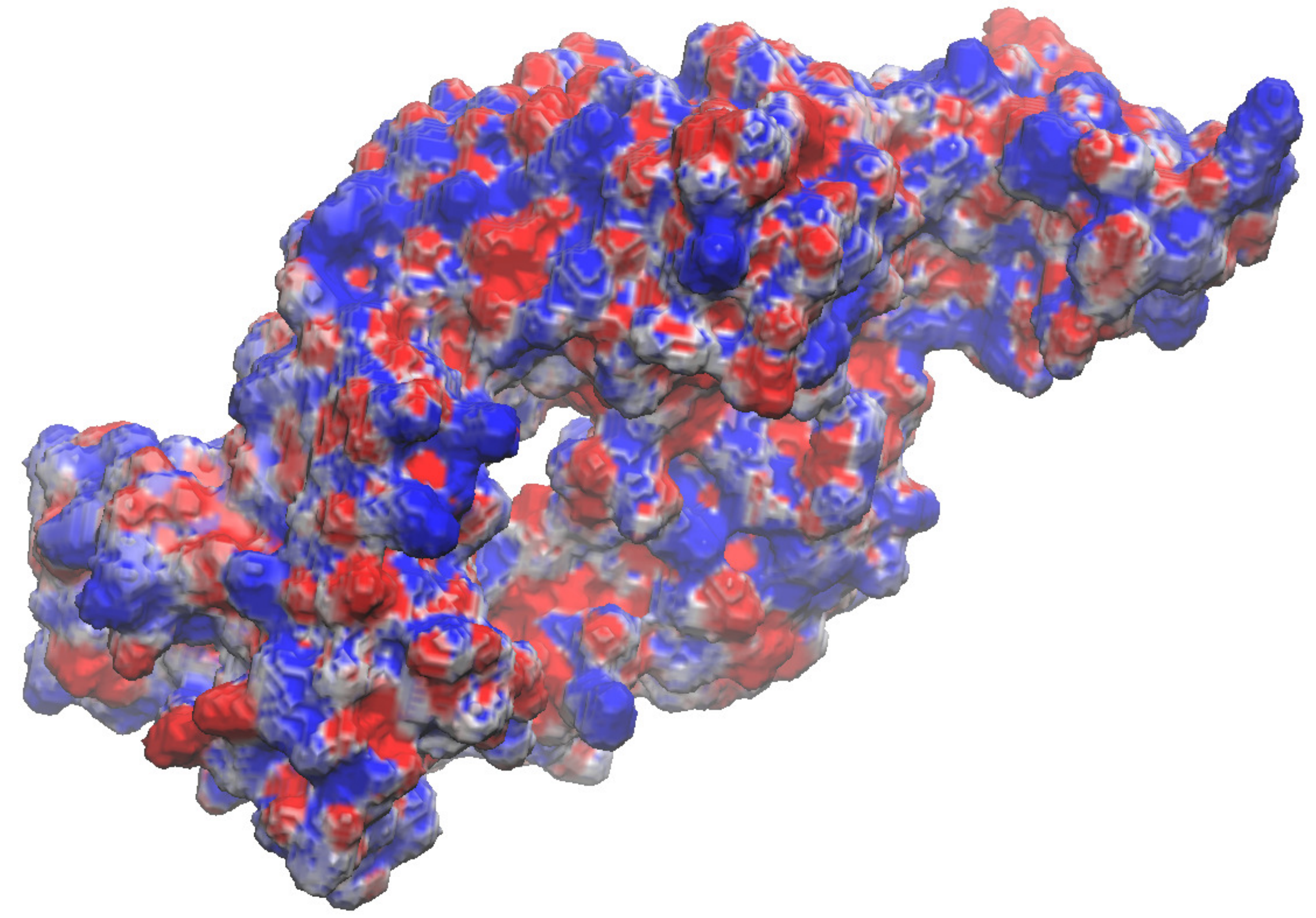} &
\includegraphics[width=0.5\linewidth]{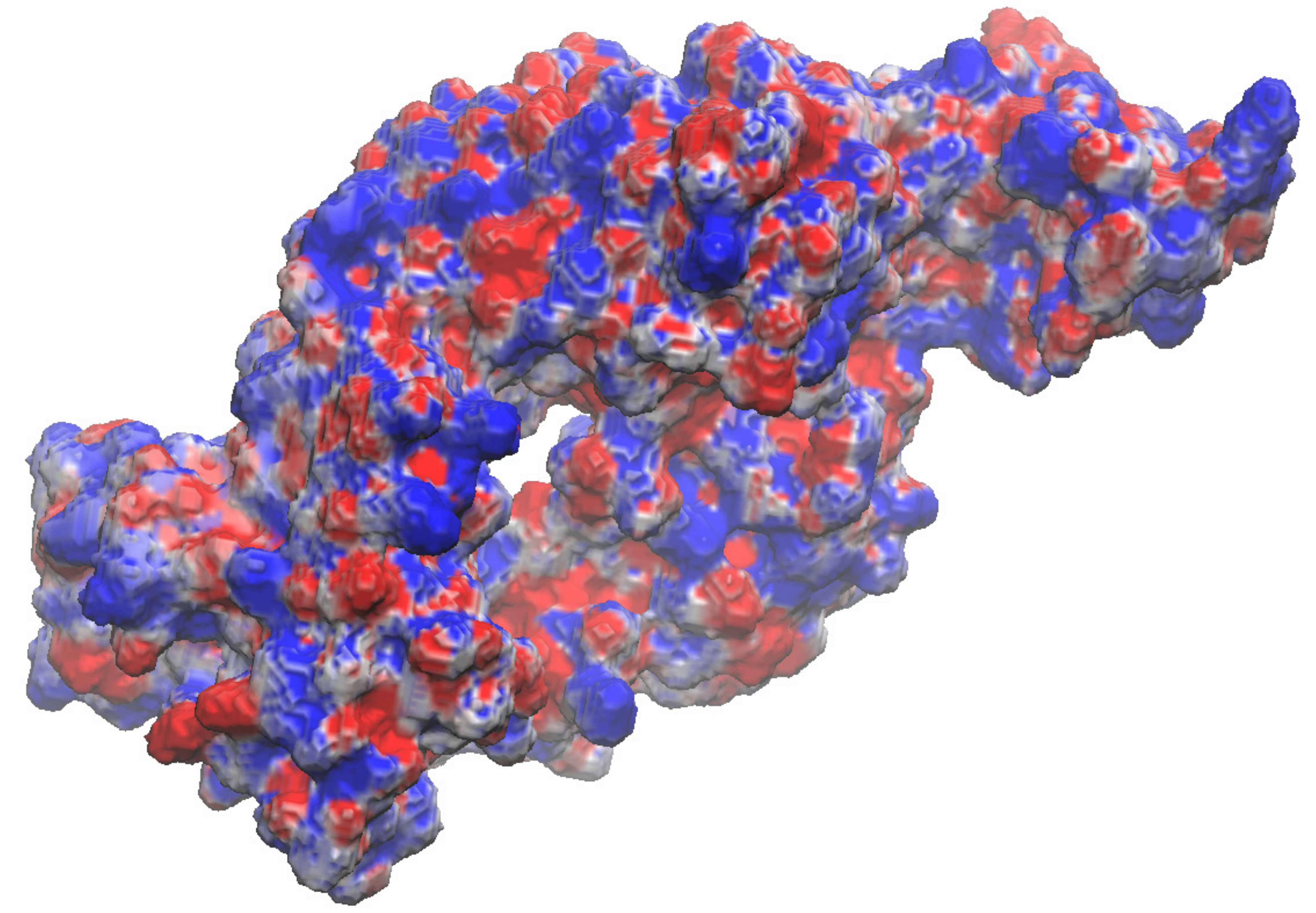}\\
(a) & (b) \\
\includegraphics[width=0.5\linewidth]{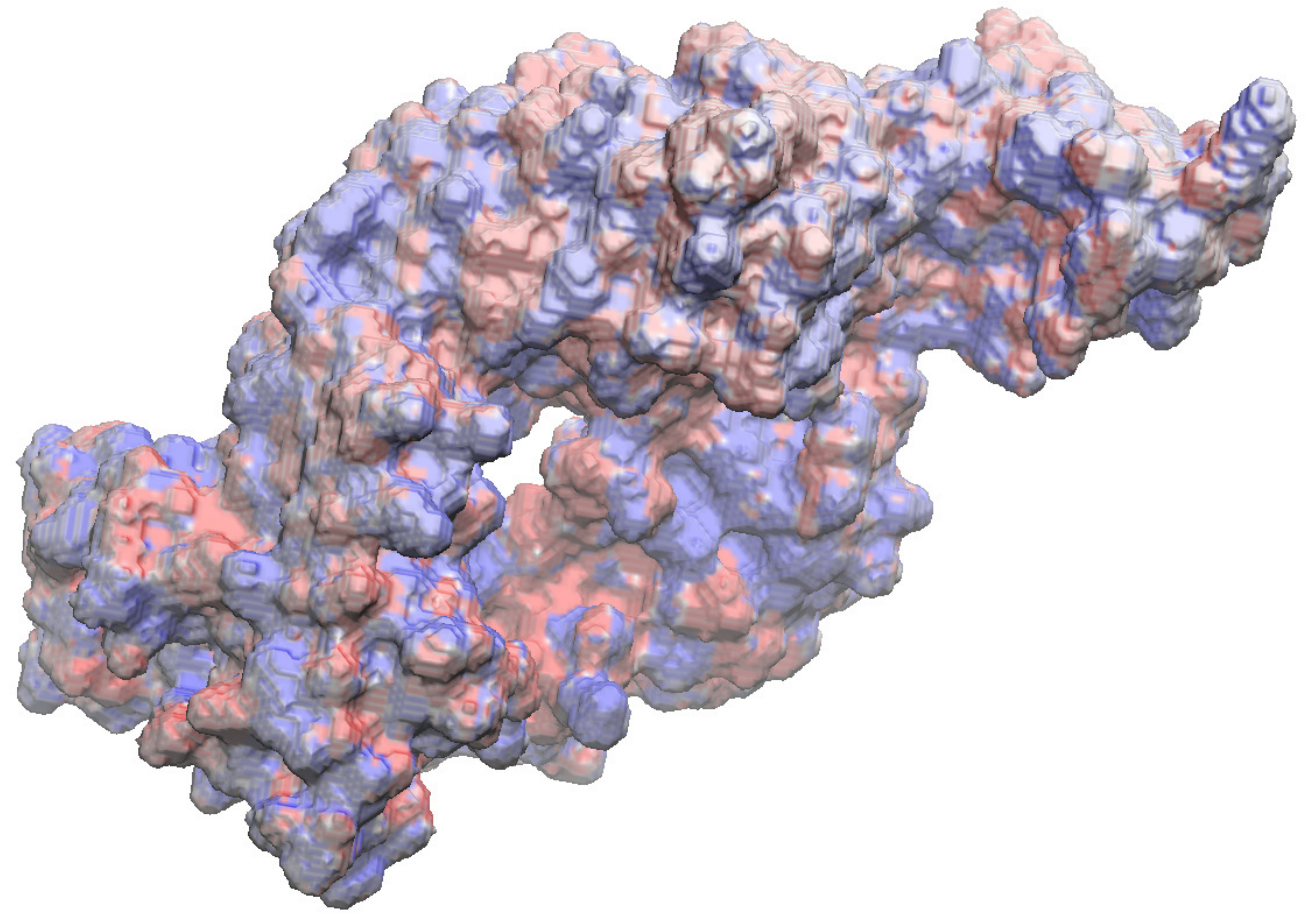} &
\includegraphics[width=0.5\linewidth]{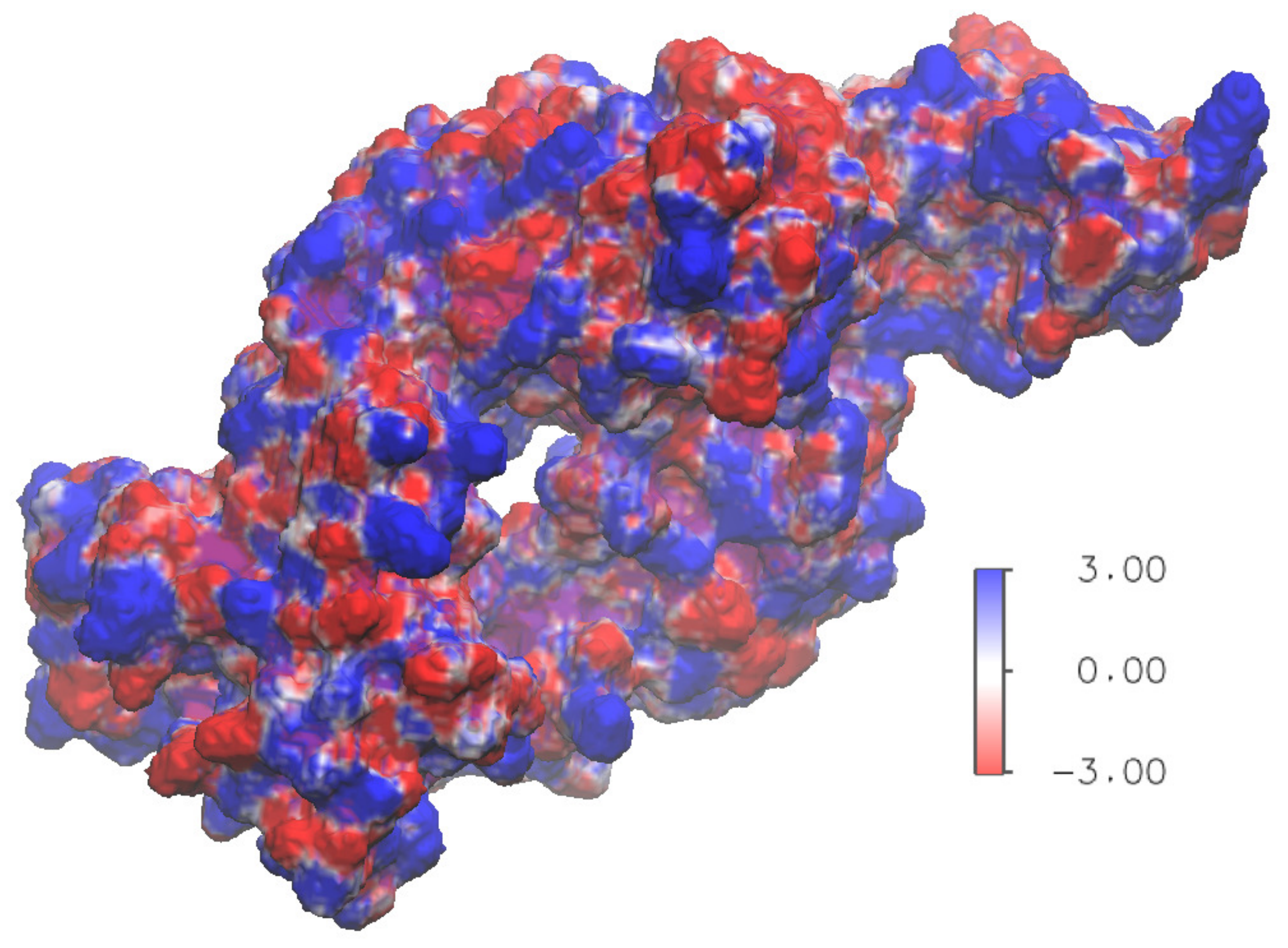}\\
(c) & (d)
\end{tabular}
\caption{Plots of the surface potential of the protein 4DN4 using different potential solutions. (a) ADI1; (b) ADI2; (c) LODIE2 RE+V; (d) LODIE2 RE+V with recovery.}
\label{fig.4DN4}
\end{figure}

Finally, we plot the surface potentials of a protein (PDB ID: 4DN4) in Fig. \ref{fig.4DN4}. Since no reference solution is offered in Table \ref{table.proteins2}, it is of interest to further study the reliability of the LODIE2 RE+V method for processing the proteins listed in that table, such as the 4DN4. The electrostatic potentials generated by the ADI1, ADI2, and LODIE2 RE+V methods are depicted in parts (a), (b), and (c), respectively. Recall that the energies of the LODIE2 RE+V method in Table \ref{table.proteins2} are after recovery. But the original potential $\phi_m^{\rm LOD}(x_i,y_j,z_k)$ is used to generate the surface plot in (c). Hence, the part (c) is significantly different from those of the ADI methods. The colors of this picture are much lighter, which stands for much weaker projected potentials. Just like in calculating the free energy, we can improve the result by considering the vacuum potentials generated by the FFT and LOD methods. In part (d), the recovered solution defined as $\phi_m^{\rm LOD}(x_i,y_j,z_k)+\phi_0^{\rm LOD}(x_i,y_j,z_k)-\phi_0^{\rm FFT}(x_i,y_j,z_k)$ is utilized. Visually, this new result is close to those of the ADI methods. Hence, the recovered LODIE2 RE+V potential can also be similarly used in analyzing the fast/slow electrostatic potential changing region on the solute-solvent boundary for biomolecular studies. The present study demonstrates that besides providing a reliable free energy estimation, the proposed LODIE method can also generate a usable pointwise potential.

\section{Conclusion}
This paper presents the first unconditionally stable numerical scheme for solving the nonlinear Poisson-Boltzmann (NPB) equation in a pseudo-transient approach. Using this approach, the solution to the NPB equation is recovered from the steady state solution to the time dependent pseudo-transient form. However, due to the long time integration of the process, it is necessary to develop numerical schemes that are stable and accurate for large time increments. The alternating direction implicit (ADI) methods previously developed for solving the time dependent NPB equation \cite{Zhao13,Zhao14} are fully implicit, but are still conditionally stable. To construct unconditionally stable solution of the unsteady NPB equation, we consider various operator splitting procedures to all five terms of the NPB equation, in both multiplicative and additive styles. This gives rise to locally one-dimensional (LOD) schemes, additive operator splitting (AOS) schemes, and multiplicative-additive operator splitting (MAOS) scheme. In these schemes, the nonlinear term is analytically integrated, and both implicit Euler and Crank-Nicholson time integrations are formulated. A standard finite difference scheme is utilized for spatial discretization in all schemes. Extensive numerical experiments are conducted to verify the unconditional stability and accuracy of the proposed time splitting schemes. One LOD scheme is found to outperform other schemes in terms of both stability and accuracy, and is recommended for electrostatic free energy analysis of real proteins. Further improvements are introduced to enhance the accuracy and efficiency of this LOD scheme for a fast biomolecular simulation. 

At last, we summarize the major numerical features of the proposed operator splitting schemes as follows:
\begin{itemize}
\item{\bf Stability.}
For smooth solution, all proposed time splitting methods should be unconditionally stable, because each individual implicit Euler or Crank-Nicolson time integration with central difference approximation is unconditionally stable. A major task in our numerical experiments is to verify the stability for non-smooth solutions. In case of a simple spherical cavity test, these methods remain unconditionally stable. For real protein systems, only the LODIE methods are stable for any $\Delta t$. Some proposed methods, such as AOSIE1, AOSCN1, and LODCN methods, become conditionally stable, while other methods, such as AOSIE2, AOSCN2, MAOS methods, are stable for all $\Delta t < 1$. Overall, the stability constraints of the proposed time splitting methods are much improved in comparing with the existing ADI methods. The LODIE method is the most stable method. 

\item{\bf Accuracy.} 
Because the charge singularities (the delta functions) and the non-smoothness of the solution across the interface $\Gamma$ are treated in an approximate sense in the present central difference discretization, all numerical schemes achieve roughly the first-order convergence in space. All proposed time splitting methods are of first order accuracy in time, because the underlying time splitting is first order and the discretization error of the implicit Euler or Crank-Nicolson integration is at least first order. However, as shown in protein tests, a very small $\Delta t$ has to be used in order to produce an accurate estimate of electrostatic free energy. A Richardson extrapolation technique and a recovery based on the replacement of the fast Fourier transform method by the operator splitting method are proposed to significantly improve the precision in energy calculations. The accelerated LODIE method becomes biologically useful, with a large $\Delta t=0.4$ and a stopping time $T=10$. 

\item{\bf Efficiency.}
By using central difference discretization, the resulting linear systems in each alternating direction are tridiagonal, and can be efficiently solved via the Thomas algorithm \cite{Mitchell}. Without the loss of generality, we can assume that the numbers of grid nodes along $x$, $y$ and $z$ directions are all on the order $N$. Then, the linear system in each subsystem of the time splitting procedure has the dimensions $N \times N$. The Thomas algorithm solution of such a system has a complexity on the order of $O(N)$, so that the operations for calculating one time step in the proposed time splitting methods are on the order of $O(N^3)$. Since a fixed number of total time steps is used in the LODIE solution, the overall complexity of the proposed LODIE method is still $O(N^3)$, making this method very attractive to large protein systems. As demonstrated in our CPU tests, the LODIE method scales linearly with respect to the number of atoms of the protein. On average, the unconditionally stable LODIE method is over 20 times more efficient than the conditionally stable ADI methods for large proteins. 

\end{itemize}

It is of great interest to further generalize the proposed time splitting methods by considering other numerical difficulties associated with the NPB equation. A regularization procedure to treat the source singularities is under our investigation. The use of a different molecular surface definition will also be examined. The stability proof of the proposed operator splitting methods remains to be an open question, mainly due to the nonlinear hyperbolic sine term.

\vspace{1cm}

\centerline{\bf Acknowledgment}
\noindent
This work was supported in part by NSF grants DMS-1016579 and DMS-1318898, and the University of Alabama Research Stimulation Program (RSP) award.

\end{document}